\documentclass[11pt]{article}
\usepackage[applemac]{inputenc}
\usepackage{geometry}
\usepackage{ulem}
\usepackage[square,authoryear]{natbib}
\usepackage[colorlinks]{hyperref}
\usepackage[T1]{fontenc}
\usepackage{latexsym,amsmath,amssymb,amsfonts, amscd,theorem, graphics}
\usepackage{array}
\usepackage{graphicx}
\usepackage{float}
\usepackage{multicol}
\usepackage{fancyhdr}
\usepackage[arrow, matrix, curve]{xy}
\usepackage{url}
\usepackage{colordvi}
\usepackage{framed}
\usepackage{color}
\usepackage{eucal}

\usepackage{mathrsfs}

\definecolor{lila}{rgb}{1,0.2,0.9}
\definecolor{brown}{rgb}{0.5,0.3,0.3}
\definecolor{turquoise}{rgb}{0.2,0.9,0.7}
\definecolor{Orange}{rgb}{0.93,0.44,0}           
\definecolor{GrayBlue}{rgb}{0.25,0.35,0.72}       
\definecolor{SeafoamGreen}{rgb}{0.54,0.71,0.50}  
\definecolor{DarkRed}{rgb}{0.7,0.15,0.3}       

\definecolor{darkorange}{cmyk}{.10,.50,.80,0}
\definecolor{lightorange}{cmyk}{.07,.37,.65,0}
\definecolor{darkpeagreen}{cmyk}{.50,.30,.50,0}
\definecolor{lightpeagreen}{cmyk}{.22,.20,.40,0}

\newcommand{\mathsym}[1]{{}}

\newtheorem{theorem}{Theorem}[section]

\newtheorem{lemma}[theorem]{Lemma}
\newtheorem{remark}[theorem]{Remark}
\newtheorem{proposition}[theorem]{Proposition}

\begin{document}
\title{Multisymplectic variational integrators and space/time symplecticity}

\author{
Fran\c{c}ois Demoures$^{1}$, Fran\c{c}ois Gay-Balmaz$^{2}$, Tudor S. Ratiu$^{3}$
}
\addtocounter{footnote}{1}
\footnotetext{Section de
Math\'ematiques and Civil Engineering, \'Ecole Polytechnique F\'ed\'erale de
Lausanne,
CH--1015 Lausanne, Switzerland. Partially supported by Swiss NSF grant 200020-137704. 
\texttt{francois.demoures@epfl.ch}
\addtocounter{footnote}{1} }
\footnotetext{Laboratoire de M\'et\'eorologie Dynamique, \'Ecole Normale Sup\'erieure/CNRS, Paris, France. Partially partially supported by a ``Projet Incitatif de Recherche'' contract from the Ecole Normale Sup\'erieure de Paris and by the government grant of the Russian Federation for support of research projects implemented by leading scientists, Lomonosov Moscow State University under the agreement No. 11.G34.31.0054. 
\texttt{francois.gay-balmaz@lmd.ens.fr}
\addtocounter{footnote}{1} }
\footnotetext{Section de
Math\'ematiques and Bernoulli Center, \'Ecole Polytechnique F\'ed\'erale de
Lausanne,
CH--1015 Lausanne, Switzerland. Partially supported by Swiss NSF grants 200020-137704,
200021-140238 and by the government grant of the Russian Federation for support of research projects implemented by leading scientists, Lomonosov Moscow State University under the agreement No. 11.G34.31.0054. 
\texttt{tudor.ratiu@epfl.ch}
\addtocounter{footnote}{1}}

\date{}

\maketitle

\begin{abstract} 
Multisymplectic variational integrators are structure preserving numerical schemes especially designed for PDEs derived from covariant spacetime Hamilton principles. The goal of this paper is to study the properties of the temporal and spatial discrete evolution maps obtained from a multisymplectic numerical scheme. Our study focuses on a 1+1 dimensional spacetime discretized by triangles, but our approach carries over naturally to more general cases.
In the case of Lie group symmetries, we explore the links between the discrete Noether theorems associated to the multisymplectic spacetime discretization and to the temporal and spatial discrete evolution maps, and emphasize the role of boundary conditions. We also consider in detail the case of multisymplectic integrators on Lie groups. Our results are illustrated with the numerical example of a geometrically exact beam model.
\end{abstract}

\tableofcontents

\section{Introduction}
\label{sec:introduction}

Multisymplectic variational integrators are structure preserving numerical schemes especially designed for solving PDEs arising from covariant Euler-Lagrange equations. These schemes are derived from a discrete version of the covariant Hamilton principle of field theory and preserve, at the discrete level, the associated multisymplectic geometry.

Multisymplectic variational integrators can be seen as the spacetime generalization of the well-known variational integrators for classical mechanics, \cite{MaWe2001}. Recall that the discrete Lagrangian flow obtained through a classical variational integrator preserves a symplectic form. From this property, it follows, by backward error analysis, that the energy is approximately preserved. 
For multisymplectic integrators, however, the situation is much more involved, the analogue of the symplectic property being given by a discrete version of the multisymplectic formula derived in \cite{MaPaSh1998}, that we will be recalled in the paper.
This formula is the spacetime analogue of the symplectic property of the discrete flow associated to variational integrators in time. The continuous multisymplectic form formula is a property of the solution of the covariant Euler-Lagrange equations in field theory, see \cite{GoIsMaMo2004} to which we refer for the multisymplectic geometry of classical field theory. In particular, several important concepts, such as covariant momentum maps associated to symmetries and the covariant Noether theorem, are naturally formulated in terms of multisymplectic forms.

In the continuous case, the articles \cite{MaSh1999}, \cite{MaPeShWe2001}, \cite{FeMaWe2003}, \cite{YaMaOr2006}, \cite{ElGa-BaHoPuRa2010} are examples of papers in which multisymplectic geometry has been further developed and applied in the context of continuum mechanics. We refer to \cite{CaRaSh2000}, \cite{CaRa2003}, \cite{ElGa-BaHoRa2011}, \cite{GB2013} for the development and the use of the techniques of reduction by symmetries for covariant field theory.

Discrete multisymplectic geometry of field theory has seen its first
development in \cite{MaPaSh1998}. The discrete Cartan forms, the 
discrete covariant momentum map, and the discrete covariant Noether 
theorems are introduced and the discrete multisymplectic form formula 
is established. This work was further developed in \cite{LeMaOrWe2003} to treat more general spacetime discretizations. This allowed the development of asynchronous variational integrators which permit the selection of independent time steps in each element, while exactly 
preserving the discrete Noether conservation and the multisymplectic structure, and offering the possibility of imposing discrete energy conservation.

\medskip

The main goal of the present paper is to study and exploit the \textit{symplectic} properties verified by the solutions of a multisymplectic variational integrator. Note, however, that the 
solution of a discrete multisymplectic scheme is a discrete spacetime section and not a discrete curve. Assuming that the spacetime is 1+1 dimensional, the discrete spacetime section can be organized either as a discrete \textit{time-evolutionary flow} or a discrete \textit{space-evolutionary flow}.
More precisely, given a discrete field that is a solution of the multisymplectic scheme, we can first construct from it a vector of  discrete positions of all spacetime nodes at a given time, and then consider the sequence of these vectors indexed by the discrete time.
Conversely, the discrete field can be organized in a space-evolutionary fashion, by first forming a vector of the discrete positions of all spacetime nodes at a given space index, and then considering the sequence of these vectors.

In order to study the symplectic character of these time-evolutionary and space-evolutionary discrete flows, we construct from the discrete covariant Lagrangian, two discrete Lagrangians ($\mathsf{L}_d$ and $\mathsf{N}_d$) associated to the temporal and spatial evolution, respectively. From this point of view, it follows naturally that a multisymplectic integrator gives rise to a variational integrator in time and a variational integrator in space. This also allows us to relate the discrete multisymplectic forms with the discrete symplectic forms associated to the time and space evolutionary descriptions. The type of boundary conditions imposed on the discrete spacetime domain are crucial
and we shall consider several type of boundary conditions. In the unrealistic situation of a spacetime without boundary, such a study 
would be essentially trivial.

Let us recall that, in the continuous setting, when the configuration field is not prescribed at the boundary, the variational principle 
yields natural boundary conditions, such as zero-traction boundary conditions. These conditions will be discretized in a structure preserving way, by means of the discrete covariant variational principle.

A main property of variational integrators (both in the multisymplectic or in the symplectic case) is that they allow for a consistent definition of discrete momentum maps and a discrete version of Noether conservation theorem in presence of Lie group symmetries. Our goal in this direction is to study and relate the discrete covariant Noether theorem associated to the discrete multisymplectic  formulation with the discrete Noether theorems associated to the time-evolutionary and space-evolutionary discrete flows built from the discrete field. Here again, this study highly depends on the type of boundary conditions involved.

Another important goal of the paper is the derivation of multisymplectic variational integrators on Lie groups. These schemes adapt to the covariant spacetime situation, such as the variational integrators on Lie groups developed in \cite{Bo-RaMa2009}, \cite{KoMa2011} which are based on the Lie group methods of \cite{IsMu-KaNoZa2000}. This approach involves the choice of a retraction map to consistently encode in
the Lie algebra the discrete displacement made on the Lie group.

The theory developed here, especially the space-evolutionary point of view, will be illustrated with the case of a geometrically exact beam. Geometrically exact models, developed in \cite{Si1985} and \cite{SiMaKr1988}, are formulated as $SE(3)$-valued covariant field theories in \cite[\S6, \S7]{ElGa-BaHoPuRa2010}. Using this covariant formulation, we will derive a multisymplectic variational integrator for geometrically exact models. As explained later, this approach further develops the variational Lie group integrators developed in \cite{DeGBLeOBRaWe2013} for geometrically exact beams.

\paragraph{Plan of the paper.} We begin by reviewing below some basic 
facts on discrete Lagrangian mechanics, following \cite{MaWe2001}. 
In Section \ref{sec_CLF}, we give a quick account of the geometry of the covariant Euler-Lagrange (CEL) equations and the multisymplectic form formula. We also consider the special case when the fields are Lie group valued and present the trivialized CEL equations.
In Section \ref{sec_MVI}, we first present the main facts about multisymplectic variational integrators on a spacetime discretized by triangles. In particular, we write the discrete multisymplectic form formula, the discrete covariant momentum maps, and the discrete covariant Noether theorem. We also derive the discrete zero traction and zero momentum boundary conditions via the discrete covariant variational principle. Then we describe systematically the symplectic properties of the time-evolutionary and space-evolutionary discrete flows built from the discrete field solution of the discrete covariant Euler-Lagrange (DCEL) equations. Several cases of boundary conditions are considered.
We also study the link between the covariant discrete Noether theorem and the discrete Noether theorems associated to the time-evolutionary and space-evolutionary discrete flows. We will see that, while the covariant discrete Noether theorem holds independently of the imposed boundary conditions, this is not the case for the discrete Noether theorems associated to the time or space discrete evolutions.
Section \ref{mult_var_int} is devoted to the particular situation when the configuration field takes values in a Lie group. In this case, it is possible to trivialize the DCEL. This is a serious advantage in the discrete setting, since it allows us to make use of the vector space structure of the Lie algebra via the use of a time difference map.
Finally, in Section \ref{sect_num_ex}, we illustrate the properties of the multisymplectic variational integrator with the example of the geometrically exact beam model. The symplectic property of the space-evolutionary discrete flow is exploited to reconstruct the trajectory of the beam, knowing the evolution of the position and of the strain of one of its extremities.

\paragraph{Review of discrete Lagrangian dynamics.}
Let $Q$ be the configuration manifold of a mechanical system. 
Suppose that the dynamics of this system is described by the 
Euler-Lagrange (EL) equations associated to a Lagrangian 
$\mathsf{L}:TQ \rightarrow \mathbb{R}$ defined on the tangent 
bundle $TQ$ of the configuration manifold $Q$. Recall that 
these equations characterize the critical curves of the 
action functional associated to $\mathsf{L}$, namely
\[
\frac{d}{dt} \frac{\partial \mathsf{L}}{\partial \dot q}- 
\frac{\partial \mathsf{L}}{\partial q}=0\quad \Longleftrightarrow 
\quad \delta \int_0^T \mathsf{L}(q(t),\dot q(t))dt=0,
\]
for variations $\delta q(t) \in T_{q(t)}Q$ of the curve $q(t)$ 
vanishing at the endpoints, i.e., $\delta q(0)=\delta q(T)=0$.
The Legendre transform associated to $\mathsf{L}$ is the 
locally trivial fiber (\textit{not} vector) bundle morphism
$\mathbb{F} \mathsf{L} :TQ \rightarrow T^*Q$ covering the identity 
that associates to a velocity its corresponding conjugate momentum, 
where $T^*Q$ denotes the cotangent bundle of $Q$. In canonical
tangent and cotangent bundle charts induced by an atlas on $Q$,
it has the expression $\left(q, \dot{q} \right)\mapsto 
\left(q,\frac{\partial \mathsf{L}}{\partial \dot{q}}\right)$.

\medskip

We recall the discrete version of this approach, following 
\cite{MaWe2001}. Suppose that a time step $h$ has been fixed, 
denote by $\{t^j =jh\mid j=0,...,N\}$ the sequence of times 
discretizing $[0,T]$, and by $q_d:\{t^j\}_{j=0}^N \rightarrow Q$, 
$q^j:=q_d(t^j)$ the corresponding discrete curve. Let 
$\mathsf{L}_d: Q \times Q \rightarrow \mathbb{R}$, 
$\mathsf{L}_d=\mathsf{L}_d(q ^j , q ^{j+1})$ be a discrete 
Lagrangian which we think of as approximating the action 
integral of $\mathsf{L}$ along the curve segment between 
$q^j$ and $q^{j+1}$, that is, we have
\[
\mathsf{L}_d( q ^j , q^{j+1})\approx
\int_{t ^j }^{t^{j+1}}\mathsf{L}( q(t), \dot q(t))dt,
\]
where $q(t^j )= q^j $ and $q(t^{j+1})= q^{j+1}$.
The \textit{discrete Euler-Lagrange (DEL) equations} are obtained 
by applying the discrete Hamilton principle to the discrete 
action
\[
\mathfrak{S}_d(q_d):=\sum_{j=0}^{N-1}\mathsf{L}_d(q^j , q^{j+1})
\]
for variations $ \delta q ^j $ vanishing at the endpoints. We 
have the formula
\begin{align*} 
\delta \mathfrak{S}_d(q_d)&=\sum_{j=1}^{N-1}
\left( D_1 \mathsf{L}_d(q^j, q^{j+1})  + 
D_2 \mathsf{L}_d (q^{j-1},q^{j}) \right) \delta q ^j \\
& \qquad \qquad \qquad +\Theta ^+_{\mathsf{L}_d}( q ^{N-1}, q ^N) ( \delta q ^{N-1}, \delta q ^N ) - 
\Theta ^-_{\mathsf{L}_d}( q ^0, q ^1 ) (\delta q^0, \delta q^1),
\end{align*} 
where
\[
\Theta^-_{\mathsf{L}_d}(q^0, q^1):= - D_1\mathsf{L}_d(q^0, q^1)
\mathbf{d} q^0 \quad\text{and}\quad \Theta^+_{\mathsf{L}_d}
(q^0, q^1):=D_2 \mathsf{L}_d(q^0, q^1)\mathbf{d} q^1 
\]
are the discrete Lagrangian one-forms and $D_1$, $D_2$ denote the
first and second partial derivatives of a function on the manifold 
$Q \times Q$. The DEL equations are thus given by
\[
D_1 \mathsf{L}_d\left(q^j, q^{j+1} \right)  + D_2 \mathsf{L}_d \left(q^{j-1},q^{j} \right) = 0, \quad \text{for all} 
\quad j=1, \ldots, N-1.
\]
The \textit{discrete Legendre transforms} associated to 
$\mathsf{L}_d$ are the two maps $\mathbb{F}^+ \mathsf{L}_d, 
\mathbb{F}^-\mathsf{L}_d : Q \times Q \rightarrow T^*Q$ defined 
by
\begin{equation}
\label{discret_Legendre}
\begin{aligned} 
\mathbb{F}^+\mathsf{L}_d (q^j, q^{j+1}) &:= 
D_2\mathsf{L}_d(q^j, q^{j+1}) 
\in T_{q^{j+1}} ^*Q \\
\mathbb{F}^-\mathsf{L}_d (q^j, q^{j+1})&:= 
-D_1\mathsf{L}_d(q^j, q^{j+1}) 
\in T_{q^{j}}^*Q.
\end{aligned}
\end{equation} 
Note that the DEL equations can be 
equivalently written as
\begin{equation}\label{Leg_form}
\mathbb{F}^+\mathsf{L}_d (q^{j-1}, q^{j})=
\mathbb{F}^-\mathsf{L}_d (q^j, q^{j+1}), \quad 
\text{for $j=1,...,N-1$}
\end{equation} 
and that we have $ \Theta ^\pm_{\mathsf{L} _d }= 
(\mathbb{F}^\pm\mathsf{L}_d)^\ast \Theta_{can}$, where 
$\Theta_{can}$ is the canonical one-form on $T^*Q$, defined
by $\Theta_{can}(\alpha_q)\left(U_{\alpha_q}\right): = 
\alpha_q\left(T_{\alpha_q}\rho\left(U_{\alpha_q}\right)\right)$,
where $\alpha_q \in T^*_qQ$, $U_{\alpha_q}\in T_{\alpha_q}(T^*Q)$,
and $\rho: T^*Q \rightarrow Q$ is the cotangent bundle projection.

\paragraph{Approximate energy conservation.} The main feature of 
the numerical scheme $(q^{j-1}, q^j) \mapsto (q^j , q^{j+1})$, 
obtained by solving the DEL equations, is 
that the associated scheme $(q^j , p^j) \mapsto (q ^{j+1},p^{j+1})$ 
induced on the phase space $T^*Q$ through the discrete 
Legendre transform, defines a \textit{symplectic integrator}. 
Here we assumed that the discrete Lagrangian $\mathsf{L}_d$ 
is regular, that is, both discrete Legendre transforms 
$\mathbb{F}^+ \mathsf{L}_d, \mathbb{F}^-\mathsf{L}_d : 
Q \times Q \rightarrow T^*Q$ are local diffeomorphisms 
(for nearby $q^j$ and $q^{j+1}$). The symplectic character 
of the integrator is obtained by showing that the scheme 
$(q^{j-1}, q^j) \mapsto (q^j , q^{j+1})$ preserves the discrete symplectic  two-forms $\Omega_{\mathsf{L}_d}^\pm:= 
(\mathbb{F}^{\pm} \mathsf{L}_d)^\ast \Omega_{can}$, 
so that $(q^j , p^j ) \mapsto (q^{j+1},p^{j+1})$ preserves 
$ \Omega_{can}$ and is therefore symplectic; see \cite{MaWe2001}, \cite{LeMaOrWe2004}. Here  
$\Omega_{can}: =- \mathbf{d} \Theta_{can}$ is the canonical 
symplectic two-form on $T^*Q$; in standard cotangent bundle
coordinates it has the expression $\Omega_{can} = \mathbf{d}q \wedge \mathbf{d}p$.

It is known (see \cite{HaLuWa2006}), that given a Hamiltonian $H$, a symplectic integrator for $H$ corresponds to solving a modified Hamiltonian system for a Hamiltonian $\bar H$ which is close to $H$. 
Hence, the discrete trajectory has all of the properties of a conservative Hamiltonian system, such as energy conservation. The same conclusion holds  on the Lagrangian side for variational integrators 
(see, e.g., \cite{LeMaOrWe2004}). This explains why energy is approximately conserved for variational integrators and typically oscillates about the true energy value. We refer to  \cite{HaLuWa2006} 
for a detailed account and a full treatment of backward error 
analysis for symplectic integrators.

\color{black}

\paragraph{Lagrange-d'Alembert principle.} In the presence of an 
external force field, given by a (in general, nonlinear) fiber 
preserving map 
$F_{\mathsf{L}}: TQ \rightarrow T^*Q$, Hamilton's principle is replaced by the \textit{Lagrange-d'Alembert principle}
\begin{equation}\label{LdA_principle} 
\delta \int_0^T \mathsf{L}(q(t),\dot q(t))dt+ \int_0 ^T F_{\mathsf{L}}(q(t),\dot q(t))\cdot \delta q\,dt=0, 
\end{equation} 
where $F_{\mathsf{L}}(q, \dot q) \cdot \delta q$ is the virtual work done by the force field $F_{\mathsf{L}}$ with a virtual displacement  $ \delta q$. This principle yields the \textit{Lagrange-d'Alembert equations}
\[
\frac{d}{dt} \frac{\partial \mathsf{L}}{\partial \dot q}- \frac{\partial \mathsf{L}}{\partial q}= F_{\mathsf{L}}(q, \dot q).
\]

In the discrete case, given such forces, the discrete Hamilton principle has to be modified to the 
\textit{discrete Lagrange-d'Alembert principle}, which seeks discrete curves $\{ q ^j \}_{j=0}^N$ that satisfy
\begin{equation}\label{DLdA} 
\delta \sum_{j=0}^{N-1}\mathsf{L}_d( q ^j , q^{j+1})+ \sum_{j=0}^{N-1}\left( F_d ^{-}( q ^j , q^{j+1}) \cdot \delta q ^j + F_d ^{+}( q ^j , q^{j+1}) \cdot \delta q ^{j+1} \right) =0,
\end{equation} 
for variations $ \delta q ^j $ vanishing at endpoints, where the two discrete Lagrangian forces $F^{\pm}_d: Q \times Q \rightarrow T^*Q$ are fiber preserving maps such that the second term above is an approximation of the integral $\int_0^T F(q, \dot q) \cdot \delta q \, dt$ of the virtual work. One gets the \textit{DEL equations with forces}
\[
D_1 \mathsf{L}_d\left(q^j, q^{j+1} \right)  + D_2 \mathsf{L}_d \left(q^{j-1},q^{j} \right) +F_d^+(q^{j-1},q ^j )+ F^-( q ^j , q^{j+1})= 0, \quad j=1, \ldots, N-1.
\]
In the forced case, the discrete Legendre transforms \eqref{discret_Legendre} have to be modified to 
\[
\begin{aligned} 
\mathbb{F}^{F+}\mathsf{L}_d (q^j, q^{j+1}) &:= D_2\mathsf{L}_d(q^j, q^{j+1})+F^+( q^j, q^{j+1})
\in T_{q^{j+1}} ^*Q \\
\mathbb{F}^{F-}\mathsf{L}_d (q^j, q^{j+1})&:= -D_1\mathsf{L}_d(q^j, q^{j+1}) - F^-(q^j, q^{j+1})
\in T_{q^{j}} ^*Q.
\end{aligned}
\]
As in \eqref{Leg_form}, the DEL equations with forces can be equivalently written as
\[
\mathbb{F}^{F+}\mathsf{L}_d (q^{j-1}, q^{j})=\mathbb{F}^{F-}\mathsf{L}_d (q^j, q^{j+1}), \quad \text{for $j=1,...,N-1$}.
\]

\section{Covariant Lagrangian formulation}\label{sec_CLF}

In this section, we recall basic facts about the geometry of covariant field theory such as the covariant Euler-Lagrange (CEL) equations, the Cartan forms, covariant momentum maps, and the multisymplectic form formula, following \cite{GoIsMaMo2004} and \cite{MaPaSh1998}. We also consider the case when the fields are Lie group valued and present the trivialized CEL equations.
Although the corresponding discrete multisymplectic integrators will be considered only on trivial fiber bundles $Y:= X \times M \rightarrow X$, we often start with the general theory written on arbitrary fiber bundles, because it yields the correct guide to write geometrically consistent formulas, both at the continuous and discrete levels.

\subsection{Preliminaries on covariant Lagrangian formulation}

In classical Lagrangian mechanics, the dynamic evolution of a system is described by a curve $q(t) \in Q$ in the configuration space $Q$ of 
the system. This curve is a solution of the EL equations obtained by Hamilton's principle:
\[
\frac{d}{dt} \frac{\partial \mathsf{L}}{\partial \dot q}- 
\frac{\partial \mathsf{L}}{\partial q}=0  \quad \Longleftrightarrow 
\quad \delta \int_{0 }^{T}\mathsf{L}(q (t) , \dot q (t) ) dt=0.
\]
In continuum mechanics, the configuration space is usually a space of maps $ \varphi : B \rightarrow M$ (such as embeddings) defined on the base manifold $B$ or parametrization space, with values in the space $M$ of allowed deformations. In this situation, $Q$ is therefore an infinite dimensional space. For example, in the case of the geometrically exact beam that we will treat in \S\ref{sect_num_ex}, we have $ \varphi :B=[0,L] \rightarrow M= SE(3)$, where
$SE(3)$ is the special Euclidean group consisting of orientation
preserving rotations and translations; thus, 
$Q= \mathcal{F} ([0,L], SE(3))$ is the space of all such maps.
In many situations, the Lagrangian $\mathsf{L}:TQ \rightarrow 
\mathbb{R}$ of the system can be written in terms of a Lagrangian 
density $\mathcal{L}$ such that
\begin{equation}\label{covariant_point_of_view} 
\mathsf{L}(\varphi  , \dot \varphi )= 
\int_B \mathcal{L} (s, \varphi (s), \dot \varphi (s), 
\nabla \varphi (s)) d^n s,
\end{equation}
where $\dot{\varphi} \in T_ \varphi Q$ and $\nabla\varphi(s):T_sB \rightarrow T_{\varphi(s)}M$ denotes the derivative (tangent map) of $\varphi$ relative to the variable $s$. If the Lagrangian $\mathsf{L}$ is defined in terms of a Lagrangian
density $\mathcal{L}$ as in \eqref{covariant_point_of_view},
one can alternatively formulate the dynamics and all its properties in terms of the Lagrangian density $\mathcal{L}$ instead of the 
Lagrangian $\mathsf{L}$. This is the covariant, or field theoretic, 
point of view. In this description, the configuration $\varphi$ 
is seen as a spacetime dependent map $\mathbb{R} \times B \ni (t,s) 
\mapsto \varphi(t,s) \in M$, rather than a curve $[0,T]\ni t\mapsto \varphi (t) \in Q: =\mathcal{F} ([0,L], M)$.

Abstractly, the maps $ \varphi $ have to be interpreted as sections of the trivial fiber bundle $\pi :Y:=X \times M \rightarrow X$, $ \pi (t,s,m)=(t,s)$, where $X:= \mathbb{R}  \times B$ is the base and $M$ is the fiber. The Lagrangian density is defined on the first jet bundle $J^1Y=J ^1 (X \times M)$ of the fiber bundle $Y=X \times M$ and takes values in the space $ \Lambda ^{n+1}X$ of $(n+1)$-forms on $X= \mathbb{R}  \times B$, where $n= \operatorname{dim}B$ . 
The fiber of the first jet bundle at $y \in Y_x:= \{x\} \times M$ is 
\begin{equation}\label{fiber_jet_bundle}
J^1_yY= \{ \gamma _y \in L(T_xX, T_yY)\mid T_y \pi \circ \gamma _y = id_{T_xX}\},
\end{equation}
where $T_y \pi :T_yY \rightarrow T_xX$ is the tangent map at $y$ of 
the bundle projection $ \pi : Y \rightarrow X$. We note that $\mathrm{dim}(J^1Y) = (n+1)+m+(n+1)m$, with $\mathrm{dim}M=m$. The first jet extension of a section $ \varphi $ is $j^1_x \varphi :=T_x \varphi \in J^1_{ \varphi (x)}Y$, so that the action functional associated to $ \mathcal{L} $ can be simply written as 
$\int_X \mathcal{L}(j^1\varphi)$.

Note that a Lagrangian density defined on $J ^1 Y$ may depend explicitly on time. In \eqref{covariant_point_of_view}, however, we assumed that there is no such dependence, since this is the case in most examples in continuum mechanics.
An explicit time dependence in $ \mathcal{L} $ in \eqref{covariant_point_of_view} would induce an explicit time dependence in the Lagrangian $\mathsf{L}$. 
 
Since the fiber bundle $\pi:X \times M \rightarrow X$ is trivial, the 
first jet bundle can be identified with the bundle $L(TX, TM) \rightarrow X \times M$ over $X \times M$, whose fiber at $(x,m) \in X \times M$ is given by linear maps $L(T_xX, T_mM)$. The first jet extension 
$ j ^1 \varphi (x)$, $x=(t,s)$ is identified with 
the linear map $(t,s, \varphi (t,s),\dot \varphi (t,s), 
\nabla \varphi (t,s)): T_{(t,s)}(\mathbb{R}\times B) \rightarrow 
T_{\varphi(t,s)}M$ in the following manner: if $a\in\mathbb{R}$,  
$u_{s} \in T_{s}B$, so that $(a, u_s) \in 
T_{(t,s)}X = T_{(t,s)}(\mathbb{R} \times B)$, then
\begin{align*}
(t,s, \varphi (t,s),\dot \varphi (t,s), 
\nabla \varphi (t,s))\cdot  (a, u_s)  : =
a \dot{\varphi}(t,s) + \nabla\varphi(t,s)(u_s) \in T_{\varphi (t,s)}M.
\end{align*} 

\noindent\textsf{Notation.} We denote by $s^i $ the coordinates 
on $B$, by $y^A $ the coordinates on $M$. We use the notation 
$d^{n+1}x:= dt \wedge d^n s=dt \wedge ds^1 \wedge ... \wedge ds^n$. 
We write locally the Lagrangian density as 
$\mathcal{L} (t, s^i, v_0^A, v_i^A) = L (t, s^i, v_0^A, v_i^A) 
dt \wedge d^n s$.

\paragraph{Covariant Euler-Lagrange equations.} The Hamilton principle reads
\[
\delta \int_{0}^{T}\!\!\! \int_B \mathcal{L} (j^1 \varphi (t,s))= 
\delta \int_{0}^{T}\!\!\! \int_B L (t,s,\varphi (t,s), 
\dot{\varphi} (t,s), \nabla \varphi (t,s))dt\wedge d^n s =0,
\]
for variations $\delta \varphi$ with 
$\delta \varphi |_{ \partial ([0,T] \times B)}=0$. This principle yields the CEL equations, locally given by
\begin{equation}\label{cov_EL} 
\frac{\partial }{\partial t}\frac{\partial L }{\partial \dot \varphi ^A }+\partial _i \frac{\partial L}{\partial  \varphi_{,i} ^A }- \frac{\partial L }{\partial \varphi ^A }=0. 
\end{equation}
For completeness, we present the derivation of these
equations for $B\subset \mathbb{R}  ^n $ an open subset with compact closure and smooth boundary. We have
\begin{equation}\label{computation_EL}
\begin{aligned}
\delta \int_{0}^{T}\!\!\! \int_B \mathcal{L} (j^1 \varphi (t,s))
&= \int_0^T\!\!\! \int_B \left(\frac{\partial L}{\partial \varphi} 
\cdot \delta \varphi + \frac{\partial L}{\partial \dot{\varphi}} 
\cdot \delta \dot{\varphi} + \frac{\partial L}{\partial \nabla \varphi} 
\cdot \delta  \nabla \varphi  \right) dt \wedge d^n s \\
& = \int_0^T\!\!\! \int_B \left( \left[ 
\frac{\partial L }{\partial \varphi} - 
\frac{\partial}{\partial t}\frac{\partial L}{\partial \dot \varphi} - \operatorname{div}\frac{\partial L}{\partial \nabla\varphi} \right] \delta \varphi \right) dt \wedge d^n s \\
& \quad + \int_B \left[\frac{\partial L}{\partial \dot{\varphi}} 
\cdot \delta \varphi \right]_{0}^T d^n s  +
\int_0^T \!\!\!\int_ {\partial B}  
\frac{\partial L}{\partial \varphi_{,i}^A}\mathbf{n} ^i 
\delta \varphi^Ad^{n-1}a\wedge dt,
\end{aligned}
\end{equation} 
where $ \mathbf{n} $ is the outward pointing unit normal to the boundary $ \partial B$ and $d^{n-1}a$ is the volume form induced on $ \partial B$.
Since $ \delta \varphi |_{ \partial ([0,T] \times B)}=0$ and $\partial ([0,T] \times B)=([0,T] \times \partial B) \cup (\{0\} \times B) \cup( \{T\} \times B)$, the boundary terms vanish, thus yielding the CEL equations.

\begin{remark}[Boundary conditions]{\rm In the above situation it is assumed that the configuration $ \varphi $ is known at $t=0, T$ and is prescribed at the boundary for all times, which corresponds to \textit{pure displacement boundary conditions}. If the configuration at the boundary is not prescribed, then Hamilton's principle yields the boundary condition
\begin{align}\label{BC_traction} 
\frac{\partial \mathcal{L}}{\partial \varphi_{,i}^A  }\mathbf{n} ^i=0, \quad \text{for all $A=1,...,m$},
\end{align}
known as \textit{zero traction boundary condition}. Note that the treatment of nonzero traction $ \frac{\partial \mathcal{L}}{\partial \varphi_{,i}^A  }\mathbf{n} ^i=\boldsymbol{\tau }_A $ requires the addition of a term in the Lagrangian; see, e.g., \cite{MaHu1983}.

Other conditions could be used in the variational principle, such as the assumption of pure displacement boundary conditions but without the assumption that the configuration is known at $ t=0, T$. In this case, 
the variational principle would yield the conditions
\begin{equation}\label{BC_wird} 
\frac{\partial \mathcal{L} }{\partial  \dot \varphi ^A }(0,s)=0= \frac{\partial \mathcal{L} }{\partial  \dot \varphi ^A }(T,s), \;\; \text{for all $s \in B$}. 
\end{equation}
known as \textit{zero momentum boundary conditions}.}
\end{remark}

\paragraph{Covariant Euler-Lagrange operator.} For future use, we recall here an intrinsic way of writing the CEL equations. Let $VY$ be the vertical vector subbundle of $TY$ whose fibers are defined by 
\begin{equation}\label{vertical_subbundle}
V_{ y} Y := \{v_y \in T_{y}Y \mid T_y\pi(v _y ) =0 \}.
\end{equation}
Let $V ^\ast Y$ be its dual vector bundle. In the case of a trivial bundle $Y=X \times M$, the fiber $V_yY$, $y=(x,m)$, is identified with the tangent space $T_mM$. 

There is a unique bundle morphism $\mathscr{EL}(\mathcal{L}) :  
J^1Y \rightarrow V^*Y \otimes \Lambda ^{n+1}X$ covering the 
identity on $Y$, called the \textit{covariant Euler-Lagrange 
operator}, such that
\begin{equation}\label{intrinsic_form} 
\left.\frac{d}{d\varepsilon}\right|_{\varepsilon=0} \int_ X 
\mathcal{L} (j ^1 \varphi _\varepsilon (x))= 
\int_X \mathscr{EL}(\mathcal{L} ) \left(j^1\varphi(x) \right)\cdot \delta  \varphi(x),
\end{equation} 
for all variations $ \varphi _\varepsilon $ of $ \varphi $, among sections of $\pi:X \times M \rightarrow X$ satisfying $ \delta \varphi |_{ \partial X}=0$, where $ \delta \varphi := \left.\frac{d}{d\varepsilon}\right|_{\varepsilon=0} \varphi _\varepsilon $.

In the examples treated in this paper, the bundle $\pi:Y \rightarrow X$
is trivial and we have $X=[0,T] \times B\ni x=(t,s)$ so that the CEL operator recovers locally the expression 
of the CEL equations \eqref{cov_EL}.

\paragraph{Forced covariant Euler-Lagrange equations.} Recall that 
in the presence of  a Lagrangian force field $F_{\mathsf{L}} : 
TQ \rightarrow T^*Q$ (a fiber preserving map covering the identity,
not necessarily linear on the fibers), Hamilton's principle has to 
be replaced by the Lagrange-d'Alembert principle \eqref{LdA_principle}.
Analogously to \eqref{covariant_point_of_view}, in the covariant formulation of continuum mechanics, we assume that the Lagrangian force can be written in terms of a Lagrangian force density 
$\mathfrak{F}_{\mathcal{L}}$ as
\begin{equation}\label{force_field_density}
F_\mathsf{L} (\varphi(t),\dot{\varphi}(t)) = 
\int_B \mathfrak{F}_{\mathcal{L}}(j^1\varphi(t,s))d^n s.
\end{equation}

In general, on an arbitrary locally trivial fiber bundle 
$\pi :Y \rightarrow X$, the Lagrangian force density is a 
bundle map $\mathfrak{F}_{\mathcal{L}}:J^1Y \rightarrow V^*Y
\otimes \Lambda^{n+1}X$ covering the identity on $Y$ and 
the covariant Lagrange-d'Alembert principle may be written as
\begin{equation}\label{cov_LdA} 
\left.\frac{d}{d\varepsilon}\right|_{\varepsilon=0} \int_X 
\mathcal{L} (t, j ^1 \varphi _\varepsilon(x)  )+ 
\int_X \mathfrak{F}_\mathcal{L} (j^1 \varphi (x)) \cdot \delta \varphi =0,
\end{equation} 
for all variations $ \delta \varphi $ with $ \delta \varphi |_ { \partial X}=0$. This yields the forced CEL equations in intrinsic form
 $\mathscr{EL}(L)+ \mathfrak{F} _L=0$. In our case, since 
 $X=[0,L] \times B$, we have
\begin{equation}\label{forced_cov_EL} 
\frac{\partial }{\partial t}\frac{\partial \mathcal{L} }{\partial \dot \varphi }+ \operatorname{div}\frac{\partial \mathcal{L} }{\partial \nabla  \varphi}- \frac{\partial \mathcal{L} }{\partial \varphi }= \mathfrak{F} _ \mathcal{L} ( \varphi , \dot \varphi , \nabla \varphi ). 
\end{equation} 

\subsection{Multisymplectic forms and covariant momentum maps} \label{com_mult_form}

The goal of this subsection is to provide a quick review concerning the multisymplectic form formula.
This formula is of central importance since it generalizes to the covariant case the symplectic property of the flow of the EL equations in classical mechanics. This formula has a discrete analogue that characterizes multisymplectic integrators (\cite{MaPaSh1998}). The formula is more easily formulated by staying on an arbitrary fiber bundle and using the geometry of jet bundles rather than focusing on the case of trivial bundles.

\paragraph{Dual jet bundles.} On the Hamiltonian side, the covariant analogue of the phase space of classical mechanics is given by the dual jet bundle $J ^1 Y^\star\rightarrow Y$. Abstractly, the fiber of the dual jet bundle at $y\in Y_x $ consists of affine maps from $J^1Y_y$ to $\Lambda ^{n+1}_xX$, i.e.,
\[
J^1Y^\star_y := \mathrm{Aff}\left( J^1Y_y, \Lambda^{n+1} _xX \right).
\]
The momentum bundle is, by definition, the vector bundle $ \Pi \rightarrow Y$, whose fiber at $y$ is $ \Pi _y=L(L(T_xX, V_yY), \Lambda ^{n+1}_xX)$. There is a line bundle $ \mu : J^1Y^\star \rightarrow \Pi $ locally given by $\mu (x ^\mu ,y ^A , p ^\mu _A , p)=( x ^\mu ,y ^A , p ^\mu _A )$.

In our case, since the bundle $ Y \rightarrow X$ is trivial, the dual jet bundle can be identified with the vector bundle $T^\ast M \otimes TX \otimes \Lambda ^{n+1}X \times _{X \times M} \Lambda ^{n+1}X$ over $ X \times M$. Coordinates on the dual jet bundle are denoted $(t, s ^i , y ^A , p _A ^0 , p _A ^i , p)$ and correspond to the affine map
\[
( v _0 ^A , v _i ^A ) \mapsto (p+ p _A ^0 v ^A _0 +p _A ^i v ^A _i)dt \wedge d ^n s.
\]
Similarly, the momentum bundle can be identified with the vector bundle $T^\ast M \otimes TX \otimes \Lambda ^{n+1}X$ over $ X \times M$.

\paragraph{The Legendre transforms.} Given a Lagrangian density $ \mathcal{L} :J^1Y \rightarrow \Lambda ^{n+1}X$, the associated covariant Legendre transform is the fiber-preserving map $\mathbb{F}  \mathcal{L} : J^1Y \rightarrow J^1Y^\star$, given locally by 
\[
p _A ^0 = \frac{\partial L}{\partial v ^A _0 }, \quad  p _A ^i = \frac{\partial L}{\partial v ^A _i}, \quad p= L- \frac{\partial L}{\partial v ^A _0 }v _0 ^A -  \frac{\partial L}{\partial v ^A _i}v _i ^A.
\]
In the case of a trivial bundle $Y=X \times M$, $X= \mathbb{R}  \times B$ and in terms of a given field $ \varphi :X \rightarrow M$, we can write
\[
\mathbb{F}  \mathcal{L}(\dot \varphi , \nabla \varphi )= \left(\frac{\partial \mathcal{L} }{\partial \dot \varphi }, \frac{\partial \mathcal{L} }{\partial \nabla \varphi }, \mathcal{L} - \frac{\partial \mathcal{L} }{\partial \dot \varphi }\!\cdot\! \dot \varphi - \frac{\partial \mathcal{L} }{\partial \nabla \varphi }\!:\!\nabla \varphi  \right),
\]
where $ \cdot $ and $:$ denote contractions (on one, respectively,
two indices). Note that since $\mathrm{dim}(J^1Y^\star) = (n+1)+m+(n+1)m+1 \neq \mathrm{dim}(J^1Y)= \operatorname{dim} ( \Pi ) $, the Legendre transform can never be a diffeomorphism. Therefore, the Legendre transform is sometimes defined as the map $\widehat{ \mathbb{F}  \mathcal{L} }:= \mu \circ \mathbb{F}  \mathcal{L}: J^1Y \rightarrow \Pi $.

\paragraph{Cartan forms.} The dual jet bundle $J ^1 Y ^\star$ is naturally endowed with a canonical $(n+1)$-form $ \Theta $. By pulling back this $(n+1)$-form with the Legendre transform, we obtain the Cartan $(n+1)$-form $ \Theta _ \mathcal{L} := \mathbb{F}  \mathcal{L} ^\ast \Theta $ on $J ^1 Y$, locally given by
\begin{align*}
\Theta _ \mathcal{L} & = \frac{\partial L}{\partial v _0 ^A  } d y^{A}  \wedge d^nx_0 +  \frac{\partial L}{\partial v _i ^A  } d y^{A}  \wedge d^nx_i +  \left( L - \frac{\partial L }{\partial v _0 ^A }v _0 ^A   - \frac{\partial L }{\partial v _i ^A  }v _i ^A  \right)d^{n+1}x \\
&=\frac{\partial L}{\partial v _0 ^A  } d y^{A}  \wedge d^ns +  \frac{\partial L}{\partial v _i ^A  } d y^{A}  \wedge dt \wedge d^{n-1}s_i +  \left( L - \frac{\partial L }{\partial v _0 ^A }v _0 ^A   - \frac{\partial L }{\partial v _i ^A  }v _i ^A  \right)dt \wedge d^n s,
\end{align*}
where
\begin{align*}
 d^nx_0 & = \mathbf{ i}_{ \frac{\partial}{\partial t}} d^{n+1}x = d^n s, \quad\text{and} \\
  d^nx_i & = \mathbf{ i}_{ \frac{\partial}{\partial x^i}} d^{n+1}x= dt \wedge  d^{n-1} s _i, \quad \text{when} \ d^{n+1}x= dt\wedge d ^n s.
\end{align*}
The Cartan form allows to write the Lagrangian density evaluated on a first jet extension $j ^1 \varphi $ as
\[
\mathcal{L}(j^1\varphi) = (j^1 \varphi)^* \Theta_{\mathcal{L}}.
\]
The Cartan form naturally appears in the covariant Hamilton principle when the variations are not necessarily vanishing at the boundary. Writing $\delta \varphi (x)=V( \varphi (x))$, where $V$ is a vertical vector field on $\pi :Y\rightarrow X$, we have $\delta j ^1 \varphi(x) = j ^1 V( j^1 \varphi (x))$, where the vertical vector field $j^1V$ on $J^1Y\rightarrow X$ is the first jet extension of $V$. With these abstract notations, \eqref{computation_EL} can be written as
\begin{equation}\label{intrinsic_form_Theta} 
\left.\frac{d}{d\varepsilon}\right|_{\varepsilon=0} \int_ X \mathcal{L} ( j ^1 \varphi _\varepsilon (x))= \int_X \mathscr{EL}(\mathcal{L} ) \left(j^1\varphi(x) \right)\cdot V( \varphi(x))+\int_{ \partial X} (j^1 \varphi )^\ast \left( \mathbf{i} _{ j^1 V} \Theta _ \mathcal{L} \right).
\end{equation} 
Finally, the CEL operator can be rewritten in terms of the $(n+2)$-form $\Omega _{ \mathcal{L} }=- \mathbf{d} \Theta _ \mathcal{L} $ as
\begin{equation}\label{EL_op_Cartan} 
\mathscr{EL}(\mathcal{L} ) \left(j^1\varphi \right)\cdot V\circ \varphi= -( j ^1 \varphi ) ^\ast \left( \mathbf{i} _{j ^1 V} \Omega _ \mathcal{L} \right).
\end{equation}

\paragraph{Multisymplectic form formula.} It is well-known in classical Lagrangian mechanics that the flow $F_t$ of the EL equations is symplectic relative to the symplectic form $ \Omega_\mathsf{L}= \mathbb{F} \mathsf{L} ^\ast \Omega_ {can} $ on $TQ$, that is, we have
\[
F_t ^\ast \Omega  _\mathsf{L}= \Omega  _\mathsf{L},
\]
where we supposed that $\mathsf{L}$ is regular. In order to generalize this fact to the case of field theory, this property has to be reformulated. We follow \cite{MaPaSh1998}. Consider the action functional $S(q(\cdot ))=\int_{t _0 }^{ t _1 }\mathsf{L}(q(t),\dot q(t))dt$. We have the formula
\begin{equation}\label{formula_dS} 
\mathbf{d} S(q(\cdot )) \cdot \delta q(\cdot )= \int_{t _0 }^{ t _1 }\left( \frac{\partial \mathsf{L}}{\partial q}- \frac{d}{dt} \frac{\partial \mathsf{L}}{\partial \dot q} \right) \delta q\, dt  + \theta _\mathsf{L}( \dot q(t)) \cdot \delta \dot q (t)\Big|_{t _0 }^{ t _1 }.
\end{equation} 
Consider now the function $S_t$ defined on the space of solutions $\mathcal{C} _\mathsf{L}$ of the EL equations, which can be identified with initial conditions $v_q \in TQ$, defined by
\[
S_t( v _q):= \int_{t _0 }^t \mathsf{L}(q(s),\dot q(s))ds,
\]
where $(q(s),\dot q(s))=F_s( v _q )$.
In this case, \eqref{formula_dS} becomes
\[
\mathbf{d} S_t(v _q ) \cdot \delta v _q = \Theta _\mathsf{L}(F_t(v _q )) \cdot \delta (F_t( v _q ))- \Theta _\mathsf{L}(v _q )\cdot \delta v _q = (F_t ^\ast \Theta _\mathsf{L}- \Theta _\mathsf{L})(v _q) \cdot \delta v _q .
\]
From the formula $ \mathbf{d} S_t= F_t ^\ast \Theta _\mathsf{L}- \Theta _\mathsf{L} $, we can deduce the symplecticity of the flow since we have $0=\mathbf{d} \mathbf{d} S_t= - F_t ^\ast \Omega _\mathsf{L}+ \Omega _\mathsf{L}$ . This formula also tells us that the symplecticity of the flow is equivalent to the formula $ \mathbf{d} \mathbf{d} S_t =0$.

Going back to \eqref{formula_dS}, we observe that on the space $ \mathcal{C} $ of curves defined on $[t _0 , t _1 ]$ the formula can be rewritten as
\[
\mathbf{d} S(q( \cdot )) \cdot V= \alpha_1 (q( \cdot )) \cdot V+ \alpha _2 (q( \cdot )) \cdot V,
\]
where $V$ is an arbitrary variation of the curve $q( \cdot )$, and the one-forms $ \alpha _1 $ and $\alpha _2 $ on $ \mathcal{C} $ are defined by
\[
\alpha _1 ( q( \cdot ))\cdot \delta q( \cdot ):= \int_{t _0 }^{ t _1 }\left( \frac{\partial \mathsf{L}}{\partial q}- \frac{d}{dt} \frac{\partial \mathsf{L}}{\partial \dot q} \right) \delta q(t) dt\qquad \alpha _2 (q( \cdot )) \cdot \delta q( \cdot )= \theta _\mathsf{L}( \dot q(t)) \cdot \delta \dot q (t)\Big|_{t _0 }^{ t _1 }.
\]
From the above formula, one deduces $0= \mathbf{d} \mathbf{d} S= \mathbf{d} \alpha_1  + \mathbf{d} \alpha _2 $. Given a solution $q(t)$ of the EL equations, a first variation at $q(t)$ is a vector field $V$ on $Q$ such that $t \mapsto \left(F^V _\varepsilon \circ q\right)(t)$ is also a solution curve, where $F_ \varepsilon ^V$ is the flow of $V$. One can associate the vectors $V_{q( \cdot )}:= V \circ q( \cdot )$ at $ q( \cdot )$ on $\mathcal{C} $ also called first variations, and deduce the formula
\begin{equation}\label{symplectic_formula_mechanics} 
\mathbf{d} \alpha _2 (V_{q( \cdot )},W_{q( \cdot )})=0.
\end{equation} 
It is this formulation of symplecticity that is generalized to the case of field theories.

In the case of field theories, \eqref{intrinsic_form_Theta} can be written as
\begin{equation}\label{def_alpha} 
\mathbf{d} S( \varphi ) \cdot \delta \varphi = \alpha _1 ( \varphi ) \cdot \delta \varphi +\alpha _2 ( \varphi ) \cdot \delta \varphi,
\end{equation} 
where $ \alpha _1$, $ \alpha _2 $ are the one-forms on sections defined as $\alpha _1 ( \varphi ) \cdot \delta \varphi= -\int_X( j ^1 \varphi ) ^\ast \left( \mathbf{i} _{j ^1 V} \Omega _ \mathcal{L} \right)$ and $\alpha _2 ( \varphi ) \cdot \delta \varphi = \int_{ \partial X} (j^1 \varphi )^\ast \left( \mathbf{i} _{ j^1 V} \Theta _ \mathcal{L} \right)$. In the case of field theories, a first variation at a given solution $\varphi $ of the CEL equations, is a vertical vector field $W \in \mathfrak{X}  (Y)$ whose flow $F _\varepsilon ^W $ is such that $ F_ \varepsilon \circ \varphi $ is still a solution of the CEL equations, that is, by \eqref{EL_op_Cartan}, $j^1(F^W_ \varepsilon \circ \varphi ) ^\ast \mathbf{i} _{ j ^1 V} \Omega _ \mathcal{L} =0$, for all vertical vector field $V \in \mathfrak{X}  ^V (Y)$. Taking the $ \varepsilon $-derivative, we obtain that $W$ verifies the equation $(j^1 \varphi ) ^\ast \pounds _{j ^1 W }\mathbf{i} _{ j ^1 V}\Omega_ \mathcal{L}=0$, for all $V \in \mathfrak{X}  ^V (Y)$. From this and \eqref{intrinsic_form_Theta}-\eqref{EL_op_Cartan}, it can be shown that if $ \varphi $ is a solution of the CEL equations, then, for all first variations $V$, $W$ at $ \varphi $, we have
\begin{equation}\label{mult_form_formula}  
\mathbf{d} \alpha _2 ( \varphi ) (V\circ \varphi , W\circ \varphi )=0\quad\text{or, equivalently,}\quad \int_{ \partial X} (j^1 \varphi ) ^\ast \mathbf{i} _{ j ^1 V} \mathbf{i} _{ j ^1 W} \Omega _ \mathcal{L} =0,
\end{equation}
as shown (in a slightly more general situation) in \cite{MaPaSh1998}.
This formula is the analogue of \eqref{symplectic_formula_mechanics} for the case of field theories and is referred to as the \textit{multisymplectic form formula}. 

\color{black}

\paragraph{Covariant momentum map, and Noether theorem.} Let $G$ be a Lie group acting on $Y$ by $ \pi$-bundle automorphisms $ \eta _Y: Y \rightarrow Y$.
A Lagrangian density $\mathcal{L} $ is said to be $G$-equivariant if
\[
\mathcal{L} (j^1 \eta_Y (\gamma)) = (\eta_X^{-1})^\ast \mathcal{L} (\gamma ), \quad \text{for all} \quad  \gamma \in J ^1 Y, \quad \eta _Y \in G,
\]
where $ \eta _X :X \rightarrow X$ denotes the diffeomorphism of $X$ induced by $ \eta _Y $ and $ j ^1 \eta _Y:J ^1 Y \rightarrow J ^1 Y $ is the lifted diffeomorphism of $J ^1 Y$.
The Lagrangian momentum map associated to this action and to $\mathcal{L} $ is the map $J^ \mathcal{L} :J^1Y \rightarrow \mathfrak{g}^* \otimes \Lambda^nJ^1Y$ defined by
\begin{equation}\label{mom_map}
J^ \mathcal{L} (\xi) = \mathbf{i}_{\xi_{J^1Y}} \Theta_L,
\end{equation}
for $ \xi \in \mathfrak{g}  $ and where $ \xi _{ J ^1 Y}$ is the infinitesimal generator associated to the lifted action of $G$ on $ J ^1 Y$.
The Noether theorem can be proved by using formula \eqref{intrinsic_form_Theta} together with the $G$-equivariance of $ \mathcal{L} $. It is recalled in the following theorem. 

\begin{theorem}\label{Noether_th} Let $\mathcal{L} :J^1Y \rightarrow \Lambda ^{n+1}X$ be a $G$-equivariant Lagrangian density and let $J^ \mathcal{L}  : J^1Y \rightarrow \mathfrak{g}^* \otimes \Lambda^nJ ^1 Y$ be the associated Lagrangian momentum map. If the section $ \varphi : X\rightarrow Y$ is a solution of the CEL equations, then, for any subset $U  \subset X$ with smooth boundary, we have
\[
\int_{\partial U} (j^1\varphi )^* J^ \mathcal{L} (\xi)=0,\quad \text{for all $ \xi \in \mathfrak{g}  $}. 
\]
The associated local conservation law is
\[
\mathbf{d}\left[  (j^1\varphi )^* J^ \mathcal{L} (\xi) \right] = 0,\quad \text{for all $ \xi \in \mathfrak{g}  $}.
\]
\end{theorem}

\subsection{Covariant Euler-Lagrange equations on Lie groups}\label{Cov_EL_LG}

In this section we suppose that the fiber is a Lie group $M=G$ and use the notation $ \varphi (t,s)=g(t,s) \in G$. 
We rewrite the CEL in a trivialized formulation, since it is this form of the CEL equations that will be discretized on Lie groups.

\paragraph{Trivialization of Lie groups.} We can rewrite the CEL equations in a trivialized form by using the vector bundle isomorphism
\begin{equation} \label{trivialization}
J^1(X \times G) = L(TX, TG) \stackrel{\sim}\longrightarrow 
L(TX, \mathfrak{g}) \times G = (T^*X\otimes \mathfrak{g}) \times G,
\end{equation}
over $X \times G$, induced by the (left) trivialization $TG \ni v_g
\stackrel{\sim}\longleftrightarrow (g,g^{-1}v_g) \in G \times \mathfrak{g}$ of the tangent bundle $TG$ of $G$. Coordinates on the trivialized jet bundle are denoted $(t, s^i, g^A, \xi_0^A, \eta^A_i)$ 
and the above vector bundle isomorphism reads $(t, s ^i , 
g^{A}, v_0^A , v_i^A) \mapsto (t, s^i, g^{A}, g^{-1} v_0^A , 
g^{-1} v_i^A)$.
The induced trivialized Lagrangian density $\bar{\mathcal{L}}$ on 
$L(TX, \mathfrak{g}) \times G$ verifies
\begin{align*} 
\mathcal{L} (t,s, g(t,s), \dot g(t,s), \nabla g(t,s))&= \bar{\mathcal{L}} (t,s, g(t,s), g(t,s)^{-1} \dot{g}(t,s), g(t,s)^{-1}\nabla g(t,s))\\
&=: \bar{ \mathcal{L} }(t,s,g(t,s), \xi (t,s) , \eta (t,s)),
\end{align*} 
where $\xi(t,s): = g(t,s)^{-1} \dot{g}(t,s)$, $\eta(t,s) : = g(t,s)^{-1}\nabla g(t,s)$.

The trivialized CEL equations are obtained by applying Hamilton's principle to $ \bar{ \mathcal{L} }$ and using the variations induced on $ \xi = g ^{-1} \dot g $ and $ \eta = g ^{-1} \nabla g$, given by
\[
\delta \xi = \dot{\zeta} + [\xi, \zeta] \quad\text{and}\quad \delta \eta  = \nabla \zeta + [\eta, \zeta],
\]
where $\zeta := g ^{-1} \delta g: X \rightarrow \mathfrak{g} $ is an arbitrary map with $ \zeta |_{ \partial X}=0$.
We get the trivialized CEL equations
\begin{equation}\label{Triv_CEL} 
\frac{\partial }{\partial t}\frac{\delta \bar {\mathcal{L} }}{\delta \xi }+ \operatorname{div}  \frac{\delta \bar {\mathcal{L} } }{\delta \eta }= \operatorname{ad}^*_ \xi \frac{\delta \bar {\mathcal{L} } }{\delta \xi }+ \operatorname{ad}^*_ \eta \frac{\delta \bar {\mathcal{L} } }{\delta \eta  }+ g ^{-1} \frac{\partial  \bar {\mathcal{L} }}{\partial  g}.
\end{equation} 
Other boundary conditions can be used in the variational principle. The analogue of \eqref{BC_traction} and \eqref{BC_wird} being, respectively,
\[
\frac{\partial \bar{\mathcal{L}}}{\partial \eta _{,i}^A  }\mathbf{n} ^i=0, \quad \text{for all $A=1,...,m$}\quad\text{and}\quad \frac{\partial \bar{\mathcal{L}}}{\partial \eta ^A  }(0,s)=\frac{\partial \bar{\mathcal{L}}}{\partial \eta ^A  }(T,s)=0, \quad \text{for all $ s \in B$.} 
\]

\begin{remark}[$G$-invariance]{\rm
The Lagrangian $\mathcal{L}$ is $G$-invariant if and only if 
$\bar{\mathcal{L}}(g, \xi, \eta)= \bar{\mathcal{L}}(e, \xi, \eta)$
for any $g \in G$. In this case, instead of working with 
$\bar{\mathcal{L}}:  L(TX, \mathfrak{g}) \times G
\rightarrow \Lambda^{n+1}X$, it suffices to consider the reduced 
Lagrangian $\ell( \xi, \eta)$ associated to $\mathcal{L}$ by 
$G$-invariance, i.e., $\ell(\xi, \eta): =
\bar{\mathcal{L}}( e, \xi, \eta)$. Since 
$\frac{\partial \bar {\mathcal{L} }}{\partial g}=0$, 
the trivialized CEL equations consistently 
recover the covariant Euler-Poincar\'e equations
\[
\frac{\partial}{\partial t}\frac{\delta \ell}{\delta \xi}+ 
\operatorname{div}\frac{\delta \ell}{\delta \eta}
= \operatorname{ad}^*_ \xi \frac{\delta \ell}{\delta \xi}+ \operatorname{ad}^*_\eta \frac{\delta \ell}{\delta \eta},
\]
obtained by reduction, see \cite{CaRaSh2000}.
}
\end{remark}

\paragraph{Legendre transforms.} Analogously to \eqref{trivialization}, 
the dual jet bundle can be trivialized by using the vector bundle isomorphism
\begin{equation*}
J^1(X \times G ) ^\star= TX \otimes T^*G \otimes \Lambda ^{n+1}X 
\times _{X \times G} \Lambda ^{n+1}X \stackrel{\sim}\rightarrow 
TX \otimes \mathfrak{g}^\ast \otimes \Lambda^{n+1}X 
\times_{X \times G} \Lambda ^{n+1}X,
\end{equation*}
induced by the (left) trivialization $T^*G\ni \alpha_g \stackrel{\sim}
\longleftrightarrow (g, T_e^*L_g \alpha_g) \simeq G \times\mathfrak{g}  ^\ast $. Local coordinates on the trivialized dual jet bundle are denoted $(t, s ^i , g , \mu _A^0  , \mu _A ^i ,p)$ and the above vector bundle isomorphism reads $(t, s ^i , g , p _A ^0 , p _A ^i ,p) \mapsto(t, s ^i , g , g ^{-1} p _A ^0 , g ^{-1} p _A ^i ,p)$. Locally, the trivialized Legendre transform $\mathbb{F}  \bar{ \mathcal{L} }$ 
is the fiber bundle map over $X \times G$
\[
\mathbb{F}\bar{\mathcal{L}} : (TX\otimes \mathfrak{g}) \times G \rightarrow (TX  \otimes \mathfrak{g}^*\otimes \Lambda^{n+1}X)  \times _{X \times G} \Lambda ^{n+1}X.
\]
given by
\[
\mu _A ^0 = \frac{\delta \bar L}{\delta \xi ^A }, \quad  \mu _A ^i = \frac{\delta \bar L}{\delta \eta ^A _i  }, \quad p= \bar L- \frac{\delta \bar L}{\delta \xi ^A } \xi ^A - \frac{\delta \bar L}{\delta \eta ^A _i  } \eta ^A _i.
\]
Given a field $g: X \rightarrow G$, and defining $ \xi := g ^{-1} \dot g$, $ \eta := g ^{-1} \nabla g \in \mathfrak{g}$, we can write
\[
\mathbb{F}  \bar{ \mathcal{L} }( g, \xi , \eta )= \left(g,\frac{\delta \bar{ \mathcal{L} }}{\delta \xi } , \frac{\delta \bar{ \mathcal{L} }}{\delta \eta  }, \bar{ \mathcal{L} }- \frac{\delta \bar{ \mathcal{L} }}{\delta \xi } \xi - \frac{\delta \bar{ \mathcal{L} }}{\delta \eta  }\eta   \right).
\]
Similarly, the trivialized version of $\widehat{ \mathbb{F}  \mathcal{L} }$ reads $\widehat{ \mathbb{F}  \bar{\mathcal{L} }}( g, \xi , \eta )= \left(g,\frac{\delta \bar{ \mathcal{L} }}{\delta \xi } , \frac{\delta \bar{ \mathcal{L} }}{\delta \eta  }\right) $.

The trivialized Cartan form $ \bar{\Theta}_{\bar{\mathcal{L}}} := \mathbb{F} \bar{\mathcal{L}}^\ast \bar{\Theta}$ is found to be
\[
\bar{\Theta} _{ \bar{\mathcal{L}}} = \frac{\delta \bar{ L }}{\delta \xi^A  } g ^{-1} d g ^A  \wedge d^ns +  
\frac{\delta \bar{L}}{\delta \eta_i ^A }g ^{-1} dg ^A  
\wedge dt \wedge d^{n-1} s_i   +  \left(  \bar{ L }- 
\frac{\delta \bar{L}}{\delta \xi_i} \xi_i - 
\frac{\delta \bar{L}}{\delta \eta_i ^A}\eta_i^A \right)dt \wedge d^n s.
\]

\begin{remark}[$G$-invariance]{\rm Recall that if the Lagrangian is
$G$-invariant, we have $\bar{ \mathcal{L} }(g, \xi , \eta )=
\ell( \xi , \eta )$. Therefore, the maps $ \mathbb{F}  \bar{ \mathcal{L} }$ and $\widehat{ \mathbb{F}  \bar{\mathcal{L} }}$ yield the reduced Legendre transforms $ \mathbb{F}\ell: T^*X
\otimes \mathfrak{g}\rightarrow (TX  \otimes \mathfrak{g}^*
\otimes \Lambda^{n+1}X)  \times _X \Lambda ^{n+1}X$ and 
$\widehat{\mathbb{F}\ell}:  T^*X\otimes \mathfrak{g}\rightarrow 
TX  \otimes \mathfrak{g}^*\otimes \Lambda^{n+1}X$
given by
\[
\mathbb{F}\ell(\xi , \eta )=  \left(\frac{\delta\ell}{\delta \xi } , \frac{\delta\ell}{\delta \eta  }, \ell- \frac{\delta \ell}{\delta \xi } \xi - \frac{\delta \ell}{\delta \eta  }\eta   \right)\quad\text{and}\quad \widehat{\mathbb{F}\ell}(\xi , \eta )=  \left(\frac{\delta\ell}{\delta \xi } , \frac{\delta \ell}{\delta \eta  }\right).
\]}
\end{remark}

\section{Multisymplectic variational integrators and space/time splitting}\label{sec_MVI}

In this section we study the symplectic properties and the conservation laws of a multisymplectic integrator on a 1+1 dimensional 
spacetime discretized by triangles. This study uses a covariant point 
of view as well as time-evolutionary and space-evolutionary approaches.

In \S\ref{subsec_prelim}, we review from \cite{MaPaSh1998}, some basic facts on multisymplectic integrators, such as the discrete covariant Euler-Lagrange (DCEL) equations, the discrete Cartan forms, the notion of multisymplecticity, the discrete covariant Legendre transform, the discrete covariant momentum map, and the discrete covariant Noether theorem. We also consider the case with external forces and write the explicit expression of the discrete Noether quantity on arbitrary rectangular subdomains.
We consider three different classes of boundary conditions: the case where the configuration is prescribed at the space and time extremities, the case when the configuration is only prescribed at the temporal extremities, and the case where the configuration is only prescribed at the spatial extremity. In the last two cases, the associated discrete zero-traction boundary conditions are derived from the discrete covariant variational principle (Proposition \ref{Discete_Boundary_Conditions}).

The solution of the discrete problem can be organized in a time-evolutionary fashion, by first forming a vector of the discrete positions of all nodes at a given time, and then considering the sequence of these vectors indexed by the discrete time. 

Conversely, the solution of the discrete problem can be organized in a space-evolutionary fashion, by first forming a vector of the discrete positions of all nodes at a given space index, and then considering the sequence of these vectors. 

In \S\ref{symplectic_properties}, we study the symplectic character of these time-evolutionary and space-evolutionary discrete flows. This is done by constructing from the discrete covariant Lagrangian, two discrete Lagrangians ($\mathsf{L}_d$ and $\mathsf{N}_d$) associated to the temporal and spatial evolution, respectively. For this construction, it is assumed that the discrete covariant Lagrangian does not depend explicitly on the discrete time, resp., on the discrete space. We will carry out this study for each of the three boundary conditions mentioned before. The corresponding results are given in the six Propositions \ref{prop1}--\ref{prop5}.

In \S\ref{cov_VS_usual_momap}, we study the various Noether conservation theorems available when the discrete covariant Lagrangian density is invariant under the action of a Lie group $G$. Indeed, $G$-invariance of the Lagrangian density induces $G$-invariance of the discrete Lagrangians $\mathsf{L}_d$ and $\mathsf{N}_d$ so, besides the discrete covariant Noether theorem  for $ \mathcal{L} _d $, one can ask if the discrete Noether theorems associated to the time-evolutionary and space-evolutionary discrete flows are also verified. This depends on the boundary conditions considered, see Theorem \ref{summary_Noether}.

\subsection{Preliminaries on multisymplectic integrators}\label{subsec_prelim} 

We present below some basic facts about multisymplectic integrators, 
following \cite{MaPaSh1998}. In view of the applications we have
in mind, we assume from now that $B=[0,L]$ and hence 
$X=[0,T] \times [0,L]$ is a two-dimensional rectangle.

\subsubsection{Discrete covariant Euler-Lagrange equations and boundary conditions}\label{DMFT} 

We consider following discretization of spacetime $X$ given by
\begin{equation*}
X_d = \left\{ (j,a) \in \mathbb{Z} \times \mathbb{Z} \mid 
 j=0,..., N-1, \;a= 0,...,A-1\right\}.
\end{equation*}
This defines the triangles $\triangle_a^j$ 
by specifying their vertices as the ordered 
triples 
\[
\triangle_a^j= ((j,a),(j+1,a),(j,a+1)).
\]

\begin{figure}[ht]
\centering
\includegraphics[width=2.3 in]{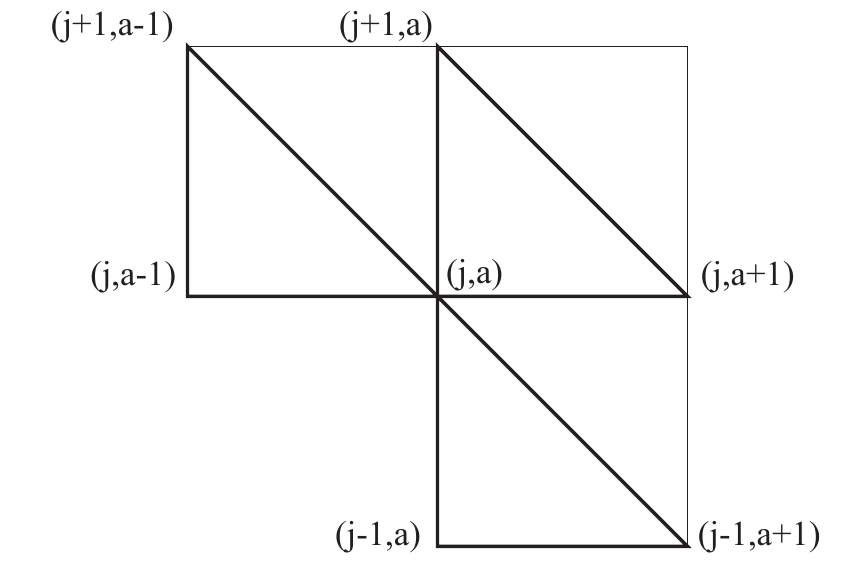}
\vspace{-3pt}
\caption{\footnotesize The triangles $\triangle_a^j$, 
$\triangle_a^{j-1}$, $\triangle_{a-1}^j$. }
\label{triangles_touching_j_a}    
\end{figure}

Denote by $ \Delta t$ and $ \Delta s$ the time and space steps. 
The discrete analogue of the bundle $Y = X \times M$ is $Y _d = 
X _d \times M$ and the discrete sections $\varphi_d$ are defined
to be maps $X _d\ni (j,a) \mapsto \varphi_d(j,a)=:\varphi _a^j \in M$.
Recall that $M$ is the space of allowed deformations;
$M = SE(3)$ for the example of the beam that we will consider in \S\ref{sect_num_ex}.
The discrete analogue of the first jet bundle is $J^1Y_d= 
X_d^\triangle \times M \times M \times M\rightarrow X_d^\triangle$, 
where $X_d^\triangle$ denotes the set of all triangles 
$\triangle_a^j$ defined above. Elements of $J^1Y_d$ are of the 
form  $(\triangle ^j _a, \varphi _a ^j , \varphi _a^{j+1}, 
\varphi _{a+1}^{j})$. The first jet extension of a discrete 
section $\varphi_d $ is the map $j^1\varphi _d : X^\triangle_d  \rightarrow X_d^\triangle \times M \times M \times M$ defined by
\[
j^1\varphi_d (\triangle) := (\triangle, \varphi_d (\triangle^{(1)}),\varphi_d (\triangle^{(2)}),\varphi_d (\triangle^{(3)})),
\]
where $\triangle^{(1)} := (j,a), \,\triangle^{(2)}:=(j+1,a), \,\triangle^{(3)}:=(j,a+1)$ and \[
\varphi_d (\triangle ^{(1)}):= \varphi_a^j, \quad 
\varphi_d (\triangle ^{(2)}):= \varphi_a^{j+1}, \quad 
\varphi_d (\triangle ^{(3)}):= \varphi_{a+1}^j.
\]

\paragraph{Discrete covariant Euler-Lagrange equations.}  A 
discrete Lagrangian density is a map 
$ \mathcal{L} _d :  
J^1Y_d \rightarrow \mathbb{R}  $ defined such that the value 
$\mathcal{L} _d (  \triangle_a^j, \varphi_a^j, \varphi_a^{j+1}, 
\varphi_{a+1}^{j})$ 
is an approximation of the integral
\begin{equation*} 
\int_{\square_a^j }L(t,s, \varphi(t,s), \partial_t  \varphi (t,s),\partial_s \varphi (t,s))dt \wedge d s,
\end{equation*}
where  $\square_a^j$ is the rectangle with vertices 
$((j,a),(j+1,a),(j,a+1), (j+1,a+1))$ and 
$\varphi : X \rightarrow M$ is a smooth map interpolating the field values $(\varphi_a^j, \varphi_a^{j+1}, \varphi_{a+1}^j)$.

The discrete action functional associated to a discrete section 
$\varphi_d$ is
\[
\mathfrak{S}_d(\varphi_d) := \sum_{j=0}^{N-1} \sum_{a=0}^{A-1} \mathcal{L} _d ( \triangle_a^j,  \varphi _a ^i , \varphi _a^{j+1}, \varphi _{a+1}^{j}).
\]

To simplify notations, we will write $\mathcal{L}_a^j := 
\mathcal{L}_d(\triangle_a^j ,\varphi_a^i , \varphi _a^{j+1}, 
\varphi_{a+1}^{j})$. Computing the derivative of the discrete 
action map (relative to $\varphi_d$) gives
\begin{align}\label{discrete_variations} 
\delta\mathfrak{S}_d(\varphi_d) & =  \sum_{j=0}^{N-1} \sum_{a=0}^{A-1} D_1\mathcal{L}_a^j \cdot \delta \varphi_a^j + D_2\mathcal{L}_a^j \cdot \delta \varphi_a^{j+1} + D_3 \mathcal{L}_a^j \cdot \delta \varphi_{a+1}^j \nonumber\\
& = \sum_{j=1}^{N-1} \sum_{a=1}^{A-1} \left[ D_1\mathcal{L}_a^j  + D_2\mathcal{L}_a^{j-1} + D_3 \mathcal{L}_{a-1}^j \right] \cdot \delta \varphi_{a}^j  \nonumber\\
& \quad + \sum_{j=1}^{N-1} \left( (D_1\mathcal{L}_{0}^j +D_2\mathcal{L}_0^{j-1})\cdot \delta \varphi_0^j + D_3 \mathcal{L}_{A-1}^j \cdot \delta \varphi_A^j \right)  \\
& \quad + \sum_{a=1}^{A-1} \left( (D_1\mathcal{L}_{a}^0 +D_3 \mathcal{L}_{a-1}^0)\cdot \delta \varphi_a^0 + D_2 \mathcal{L}_{a}^{N-1} \cdot \delta\varphi_a^N  \right) \nonumber \\
& \quad + D_1\mathcal{L}_0^0\cdot \delta \varphi_0^0 + D_2 \mathcal{L}_0^{N-1} \cdot \delta \varphi_0^N + D_3\mathcal{L}_{A-1}^0\cdot \delta \varphi_A^0 \nonumber.
\end{align}

\paragraph{(A) Discrete spacetime boundary conditions.} We shall first consider the case when the discrete configuration is known at the boundary of the spacetime domain. In this case, from the discrete 
covariant Hamilton principle, it follows tat $\delta \mathfrak{S}_d(\varphi_d) =0$, for 
all variations $\delta \varphi_a^j$ vanishing at the boundary, that 
is, such that 
\begin{align}
\label{boundary_variations}
& \delta \varphi_0^0=0, \quad \delta \varphi_a^0 = 0 , \ \ \forall a \in \{1,...,A-1\}, \quad \delta \varphi_A^0 = 0, \quad \delta \varphi_A^j = 0, \ \ \forall j \in \{1,...,N-1\}, \nonumber\\
& \delta \varphi_0^j = 0, \ \ \forall j \in \{1,..., N-1\}, \quad \delta \varphi_0^N =0, \quad \delta \varphi_{a}^N = 0, \ \ \forall a \in \{1,...,A-1\}.
\end{align} 
We thus get from \eqref{discrete_variations}  the discrete covariant
Euler-Lagrange (DCEL) equations
\begin{equation}\label{cov_DEL}
D_1 \mathcal{L} _a ^j + D_2 \mathcal{L}_a^{j-1}+ 
D_3 \mathcal{L}_{a-1}^j=0, \quad j=1,...,N-1, \;\; a=1,...,A-1,
\end{equation}
where we recall that $ \mathcal{L}_a^j := \mathcal{L}_d(
\triangle_a^j ,\varphi_a^j , \varphi_a^{j+1}, \varphi_{a+1}^{j})$ 
and the values of $\varphi_a^j$ at the boundary are prescribed. 
Note that three triangles contribute to each DCEL equation in \eqref{cov_DEL}, namely, the triangle 
$\triangle_a ^j$ associated to the vertices 
$((j,a), \, ( j+1,a), \, (j, a+1))$, the triangle 
$\triangle_a^{j-1}$ with vertices $(j-1, a), \, (j,a), \, (j-1, a+1)$, and the triangle $\triangle_{a-1}^j$ with vertices 
$(j,a-1), \, (j+1,a-1), \, (j,a)$. The intersection of the three triangles is $(j,a)$, as shown in Fig. \ref{triangles_touching_j_a}.

\paragraph{(B) Discrete boundary conditions in time.} If we assume that the discrete configuration is prescribed at $j=0$ and $j=N$, for all $a=0,...,A$, then, instead of the equalities in \eqref{boundary_variations} only the following variations vanish:
\[
\delta \varphi^0_a = \delta \varphi^N_a =0, \quad \text{for all $a=0,...,A$}.
\]
In this case, from \eqref{discrete_variations} the discrete Hamilton principle yields the boundary condition
\begin{equation}\label{resulting_BC_first} 
D_1 \mathcal{L} _0 ^j+ D_2 \mathcal{L}_0 ^{j-1}=0 \quad\text{and}\quad D_3 \mathcal{L} _{A-1}^j=0 , \quad \text{for all $j=1,...,N-1$},
\end{equation} 
referred to as the \textit{discrete zero traction boundary condition}.

\paragraph{(C) Discrete boundary conditions in space.} Conversely, if we assume that the discrete configuration is prescribed at the boundary $a=0$ and $a=A$, 
for all $j=0,...,N$, then the following variations vanish:
\[
\delta \varphi _0^j = \delta \varphi _A^j  =0, \quad \text{for all $j=0,...,N$}. 
\]
In this case, using \eqref{discrete_variations}, the discrete 
Hamilton principle yields the boundary condition
\begin{equation}\label{resulting_BC_second} 
D_1 \mathcal{L} _a ^0+ D_3 \mathcal{L}_{a-1} ^0=0 \quad\text{and}\quad D_2 \mathcal{L} _a^{N-1}=0 , \quad \text{for all $a=1,...,A-1$},
\end{equation}
referred to as the \textit{discrete zero momentum boundary condition}.

\medskip

We summarize these facts in the following proposition.

\begin{proposition}\label{Discete_Boundary_Conditions} Let $ \mathcal{L} _d : J^1Y_d \rightarrow \mathbb{R}$ be a discrete Lagrangian density. The discrete zero traction boundary conditions and zero momentum boundary conditions obtained via the covariant discrete Hamilton principle are, respectively, given by
\[
D_1 \mathcal{L} _0 ^j+ D_2 \mathcal{L}_0 ^{j-1}=0 \quad\text{and}\quad D_3 \mathcal{L} _{A-1}^j=0 , \quad \text{for all} \quad j=1,...,N-1
\]
and 
\[
D_1 \mathcal{L} _a ^0+ D_3 \mathcal{L}_{a-1} ^0=0 \quad\text{and}\quad D_2 \mathcal{L} _a^{N-1}=0 , \quad \text{for all} \quad a=1,...,A-1.
\]
\end{proposition}

\begin{remark}{\rm It is important to note that the discrete boundary conditions above are obtained exactly in the same way as their continuous counterparts, namely, they arise as boundary terms in the variational principles. These boundary terms do not contribute when the configuration is prescribed at the boundary, since the corresponding variations vanish on the boundary. However, when the configuration is not prescribed at the boundary, the variational principles yield "natural" boundary conditions, given here by \eqref{resulting_BC_first} and/or \eqref{resulting_BC_second}.
}
\end{remark}

\subsubsection{Discrete Cartan forms and multisymplecticity}\label{DCFM} 

The discrete Cartan forms, denoted $\Theta_{\mathcal{L}}^{1}, 
\Theta_{\mathcal{L}}^{2},\Theta_{\mathcal{L}}^{3}$, are the 
one-forms on $J^1Y_d$, defined by
\begin{equation}\label{DCFs}
\begin{aligned}
\Theta_{\mathcal{L}_d }^{1}(\triangle _a ^j , \varphi _a ^j , 
\varphi _a^{j+1}, \varphi _{a+1}^{j})  & := 
D_1\mathcal{L}_d(\triangle _a ^j ,\varphi_a^j , \varphi_a^{j+1}, 
\varphi_{a+1}^{j}) d \varphi_a^j  , \\
\Theta_{\mathcal{L}_d }^{2}(\triangle_a^j ,\varphi_a^j , 
\varphi_a^{j+1}, \varphi_{a+1}^{j})  & :=  
D_2 \mathcal{L}_d(\triangle_a^j ,\varphi_a^j , \varphi_a^{j+1}, 
\varphi_{a+1}^{j}) d \varphi_a^{j+1} , \\
\Theta_{\mathcal{L}_d }^{3}(\triangle_a^j ,\varphi_a^j , 
\varphi _a^{j+1}, \varphi _{a+1}^{j}) & :=  
D_3\mathcal{L}_d(\triangle_a^j ,\varphi_a^j , \varphi_a^{j+1}, 
\varphi _{a+1}^{j}) d \varphi_{a+1}^j ,
\end{aligned}
\end{equation} 
for all $(\triangle_a^j , \varphi_a^j , \varphi_a^{j+1}, 
\varphi _{a+1}^{j}) \in J^1 Y_d = X^\triangle_d \times M \times M 
\times M$, see  \cite{MaPaSh1998}. Note that one can also interpret 
the Cartan forms as $X^\triangle_d$-dependent one-forms on 
$M \times M \times M$. Viewed this way, they verify the relation
\begin{equation} \label{one_form_relationship}
\mathbf{ d} \mathcal{L}_d = \Theta_{\mathcal{L}_d }^{1} + \Theta_{\mathcal{L}_d }^{2} + \Theta_{\mathcal{L}_d }^{3}.
\end{equation}
The discrete Cartan 2-forms are defined as $ \Omega _ \mathcal{L} ^k := -\mathbf{d}  \Theta ^k _ \mathcal{L} $ and thus verify
\[
\Omega _{\mathcal{L}_d }^{1} + \Omega _{\mathcal{L}_d }^{2} + \Omega _{\mathcal{L}_d }^{3}=0.
\]
The definition of the discrete Cartan forms $ \Theta ^k _\mathcal{L} $ is motivated by the following observation.

Given a vector field $V$ on $Y_d$, we denote by $V_a^j $ its 
restriction to the fiber at $(j,a)$. Its first jet extension is 
the vector field $j^1V$ on $J^1Y_d$ defined by $j^1V(
\triangle_a^j, \varphi_a^j ,\varphi_a^{j+1}, \varphi_{a+1}^{j} )
:=(V_a^j (\varphi_a^j) ,V_a^{j+1}(\varphi_a^{j+1}),V_{a+1}^{j} 
(\varphi _{a+1}^{j}))$. Using these notations, we can now rewrite 
the variations of the discrete action \eqref{discrete_variations} 
in a way analogous to \eqref{intrinsic_form_Theta}. Namely, given 
a discrete field $\varphi_d $ with variations $\delta \varphi_d$, defining the vector field $V$ on $Y_d $ such that 
$\delta \varphi_a^j = V(\varphi_a^j)$,  and rewriting 
\eqref{DCFs} in the form
\begin{align*} 
D_1\mathcal{L}_d(\triangle _a ^j ,\varphi _a ^j , \varphi _a^{j+1}, \varphi _{a+1}^{j}) \delta  \varphi_a^j&= \left[ (j^1 \varphi _d ) ^\ast (\mathbf{i} _{ j ^1 V} \Theta _ { \mathcal{L} _d } ^1 )\right] ( \triangle_a ^j ),\\
D_2\mathcal{L}_d(\triangle _a ^j ,\varphi _a ^j , \varphi _a^{j+1}, \varphi _{a+1}^{j}) \delta  \varphi_a^{j+1}&= \left[ (j^1 \varphi _d ) ^\ast (\mathbf{i} _{ j ^1 V} \Theta _ { \mathcal{L} _d } ^2 )\right] ( \triangle_a ^j ),\\
D_3\mathcal{L}_d(\triangle _a ^j ,\varphi _a ^j , \varphi _a^{j+1}, \varphi _{a+1}^{j}) \delta  \varphi_{a+1}^j&= \left[ (j^1 \varphi _d ) ^\ast (\mathbf{i} _{ j ^1 V} \Theta _ { \mathcal{L} _d } ^3 )\right] ( \triangle_a ^j ),
\end{align*}
equality \eqref{discrete_variations} becomes
\begin{equation} \label{reformulated_discrete_action}
\begin{aligned}
\delta \mathfrak{S}_d(\varphi_d) =&\sum_{j=1}^{N-1} \sum_{a=1}^{A-1} 
\left[ D_1\mathcal{L}_a^j  + D_2\mathcal{L}_a^{j-1} + 
D_3 \mathcal{L}_{a-1}^j \right] \cdot \delta \varphi_{a}^j\\
 &+ \sum_{\{\triangle \in X_d^ \triangle\mid \triangle \cap 
 \partial X_d \neq \varnothing\}} \left(\sum_{k \in \{1,2,3\}; 
 \triangle^{(k)}\in \partial X_d} \left[(j^1\varphi_d)^* 
(\mathbf{i}_{j^1V} \Theta_{\mathcal{L}_d}^k)\right](\triangle)\right)\\
&=: \alpha_1 (\varphi_d) \cdot \delta \varphi_d+ 
\alpha _2 (\varphi_d) \cdot \delta \varphi_d.
\end{aligned} 
\end{equation} 

The one-forms $\alpha_1$ and $\alpha_2$ on the space of 
discrete sections are defined analogously with \eqref{def_alpha}, see \cite{MaPaSh1998} for details. When evaluated on first variations $V$, $W$ at a solution $ \varphi _d $, the formula $0= \mathbf{d} \mathbf{d} \mathfrak{S}_d = \mathbf{d} \alpha _1 + \mathbf{d} \alpha _2 $ yields $ \mathbf{d} \alpha _2 ( \varphi _d )(V,W)=0$, or equivalently
\begin{equation}\label{discrete_multisymp_formula} 
\sum_{\{\triangle \in X_d^ \triangle\mid \triangle \cap 
\partial X_d \neq \varnothing\}} \left(\sum_{k \in \{1,2,3\}; \triangle^{(k)}\in \partial X_d} 
\left[(j^1\varphi_d )^* (\mathbf{i}_{j^1V} \mathbf{i}_{j^1W} 
\Omega_{\mathcal{L}_d}^k) \right](\triangle)\right)=0.
\end{equation} 
This formula is referred to as the \textit{discrete multisymplectic form formula}. It is the discrete version of \eqref{mult_form_formula} and generalizes the notion of symplecticity for integrators in mechanics to the case of integrators in field theory.

\paragraph{Discrete covariant Legendre transform.} The discrete covariant Legendre transforms are the maps $\mathbb{F}^k \mathcal{L}_d:J^1Y_d \rightarrow T^*M$ given by
 \begin{equation}\label{DCLegendre_transf}
\begin{aligned}
\mathbb{F}^1 \mathcal{L}_d ( \triangle _a ^j , \varphi _a ^j , \varphi _a^{j+1}, \varphi _{a+1}^{j} ) & = \left(\varphi_a^j, \, D_1\mathcal{L}_a^j \right), \\
\mathbb{F}^2 \mathcal{L}_d (\triangle _a ^j , \varphi _a ^j , \varphi _a^{j+1}, \varphi _{a+1}^{j}) & = \left( \varphi_a^{j+1},\, D_2 \mathcal{L}_a^j \right), \\
\mathbb{F}^3 \mathcal{L}_d (\triangle _a ^j , \varphi _a ^j , \varphi _a^{j+1}, \varphi _{a+1}^{j} ) & = \left( \varphi_{a+1}^j,\, D_3\mathcal{L}_a^j\right).
\end{aligned}
\end{equation}
We note that the DCEL equations can be thus written in the form
\begin{equation*} 
\mathbb{F}^1\mathcal{L}_d (\triangle _a ^j  \varphi _a ^j , \varphi _a^{j+1}, \varphi _{a+1}^{j}) + \mathbb{F}^2\mathcal{L}_d (\triangle _a ^{j-1}  \varphi _a ^{j-1} , \varphi _a^{j}, \varphi _{a+1}^{j-1}) + \mathbb{F}^3\mathcal{L}_d (\triangle _{a-1}^j  \varphi _{a-1}^j , \varphi _{a-1}^{j+1}, \varphi _{a}^{j})=0,
\end{equation*}
which can be regarded as a matching of momenta in $T^*_{\varphi_a^j}M$.

\subsubsection{Discrete covariant momentum maps}\label{subsubsec_momap} 

We consider only vertical symmetries, that is, group actions that act trivially on the base $X$. Let $ \Phi :G \times M \rightarrow M$ be a left action of a Lie group $G$ on $M$. This action naturally induces an action on the discrete jet bundle, given by
\[
\Phi _g ^{J^1Y_d}( \triangle _a ^j ,\varphi_a^j, \varphi_a^{j+1}, \varphi_{a+1}^j):= ( \triangle _a ^j ,\Phi _g (\varphi_a^j), \Phi _g (\varphi_a^{j+1}),\Phi _g ( \varphi_{a+1}^j)),\quad g\in G,
\]
whose infinitesimal generator is
\[
\xi _{J^1Y_d}( \triangle _a ^j ,\varphi_a^j, \varphi_a^{j+1}, \varphi_{a+1}^j):= ( \triangle _a ^j ,\xi_M  (\varphi_a^j), \xi_M  (\varphi_a^{j+1}),\xi_M  ( \varphi_{a+1}^j)).
\]
We say that the discrete Lagrangian is invariant with respect to  this action if $ \mathcal{L} _d \circ \Phi ^{J ^1Y _d }_g= \mathcal{L} _d $, for all $ g \in G$. As a consequence, we have the infinitesimal invariance $\mathbf{d} \mathcal{L} _d \cdot \xi _{J ^1 Y _d }=0$.

The discrete momentum maps are defined by
\[
J ^k_{ \mathcal{L} _d } : J ^1 Y _d \rightarrow \mathfrak{g}  ^\ast , \quad \left\langle J ^k _{ \mathcal{L}_d  }, \xi \right\rangle :=\mathbf{i} _{ \xi _{ J ^1 Y _d }} \Theta ^k _{ \mathcal{L} },\quad \xi  \in \mathfrak{g},
\]
so we have the formulas
\begin{equation}\label{discrete_momap}
\begin{aligned}
\left\langle J ^1_{ \mathcal{L} }  ( \triangle _a ^j ,\varphi_a^j, \varphi_a^{j+1}, \varphi_{a+1}^j), \xi \right\rangle &=D _1 \mathcal{L} _d ( \triangle _a ^j ,\varphi_a^j, \varphi_a^{j+1}, \varphi_{a+1}^j) \cdot \xi _M ( \varphi _a  ^j )\\
\left\langle J ^2 _{ \mathcal{L} } ( \triangle _a ^j ,\varphi_a^j, \varphi_a^{j+1}, \varphi_{a+1}^j), \xi \right\rangle &=D _2 \mathcal{L} _d ( \triangle _a ^j ,\varphi_a^j, \varphi_a^{j+1}, \varphi_{a+1}^j) \cdot \xi _M ( \varphi _a  ^{j+1} )\\
\left\langle J ^3 _{ \mathcal{L} } ( \triangle _a ^j ,\varphi_a^j, \varphi_a^{j+1}, \varphi_{a+1}^j), \xi \right\rangle &=D _3 \mathcal{L} _d ( \triangle _a ^j ,\varphi_a^j, \varphi_a^{j+1}, \varphi_{a+1}^j) \cdot \xi _M ( \varphi _{a +1} ^j ).
\end{aligned}
\end{equation}  
We note that the infinitesimal invariance of $ \mathcal{L} _d $ can be rewritten as
\begin{equation}\label{local_DNT} 
\left( J ^1_{ \mathcal{L} _d }  +J ^2_{ \mathcal{L} _d } +J ^3_{ \mathcal{L} _d }  \right)(j^1  \varphi _d (\triangle_a^j))=0,
\end{equation} 
for all $j=0,...,N-1$ and $a=0,...,A-1$. This is the statement of the local discrete Noether theorem. To obtain the global discrete Noether theorem, one applies the formula \eqref{reformulated_discrete_action} for variations induced by the group action.
More generally, given a restriction 
\[
\mathfrak{S}_d^U(\varphi_d) := \sum_{j=K}^{L} \sum_{a=B}^{C} \mathcal{L} _d ( \triangle_a^j,  \varphi_a^i, \varphi _a^{j+1}, \varphi_{a+1}^{j})
\]
of the action functional to a subdomain $U$ of $X_d$
given by $U:=\{(a,j)\mid K\leq j\leq L, \;\; B\leq a\leq C\}$, that is,
$U$ is union of triangles whose lower left vertex belongs to a given
rectangular subdomain, and by applying formula \eqref{reformulated_discrete_action}, we get
the following result.

\begin{theorem}[Discrete global Noether theorem]\label{Cov_Noether}  Suppose that the discrete Lagrangian $ \mathcal{L} _d : J ^1 Y _d \rightarrow \mathbb{R}  $ is invariant under the action of a Lie group $G$ on $M$. Suppose that $ \varphi_d  $ is a solution of the DCEL equations for $ \mathcal{L} _d $. Then, for all $0\leq B<C\leq A-1$, $0\leq K<L\leq N-1$, we have the conservation law
\begin{equation}\label{DCN}  
\mathscr{J}_{B,C}^{K,L}(\varphi_d ) = 0,
\end{equation}
where,
\begin{align}\label{Def_DN}  
\mathscr{J}_{B,C}^{K,L}(\varphi_d ) :=&\sum_{j=K+1}^L \left( J ^1_{\mathcal{L}_d}( j ^1 \varphi _d (\triangle _B ^j )+J ^2_{\mathcal{L}_d}(  j ^1 \varphi _d (\triangle _B ^{j-1}) )+ J ^3_{\mathcal{L}_d}( j ^1 \varphi _d (\triangle _C ^j )) \right) \nonumber\\
&+ \sum_{a=B+1}^C \left( J ^1_{\mathcal{L}_d}( j ^1 \varphi _d (\triangle _a ^K )+J ^2_{\mathcal{L}_d}(  j ^1 \varphi _d (\triangle _a ^L) )+ J ^3_{\mathcal{L}_d}( j ^1 \varphi _d (\triangle _{a-1} ^K )) \right) \\
&+ J ^1_{\mathcal{L}_d}( j ^1 \varphi _d (\triangle _B ^K )+J ^2_{\mathcal{L}_d}(  j ^1 \varphi _d (\triangle _B ^L) )+ J ^3_{\mathcal{L}_d}( j ^1 \varphi _d (\triangle _C  ^K )).\nonumber
\end{align}
\end{theorem}

\medskip

Of course, the expression for $\mathscr{J}^{K,L}_{B,C}$ can be written in a condensed form like the one appearing in \eqref{reformulated_discrete_action} by using the discrete Cartan forms $\Theta_{\mathcal{L}_d}^k$. 

\subsubsection{Discrete covariant Euler-Lagrange equations with forces }

Given a Lagrangian force field density $\mathfrak{F}_{\mathcal{L}}$, 
the discrete Lagrangian forces are maps $F_d^k : J^1  Y _d  \rightarrow T^*M$, $k=1,2,3$,
with $F^1_d (j ^1 \varphi _d (\triangle_a^j)) \in T_{\varphi_a^j}^*M, \, F^2_d (j ^1 \varphi _d (\triangle_a^j)) \in 
T_{\varphi_a^{j+1}}^*M$, $F^3_d (j ^1 \varphi _d (\triangle_a^j)) \in T_{\varphi_{a+1}^j}^*M$, such that the approximation
\begin{align*}
& \int_{\square_a^j} \mathfrak{F}_\mathcal{L} (\varphi, \dot{\varphi}, 
\partial_s \varphi)dt \wedge ds\\
&\approx \sum_{j=0}^{N-1} 
\sum_{a=0}^{A-1}\left(F^1_d (j ^1 \varphi _d (\triangle_a^j)) \delta \varphi_a^j
+ F ^2 _d (j ^1 \varphi _d (\triangle_a^j)) \delta \varphi_a^{j+1}
+ F^3_d (  j^1 \varphi _d (\triangle_a^j) ) \delta \varphi_{a+1}^j\right)
\end{align*}
holds, where $\varphi : X \rightarrow M$ is a smooth map interpolating 
the field values $(\varphi_a^j, \varphi_a^{j+1}, \varphi_{a+1}^j)$.

The discrete version of the covariant Lagrange-d'Alembert principle \eqref{cov_LdA} is, therefore, 
\begin{align}
\label{Discrete_CLdA}
&\delta \sum_{j=0}^{N-1} \sum_{a=0}^{A-1} \mathcal{L} _d (j ^1 \varphi _d (
\triangle_a^j))  + \sum_{j=0}^{N-1} \sum_{a=0}^{A-1}\Big( 
F^1_d (j ^1 \varphi _d (\triangle_a^j))\delta \varphi_a^j +
 F^2 _d (j ^1 \varphi _d (\triangle_a^j))\delta \varphi_a^{j+1}
 \nonumber \\
 & \qquad \qquad \qquad \qquad \qquad \qquad \qquad \qquad + 
F^3_d (j ^1 \varphi _d (\triangle_a^j)) \delta \varphi_{a+1}^j\Big) =0.
\end{align} 
For variations vanishing at the boundary, this principle yields 
the forced DCEL equations
\begin{align*}& D_1 \mathcal{L} _a ^j + D_2 \mathcal{L} _a ^{j-1}+ D_3 \mathcal{L} _{a-1}^j   + F_d^1(j ^1 \varphi _d (\triangle_a^j))+ F_d^2(j ^1 \varphi _d (\triangle_a^{j-1})) + F_d^3(j ^1 \varphi _d (\triangle_{a-1}^j))=0, \\
&  \text{with} \quad j=1,...,N-1, \;\; a=1,...,A-1.
\end{align*}
 
Given a left action $ \Phi :G \times M \rightarrow M$, we say that the discrete forces are orthogonal to the $G$-group action, if
\begin{align}\label{orthogonal_forces} 
&\left\langle F_d^1(j^1 \varphi_d (\triangle_a^j)), 
\xi _M(\varphi _a^j) \right\rangle + 
\left\langle F_d^2(j^1 \varphi_d (\triangle_a^j)), 
\xi_M(\varphi_a^{j+1}) \right\rangle \nonumber \\
& \qquad + 
\left\langle F_d^3(j^1 \varphi_d (\triangle_a^j)), 
\xi_M(\varphi_{a+1}^j) \right\rangle=0,
\end{align} 
for all $ \xi \in \mathfrak{g}  $. In this case, it is easy to extend Theorem \ref{Cov_Noether} to the forced case, as follows.

\begin{theorem}[Discrete global Noether theorem with forces]\label{Cov_Noether_force}  Suppose that the discrete Lagrangian $ \mathcal{L} _d : J ^1 Y _d \rightarrow \mathbb{R}  $ is invariant under the action of a Lie group $G$ on $M$ and suppose that \eqref{orthogonal_forces} holds. Suppose that $ \varphi_d  $ is a solution of the discrete covariant Lagrange-d'Alembert equations for $ \mathcal{L} _d $. Then, for all $0\leq B<C\leq A-1$, $0\leq K<L\leq N-1$, we have the conservation law
\begin{equation}\label{DCN_force}  
(\mathscr{J}^F)_{B,C}^{K,L}(\varphi_d ) = 0,
\end{equation}
where $\mathscr{J}^F$ is defined as in \eqref{Def_DN}, except that $J_{ \mathcal{L} _d }^k$ are replaced by the forced covariant momentum maps $J_{\mathcal{L}_d}^{k,F}:J ^1 Y _d \rightarrow \mathfrak{g}  ^\ast $, $k=1,2,3$ defined by
\begin{align*} 
\left\langle J^{1,F}_{\mathcal{L}}(j^1 \varphi_d (\triangle_a^j )), 
\xi \right\rangle&:=\left\langle (D_1\mathcal{L}_d +F_d^1)
(j^1 \varphi_d (\triangle_a^j)), \xi_M(\varphi_a^j)\right\rangle\\
\left\langle J^{2,F}_{\mathcal{L}}(j^1 \varphi_d (\triangle _a ^j )), \xi \right\rangle&:=\left\langle (D _2 \mathcal{L} _d +F_d^2)(j ^1 \varphi _d (\triangle_a^j)), \xi _M( \varphi _a ^{j+1} )\right\rangle\\
\left\langle J ^{3,F}_{ \mathcal{L} }  ( j ^1 \varphi _d (\triangle _a ^j )), \xi \right\rangle&:=\left\langle (D _3 \mathcal{L} _d +F_d^3)(j ^1 \varphi _d (\triangle_a^j)), \xi _M( \varphi _{a+1} ^j )\right\rangle .
\end{align*} 
\end{theorem}

\subsection{Symplectic properties of the time and space discrete evolutions}\label{symplectic_properties}

In this subsection we study the symplectic character of the time-evolutionary and space-evolutionary discrete flows built from a discrete solution section $ \varphi _d $ of the DCEL equations associated to $ \mathcal{L} _d $. From  \S\ref{DCFM}, we already know that $ \varphi _d $ verifies the discrete multisymplectic form formula \eqref{discrete_multisymp_formula}, which is the multisymplectic generalization of symplecticity. However, this does not guarantee that the time-evolutionary and space-evolutionary discrete flows are symplectic. As we will explain in detail, the conclusion depends on the type of boundary conditions considered.

We will study symplecticity in time and in space by constructing from the discrete covariant Lagrangian, two discrete Lagrangians $\mathsf{L}_d(\boldsymbol{\varphi}^j, \boldsymbol{\varphi}^{j+1})$ and $\mathsf{N}_d ( \boldsymbol{\varphi} _a , \boldsymbol{\varphi} _{a+1})$ associated to the temporal and spatial evolution, respectively. For this construction, it is assumed that the discrete covariant Lagrangian does not depend explicitly on the discrete time, resp., on the discrete space. Knowing that the DEL equations associated to $\mathsf{L}_d$ and $\mathsf{N}_d$ yield symplectic scheme, we will study in details the relation between these two DEL equations and the CDEL for $ \mathcal{L} _d $. The answer depends on the class of boundary conditions considered.

Of course, when $ \mathcal{L} _d $ depends explicitly on discrete time, the discrete Lagrangian $\mathsf{L}_d $ can be defined in the same way, but it will be time dependent. The same remark applies to $\mathsf{N}_d$ and the dependence on discrete space.
 
\subsubsection{Discrete time evolution: the discrete Lagrangian $\mathsf{L}_d$}\label{DTE} 

The configuration space for the discrete Lagrangian 
$\mathsf{L}_d$ is $M^{A+1}$. Using the notation 
$\boldsymbol{\varphi}^j := (\varphi_0 ^j, ..., \varphi^j_A) 
\in M^{A+1}$, the discrete Lagrangian 
$\mathsf{L}_d: M^{A+1}\times M^{A+1} \rightarrow \mathbb{R}$ 
is defined by
\[
\mathsf{L}_d(\boldsymbol{\varphi}^j, \boldsymbol{\varphi}^{j+1})
:= \sum_{a=0}^{A-1} \mathcal{L}_d(\triangle_a^j, \varphi_a^j, \varphi_a^{j+1},\varphi_{a+1}^j),
\]
so that the associated discrete action is
\begin{equation}\label{discrete_action_L} 
\mathfrak{S}_d(\varphi_d) = \sum_{j=0}^{N-1} \mathsf{L}_d(\boldsymbol{\varphi}  ^j , \boldsymbol{\varphi}  ^{j+1}) = \sum_{j=0}^{N-1} \sum_{a=0}^{A-1} \mathcal{L}_d(\triangle_a^j, \varphi_a^j, \varphi_a^{j+1},\varphi_{a+1}^j).
\end{equation} 
In order to analyze the relation between the discrete Hamilton principles associated to $ \mathsf{L}_d$ and $ \mathcal{L} _d $, we shall first assume that there are no boundary conditions, so that the discrete Hamilton principle for $\mathsf{L}_d$ yields the stationarity conditions
\begin{equation}\label{stat_dynamic} 
D_1\mathsf{L}^j_d+D_2\mathsf{L}^{j-1}_d =0, \quad j=1,...,N-1,\quad\text{and}\quad D_1\mathsf{L}_d^0=0, \quad D_2\mathsf{L}_d^{N-1}=0,
\end{equation}
since the variations do not vanish at the boundary.
Computing these expressions in terms of $ \mathcal{L} _d $, we get
\begin{equation}\label{computations_derivatives}
\begin{aligned}  
\left\langle D_1\mathsf{L}_d(\boldsymbol{\varphi}  ^j , \boldsymbol{\varphi}  ^{j+1}), \delta \boldsymbol{\varphi} ^j \right\rangle& =\sum_{a=1}^{A-1}( D_1 \mathcal{L} _a ^j + D_3 \mathcal{L} _{a-1} ^j ) \delta \varphi _a ^j + D_1 \mathcal{L} _0^j \delta \varphi _0 ^j + D_3 \mathcal{L} _{A-1} ^j \delta \varphi _A ^j \\
\left\langle D_2\mathsf{L}_d(\boldsymbol{\varphi}  ^j , \boldsymbol{\varphi}  ^{j+1}), \delta \boldsymbol{\varphi} ^{j+1} \right\rangle &=\sum_{a=1}^{A-1} D_2 \mathcal{L} _a ^j \delta \varphi _a ^{j+1} + D_2 \mathcal{L} _0 ^j \delta \varphi _0 ^{j+1}.\\
\end{aligned} 
\end{equation} 
So the DEL for $\mathsf{L}_d$ in \eqref{stat_dynamic}  yield the equations
\begin{equation}\label{part1} 
\begin{aligned} 
&D_1 \mathcal{L} _a ^j  +D_2 \mathcal{L} _a ^{j-1} + D_3 \mathcal{L} _{a-1} ^j=0, \quad j=1,...,N-1,\;\; a=1,...,A-1\\
&D _1 \mathcal{L} _0 ^j + D _2 \mathcal{L} _0 ^{j-1}=0, \quad j=1,...,N-1\\
&D_3 \mathcal{L} _{A-1} ^j =0, \quad j=1,...,N-1
\end{aligned}
\end{equation} 
and the boundary conditions in \eqref{stat_dynamic} imply the equations
\begin{equation}\label{part2}
\begin{aligned} 
&D_1 \mathcal{L} _a ^0 + D_3 \mathcal{L} _{a-1} ^0=0, \quad a=1,...,A-1\\
&D_1 \mathcal{L}_0^0=0, \quad D_3 \mathcal{L} _{A-1} ^0=0, \quad D_2 \mathcal{L} _a ^{N-1}=0, \quad a=0,...,A-1.
\end{aligned}
\end{equation} 
Of course, \eqref{part1}-\eqref{part2} agree with the stationarity condition \eqref{discrete_variations} obtained from the discrete covariant Hamilton's principle when no boundary condition is imposed. Moreover, this computation shows that the DEL for $\mathsf{L}_d$ (i.e., the first equation in \eqref{stat_dynamic}) is equivalent to the DCEL equations for $ \mathcal{L} _d $ together with the discrete zero traction boundary conditions \eqref{resulting_BC_first} (the second and third lines in \eqref{part1}).

\begin{remark}[Discrete Cartan forms]\label{DCF_L} {\rm 
We now describe the relation between the two discrete Cartan forms 
$\Theta^\pm_{\mathsf{L}_d}$ associated to $\mathsf{L}_d$ and the 
three discrete Cartan forms $\Theta_{\mathcal{L}_d}^k$, 
$k=1,2,3$, associated to $\mathcal{L}_d$. On $M^{A+1}\times M^{A+1}$ 
we have
\begin{align*}
\Theta _{\mathsf{L}_d}^-( \boldsymbol{\varphi} ^j , \boldsymbol{\varphi} ^{j+1})&=-D_1\mathsf{L}_d ( \boldsymbol{\varphi} ^j , \boldsymbol{\varphi} ^{j+1})d\boldsymbol{\varphi} ^{j}\\
&= - \sum_{a=1}^{A-1}( D_1 \mathcal{L} _a ^j + D_3 \mathcal{L} _{a-1} ^j ) \mathbf{d}  \varphi _a ^j + D_1 \mathcal{L} _0^j \mathbf{d}  \varphi _0 ^j + D_3 \mathcal{L} _{A-1} ^j d \varphi _A ^j\\
&=- \sum_{a=0}^{A-1} \left( \Theta _ {\mathcal{L} _d }^1 (\triangle_a^j, \varphi_a^j, \varphi_a^{j+1},\varphi_{a+1}^j)+ \Theta _ {\mathcal{L}_d } ^3 (\triangle_a^j, \varphi_a^j, \varphi_a^{j+1},\varphi_{a+1}^j)\right) \\
\Theta _{\mathsf{L}_d}^+( \boldsymbol{\varphi} ^j , \boldsymbol{\varphi} ^{j+1})&=D_2\mathsf{L}_d ( \boldsymbol{\varphi} ^j , \boldsymbol{\varphi} ^{j+1}) d\boldsymbol{\varphi} ^{j+ 1}= \sum_{a=0}^{A-1} D_2 \mathcal{L} _a ^j d\varphi _a ^{j+1}\\
&= \sum_{a=0}^{A-1} \Theta _{ \mathcal{L}_d }^2 (\triangle_a^j, \varphi_a^j, \varphi_a^{j+1},\varphi_{a+1}^j).
\end{align*} 
We abbreviate these relations as
\[
\Theta _{\mathsf{L}_d}^-= - \sum_{a=0}^{A-1}\left(  \Theta^1 _{ \mathcal{L}_d }(\triangle_a ^j )+  \Theta ^3_{ \mathcal{L}_d }(\triangle_a ^j ) \right)  \quad\text{and}\quad \Theta _{\mathsf{L}_d}^+= \sum_{a=0}^{A-1} \Theta^2 _{ \mathcal{L}_d}(\triangle_a ^j ).
\]
Applying identity \eqref{one_form_relationship} to the forms $ \Theta _ {\mathcal{L}_d } ^k $, we consistently recover the relation $\mathbf{d} \Theta _{\mathsf{L}_d}^-= \mathbf{d} \Theta _{\mathsf{L}_d}^+$:
\begin{equation}\label{relation_thetas} 
\mathbf{d} \Theta _{\mathsf{L}_d}^- = - \mathbf{d} \sum_{a=0}^{A-1}\Theta _ {\mathcal{L}_d } ^1 (\triangle_a ^j )- \mathbf{d} \sum_{a=0}^{A-1}\Theta _{ \mathcal{L}_d } ^3 (\triangle_a ^j )= \mathbf{d} \sum_{a=0}^{A-1}\Theta _ {\mathcal{L} _d }^2 (\triangle_a ^j )= \mathbf{d} \Theta _{\mathsf{L}_d}^+.
\end{equation} 
The discrete Cartan 2-form $ \Omega _{\mathsf{L}_d}$ is thus related to the 2-forms $ \Omega  _{ \mathcal{L}_d  }$ as
\begin{equation}\label{relation_omegas} 
\Omega _{\mathsf{L}_d}=  \sum_{a=0}^{A-1}\Omega ^2_ {\mathcal{L}_d }(\triangle_a ^j ) =-  \sum_{a=0}^{A-1}(\Omega ^1 _ {\mathcal{L}_d }(\triangle_a ^j ) + \Omega ^3 _{\mathcal{L}_d }(\triangle_a ^j )).
\end{equation} 
\noindent\textsf{Notation.} Note that we maintained the index $j$ in the expressions above, for consistency with all the notations used in the paper. The index is, however, not needed here since all the triangles are identical. One can for example write the expression of 
$\Theta_{\mathsf{L}_d}^-$ as
\[
\Theta _{\mathsf{L}_d}^-(\boldsymbol{\varphi}, \boldsymbol{\psi})
= - \sum_{a=0}^{A-1}\left(\Theta_{\mathcal{L}_d}^1 (\varphi_a, 
\psi_a,\varphi_{a+1})+ \Theta_{ \mathcal{L}_d}^3 (\varphi_a, 
\psi_a,\varphi_{a+1})\right),
\]
for $ \boldsymbol{\varphi} =(\varphi_0 ,..., \varphi_A), \boldsymbol{\psi} = (\psi_0, ...,\psi_A) \in M^{A+1} $. Similarly 
for $\Theta_{\mathsf{L}_d}^+$ and $\Omega_{\mathsf{L}_d}$.
}
\end{remark}

\paragraph{(A) Boundary conditions in time.} Let us first consider 
the case when the configuration is prescribed at $j=0$ and $j=N$, 
for all $a=0,...,A$. This means that $\boldsymbol{\varphi}^0 $ and 
$ \boldsymbol{\varphi}^N$ are prescribed, therefore Hamilton's principle only yields the first equation in \eqref{stat_dynamic}, namely the DEL equations associated to $\mathsf{L}_d$. 
So, we only get the equations in \eqref{part1}. This is in complete agreement with the results obtained above in \eqref{cov_DEL} and \eqref{resulting_BC_first} via the discrete covariant variational principle when only boundary conditions in time have been assumed. This is also in complete analogy with the continuous case, where the EL equations imply the (zero-traction) boundary conditions \eqref{BC_traction}, given here in the discrete case by \eqref{resulting_BC_first}. From this discussion, we obtain that the discrete flow map
\[
F_{\mathsf{L}_d}: M^{A+1} \times M^{A+1} \rightarrow M^{A+1} \times 
M^{A+1}, \quad F_{\mathsf{L}_d}( \boldsymbol{\varphi}^{j-1}, \boldsymbol{\varphi}^j )= ( \boldsymbol{\varphi}^j, 
\boldsymbol{\varphi}^{j+1})
\]
is symplectic relative to the symplectic form $ \Omega _{\mathsf{L}_d}\in \Omega ^2 ( M^{A+1}\times M^{A+1})$: $F_{\mathsf{L}_d}^\ast \Omega _{\mathsf{L}_d}= \Omega _{\mathsf{L}_d}$. In particular, we have proven the following fact.

\begin{proposition}\label{prop1}  When boundary conditions are only imposed in time, the DEL equations for $\mathsf{L}_d$ are equivalent to the DCEL equations for $ \mathcal{L} _d $ together with the discrete zero-traction spatial boundary conditions.

Therefore, the solution $ \varphi ^j _a$, $j=0,...,N$, $a=0,...,A$ of the DCEL equations with discrete zero traction boundary conditions \eqref{resulting_BC_first} provides a symplectic-in-time discrete flow $( \boldsymbol{\varphi} ^{j-1}, \boldsymbol{\varphi} ^j ) \mapsto  ( \boldsymbol{\varphi} ^j, \boldsymbol{\varphi} ^{j+1})$ relative to the discrete symplectic form $\Omega _{\mathsf{L}_d}= \sum_{a=0}^{A-1}\Omega ^2_ {\mathcal{L}_d }(\triangle_a ^j ) =-  \sum_{a=0}^{A-1}(\Omega ^1 _ {\mathcal{L}_d }(\triangle_a ^j ) + \Omega ^3 _{ \mathcal{L}_d }(\triangle_a ^j )) $ on $M^{A+1} \times M^{A+1} $.
\end{proposition}

The equations are solved by assuming that $ \boldsymbol{\varphi} ^0$ and $ \boldsymbol{\varphi} ^1 $ are known, i.e., $ \varphi ^0 _a $ and $ \varphi ^1 _a $ for all $a=0,...,A$.

\paragraph{(B) Boundary conditions in space.} We now consider (for completeness and symmetry relative to the preceding case) the situation 
when the discrete configuration is prescribed at the boundary $a=0$ and $a=A$,
for all $j=0,...,N$, and is time independent. No boundary conditions are assumed in time. In this case, one has to incorporate these conditions in the configuration space of the discrete (dynamic) Lagrangian. Namely, we 
define the configuration space $M^{A+1}_0:=\{\boldsymbol{\varphi} \in 
M^{A+1}\mid \varphi_0 = \bar{\varphi}_0,\;\;\varphi_A=\bar{\varphi}_A\}$ 
with prescribed boundary values. This is possible since the boundary conditions at $a=0$, $a=A$ 
are assumed to be time independent. The discrete Lagrangian $\mathsf{L}_d$ is now defined as $\mathsf{L}_d: M^{A+1}_0 \times M^{A+1}_0 \rightarrow \mathbb{R}  $. The discrete Hamilton principle yields equations \eqref{stat_dynamic} (both equations) but written on $M^{A+1}_0 $ instead of $M^{A+1} $. In this case, this leads to the following slight changes in the computations of the derivatives of $\mathsf{L}_d$, namely, we have
\begin{align*}  
\left\langle D_1\mathsf{L}_d(\boldsymbol{\varphi}  ^j , \boldsymbol{\varphi}  ^{j+1}), \delta \boldsymbol{\varphi} ^j \right\rangle &=\sum_{a=1}^{A-1}( D_1 \mathcal{L} _a ^j + D_3 \mathcal{L} _{a-1} ^j ) \delta \varphi _a ^j\\
\left\langle D_2\mathsf{L}_d(\boldsymbol{\varphi}  ^j , \boldsymbol{\varphi}  ^{j+1}), \delta \boldsymbol{\varphi} ^{j+1} \right\rangle &=\sum_{a=1}^{A-1} D_2 \mathcal{L} _a ^j \delta \varphi _a ^{j+1}
\end{align*}
instead of \eqref{computations_derivatives}. In this case, equations \eqref{stat_dynamic} yield  
\begin{align*} 
&D_1 \mathcal{L} _a ^j +D_2 \mathcal{L} _a ^{j-1}+ D_3 \mathcal{L} _{a-1} ^j=0, \quad j=1,...,N-1,\;\; a=1,...,A-1\\
&D_1 \mathcal{L} _a ^0 + D_3 \mathcal{L} _{a-1} ^0 =0, \quad D_2 \mathcal{L} _a ^{N-1}=0, \quad a=1,...,A-1.
\end{align*} 
This is in agreement with equations \eqref{cov_DEL} and \eqref{resulting_BC_second} obtained earlier via the covariant discrete variational principle with boundary conditions in space only.

On $M^{A+1}_0 \times M^{A+1}_0$ the discrete one-forms $ \Theta _{\mathsf{L}_d}^\pm$ are
\begin{align*} 
\Theta _{\mathsf{L}_d}^-( \boldsymbol{\varphi} ^j , \boldsymbol{\varphi} ^{j+1})&=-\sum_{a=1}^{A-1} \Theta _{\mathcal{L}_d } ^1 (\triangle_a^j, \varphi_a^j, \varphi_a^{j+1},\varphi_{a+1}^j)+ \Theta _ {\mathcal{L} _d }^3 (\triangle_a^j, \varphi_a^j, \varphi_a^{j+1},\varphi_{a+1}^j)\\
\Theta _{\mathsf{L}_d}^+( \boldsymbol{\varphi} ^j , \boldsymbol{\varphi} ^{j+1})&= \sum_{a=1}^{A-1} \Theta _ {\mathcal{L} _d }^2 (\triangle_a^j, \varphi_a^j, \varphi_a^{j+1},\varphi_{a+1}^j).
\end{align*} 
Note the slight change in the range of summation. Relations \eqref{relation_thetas} and \eqref{relation_omegas} hold in the same way, with the same change in the summation. 

\begin{proposition}\label{prop1.5}  When boundary conditions are imposed in space, then one has to incorporate these conditions in the configuration space of the discrete Lagrangian $\mathsf{L}_d$. This yields the configuration space $M^{A+1}_0 \subset M^{A+1}$. In this case, the DEL equations for $\mathsf{L}_d$ on $M^{A+1}_0$ are equivalent to the DCEL for $\mathcal{L}_d $. The discrete covariant variational principle yields, in addition, the discrete zero momentum boundary condition in time.

Therefore, the solution $ \varphi^j_a$, $j=0,...,N$, $a=0,...,A$ of 
the DCEL equations \eqref{cov_DEL} with zero momentum boundary conditions \eqref{resulting_BC_second} provides a symplectic-in-time discrete flow $( \boldsymbol{\varphi} ^{j-1}, \boldsymbol{\varphi} ^j ) \mapsto  ( \boldsymbol{\varphi} ^j, \boldsymbol{\varphi} ^{j+1})$ relative to the discrete symplectic form $\Omega _{\mathsf{L}_d}= \sum_{a=1}^{A-1}\Omega ^2_ {\mathcal{L}_d }(\triangle_a ^j ) =-  \sum_{a=1}^{A-1}(\Omega ^1 _ {\mathcal{L}_d }(\triangle_a ^j ) + \Omega ^3 _{ \mathcal{L}_d }(\triangle_a ^j )) $ on $M_0^{A+1}\times M_0^{A+1}$.
\end{proposition}

Note that the discrete symplectic form is on $M_0^{A+1}\times M_0^{A+1}$, not on $M^{A+1}\times M^{A+1}$.

\paragraph{(C) Boundary conditions in space and time.} Of course, one 
has a similar relation between the DEL equations for $\mathsf{L}_d$ 
and DCEL equations for $\mathcal{L}_d$ in the case when both spatial 
and temporal boundary conditions are assumed. In this case, one has 
to choose, as before, the discrete configuration space $M^{A+1}_0$ 
and to consider the DEL equations for $\mathsf{L}_d$ in \eqref{stat_dynamic}, without the second boundary conditions. 
In this case, the DEL equations for $\mathsf{L}_d$ read
\[
D_1 \mathcal{L}_a^j  +D_2 \mathcal{L}_a^{j-1} + 
D_3 \mathcal{L}_{a-1} ^j=0, \quad j=1,...,N-1,\;\; a=1,...,A-1,
\]
and coincide with the DCEL equations \eqref{cov_DEL}. 
As before, we have the following result.

\begin{proposition}\label{prop2}   When boundary conditions are imposed in space and time, the DEL equations for $\mathsf{L}_d$ on $M^{A+1}_0$ are equivalent to the DCEL equations for $\mathcal{L} _d $. 

Therefore, the solution $\varphi^j _a$, $j=0,...,N$, $a=0,...,A$ of 
the DCEL equations \eqref{cov_DEL} provides a symplectic-in-time 
discrete flow $( \boldsymbol{\varphi}^{j-1}, \boldsymbol{\varphi}^j) 
\mapsto  (\boldsymbol{\varphi}^j, \boldsymbol{\varphi}^{j+1})$ 
relative to the discrete symplectic form $\Omega_{\mathsf{L}_d}= 
\sum_{a=1}^{A-1}\Omega^2_{\mathcal{L}_d}(\triangle_a^j) =
-\sum_{a=1}^{A-1}(\Omega^1_{\mathcal{L}_d}(\triangle_a^j) + 
\Omega^3_{\mathcal{L}_d}(\triangle_a^j))$ on $M^{A+1}_0  
\times M^{A+1}_0$.
\end{proposition}

In applications, it is usually assumed that $ \boldsymbol{\varphi} ^0$ and $ \boldsymbol{\varphi} ^1 $ (i.e., $ \varphi ^0 _a $ and $ \varphi ^1 _a $ for all $a=0,...,A$) are prescribed, corresponding to initial configuration and velocity (as opposed to $ \boldsymbol{\varphi} ^0 $ and $ \boldsymbol{\varphi} ^N $). 

Our discussion carries over to this case.

Note that the values at the extremities are prescribed and time independent: $ \varphi _0 ^j =\bar \varphi _0 $, $ \varphi _A ^j = 
\bar \varphi _A$ 
for all $j=0,...,N$  (or, equivalently, $ \boldsymbol{\varphi}  ^j \in M^{A+1}_0=\{ \boldsymbol{\varphi} \in M^{A+1}\mid \varphi _0 = \bar{ \varphi } _0 ,   \varphi_A=\bar{\varphi }_A\}$).

\subsubsection{Discrete spatial evolution: the discrete Lagrangian $\mathsf{N}_d$}

We consider now the converse situation to the one before, that is, 
we regard the spatial coordinate as the dynamic variable, whereas 
time is considered as a parameter. Mathematically speaking, this is 
simply a switching between the $s$- and $t$-variables.

The configuration space for the discrete "spatial-evolution" 
Lagrangian $\mathsf{N}_d$ is thus $M^{N+1}$. Using the notation 
$ \boldsymbol{\varphi}_a := (\varphi_a^0, ...,\varphi_a^N) \in M^{N+1}$, 
the discrete Lagrangian $\mathsf{N}_d: M^{N+1} \times M^{N+1} 
\rightarrow \mathbb{R}$ is defined by
\[
\mathsf{N}_d(\boldsymbol{\varphi}_a , \boldsymbol{\varphi}_{a+1}):= 
\sum_{j=0}^{N-1} \mathcal{L}_d(\triangle_a^j, \varphi_a^j, 
\varphi_a^{j+1},\varphi_{a+1}^j),
\]
so that the discrete action reads
\begin{equation}\label{discrete_action_N} 
\mathfrak{S}_d(\varphi_d) = \sum_{a=0}^{A-1} \mathsf{N}_d(\boldsymbol{\varphi}  _a  , \boldsymbol{\varphi}  _{a+1}) = \sum_{j=0}^{N-1} \sum_{a=0}^{A-1} \mathcal{L}_d(\triangle_a^j, \varphi_a^j, \varphi_a^{j+1},\varphi_{a+1}^j).
\end{equation}
In order to analyze the relation between the discrete Hamilton principles associated to $ \mathsf{L}_d$ and $\mathsf{N}_d$, we shall first assume that there are no boundary conditions, so that the discrete Hamilton principle for $\mathsf{N}_d$ yields the stationarity conditions
\begin{equation}\label{stat_spatial} 
D_1\mathsf{N}^a_d+D_2\mathsf{N}^{a-1}_d =0, \quad a=1,...,A-1,\quad\text{and}\quad D_1\mathsf{N}_d^0=0, \quad D_2\mathsf{N}_d^{A-1}=0.
\end{equation}
We compute 
\begin{align}
\label{comput_DN}
\left< D_1\mathsf{N}_d(\boldsymbol{\varphi}_a,
\boldsymbol{\varphi}_{a+1}), \delta \boldsymbol{\varphi}_a \right> & = \sum_{j=1}^{N-1}\Big[\left( D_1 \mathcal{L}_a^j +D_2\mathcal{L}_a^{j-1}\right)\delta \varphi_a^j + D_1\mathcal{L}_a^0\delta_a^0 + D_2 \mathcal{L}_a^{N-1}\delta \varphi_a^N\Big], \nonumber \\
\left< D_2 \mathsf{N}_d(\boldsymbol{\varphi}_a,
\boldsymbol{\varphi}_{a+1}),\delta\boldsymbol{\varphi}_{a+1}\right> & = \sum_{j=1}^{N-1}\left[D_3 \mathcal{L}_a^j\delta \varphi_{a+1}^j + D_3 \mathcal{L}_a^0 \delta \varphi_{a+1}^0\right].
\end{align} 
So the DEL equations for $\mathsf{N}_d$ in \eqref{stat_spatial} yield \begin{equation}\label{boundary_cond_1}
\begin{aligned}
& D_1\mathcal{L}_a^j+ D_2\mathcal{L}_a^{j-1} + D_3\mathcal{L}_{a-1}^j=0, \quad j=1,...,N-1, \ \ a= 1,..., A-1, \\
& D_1 \mathcal{L}_a^0+ D_3\mathcal{L}_{a-1}^0=0, \quad a=1,...,A-1, \\
& D_2 \mathcal{L}_a^{N-1}=0, \quad a=1,...,A-1.
\end{aligned}
\end{equation}
The boundary conditions in \eqref{stat_spatial} imply the equations
\begin{equation}\label{boundary_cond_2}
\begin{aligned}
& D_1\mathcal{L}_0^j + D_2\mathcal{L}_0^{j-1}=0, \quad j=1,...,N-1, \\
& D_1 \mathcal{L}_0^0=0, \ \ D_2\mathcal{L}_0^{N-1}= 0, \ \mathrm{and} \ \ D_3 \mathcal{L}_{A-1}^j=0, \quad j=0,..., N-1.
\end{aligned}
\end{equation}
So we recover exactly the stationarity conditions obtained from the discrete covariant Hamilton principle \eqref{discrete_variations} when no boundary condition is imposed.

Of course, \eqref{boundary_cond_1}-\eqref{boundary_cond_2} agree with the stationarity condition \eqref{discrete_variations} obtained from the discrete covariant Hamilton principle when no boundary condition is imposed. Moreover, this computation shows that the DEL for $\mathsf{N}_d$ (i.e., the first equation in \eqref{stat_spatial}) is equivalent to the DCEL equations for $ \mathcal{L} _d $ together with the discrete zero momentum boundary conditions in time \eqref{resulting_BC_second} (the second and third lines in \eqref{boundary_cond_1}).

\begin{remark}[Discrete Cartan forms]{\rm The discrete Cartan one-forms $\Theta_{\mathsf{N}_d}^\pm$ on $M^{N+1} \times M^{N+1}$ are computed as
\begin{align*}
\Theta_{\mathsf{N}_d}^-(\boldsymbol{\varphi}_a,\boldsymbol{\varphi}_{a+1}) & = -D_1\mathsf{N}_d(\boldsymbol{\varphi}_a, \boldsymbol{\varphi}_{a+1}) \mathbf{d} \boldsymbol{\varphi}_a \\
& = - \sum_{j=1}^{N-1}\left( D_1 \mathcal{L}_a^j +D_2\mathcal{L}_a^{j-1}\right)\mathbf{d}  \varphi_a^j + D_1\mathcal{L}_a^0 \mathbf{d} \varphi_a^0 + D_2 \mathcal{L}_a^{N-1}\mathbf{d} \varphi_a^N 
\\
& =- \sum_{j=0}^{N-1} \left( \Theta_{\mathcal{L}_d}^1(\triangle_a^j, \varphi_a^j, \varphi_a^{j+1}, \varphi_{a+1}^j) +  \Theta_{\mathcal{L}_d}^2(\triangle_a^j, \varphi_a^j, \varphi_a^{j+1}, \varphi_{a+1}^j) \right),\\
\Theta_{\mathsf{N}_d}^+(\boldsymbol{\varphi}_a,\boldsymbol{\varphi}_{a+1}) & = D_2 \mathsf{N}_d(\boldsymbol{\varphi}_a,\boldsymbol{\varphi}_{a+1}) \mathbf{d} \boldsymbol{\varphi}_{a+1} = \sum_{j= 0}^{N-1}D_3 \mathcal{L}_a^j  \mathbf{d}  \varphi_{a+1}^j \\
& = \sum_{j= 0}^{N-1} \Theta_{\mathcal{L}_d}^3 (\triangle_a^j, \varphi_a^j, \varphi_a^{j+1}, \varphi_{a+1}^j).
\end{align*}
The discrete Cartan 2-forms $\Omega_{\mathsf{N}_d}$ are related to 
the 2-forms $\Omega_{\mathcal{L}_d }$ by
\begin{equation*}
\Omega_{N_d} = \sum_{j=0}^{N-1} \Omega_{\mathcal{L}_d}^3(\triangle_a^j) 
= - \sum_{j=0}^{N-1}(\Omega_{\mathcal{L}_d}^1(\triangle_a^j) + \Omega_{\mathcal{L}_d}^2(\triangle_a^j)).
\end{equation*}
These formulas should be compared with those obtained in 
Remark \ref{DCF_L}.}
\end{remark}

\paragraph{(A) Boundary conditions in space.} When 
$\boldsymbol{\varphi}  _0 $ and $\boldsymbol{\varphi}_A$ 
are prescribed, we only get the first the equation in 
\eqref{stat_spatial}. These equations are equivalent to the results obtained in \eqref{cov_DEL} and \eqref{resulting_BC_second} via 
the discrete covariant variational principle when only boundary  
conditions in space have been assumed. The discrete flow map is 
now given by
\[
F_{\mathsf{N}_d}: M^{N+1}\times M^{N+1} \rightarrow M^{N+1}\times M^{N+1}, \qquad F_{N_d}(\boldsymbol{\varphi}_{a-1},\boldsymbol{\varphi}_{a}) = (\boldsymbol{\varphi}_a,\boldsymbol{\varphi}_{a+1})
\]
and is symplectic relative to the discrete symplectic form 
$\Omega_{\mathsf{N}_d}$. In complete analogy with Proposition 
\ref{prop1}, we get the following result.

\begin{proposition}\label{prop3}
When boundary conditions are imposed in space, the DEL equations for $\mathsf{N}_d$ are equivalent to the DCEL equations for $ \mathcal{L} _d $ together with the discrete zero momentum boundary condition in time.

Therefore, the solution $ \varphi ^j _a$, $j=0,...,N$, $a=0,...,A$ of the DCEL equations with discrete zero momentum boundary conditions \eqref{resulting_BC_second} provides a symplectic-in-space discrete flow $( \boldsymbol{\varphi}_{a-1}, \boldsymbol{\varphi}_a ) \mapsto  ( \boldsymbol{\varphi}_a, \boldsymbol{\varphi}_{a+1})$ relative to the discrete symplectic form $\Omega _{\mathsf{N}_d}=  \sum_{j=0}^{N-1} \Omega_{\mathcal{L}_d}^3(\triangle_a^j) = - \sum_{j=0}^{N-1}(\Omega_{\mathcal{L}_d}^1(\triangle_a^j) + \Omega_{\mathcal{L}_d}^2(\triangle_a^j)) $ on $M^{N+1} \times M^{N+1} $.
\end{proposition}

\paragraph{(B) Boundary conditions in time.} When the discrete configuration is prescribed at temporal boundaries (i.e., at $j=0$ and $j=N$, for all $a=0,..., A$), since we are working with the discrete spatial evolution, one has to include them in the discrete configuration space $M^{N+1}$, that is, we define the discrete configuration space $M^N_0=\{ \boldsymbol{\varphi} \in M^{N+1} \mid \varphi ^0 := \bar{\varphi}^0 ,   \varphi^N=\bar{\varphi }^N \}$, where 
$\bar{\varphi}^0$  and $\bar{\varphi}^N$ are given.

 This is possible, if these boundary conditions at $t=0$ and $t=T$ do not depend on the spatial index. The discrete Lagrangian is now defined as $\mathsf{N}_d: M^{N+1} _0 \times M^{N+1} _0 \rightarrow \mathbb{R}  $. The discrete Hamilton principle yields both equations in \eqref{stat_spatial} on $M^{N+1}_0$. Using \eqref{comput_DN}, with the obvious modifications due to the fact that we work on $M^{N+1}_0 \subset M^{N+1}$, we get the equations

\begin{equation*} 
\begin{aligned}
& D_1\mathcal{L}_a^j+ D_2\mathcal{L}_a^{j-1} + D_3\mathcal{L}_{a-1}^j=0, \quad j=1,...,N-1, \ \ a= 1,..., A-1, \\
& D_1 \mathcal{L}_0^j+ D_2\mathcal{L}_{0}^{j-1}=0, \quad D_3 \mathcal{L}_{A-1}^{j}=0, \quad j=1,...,N-1.
\end{aligned}
\end{equation*}
This is in agreement with the results obtained in \eqref{cov_DEL} and \eqref{resulting_BC_first} via the discrete covariant variational principle when only boundary conditions in time have been assumed.

The discrete one-forms $\Theta_{\mathsf{N}_d}^\pm $ on $M^{N+1}_0\times M^{N+1}_0$ are
\begin{align*}
\Theta_{\mathsf{N}_d}^-(\boldsymbol{\varphi}_a,\boldsymbol{\varphi}_{a+1}) & = - \sum_{j=1}^{N-1} \left( \Theta_{\mathcal{L}_d}^1(\triangle_a^j, \varphi_a^j, \varphi_a^{j+1}, \varphi_{a+1}^j) +  \Theta_{\mathcal{L}_d}^2(\triangle_a^j, \varphi_a^j, \varphi_a^{j+1}, \varphi_{a+1}^j) \right), 
\\
\Theta_{\mathsf{N}_d}^+(\boldsymbol{\varphi}_a,\boldsymbol{\varphi}_{a+1}) & =  \sum_{j= 1}^{N-1} \Theta_{\mathcal{L}_d}^3 (\triangle_a^j, \varphi_a^j, \varphi_a^{j+1}, \varphi_{a+1}^j),
\end{align*}
where we note the slight change in the range of summation.
We get the following result.
\begin{proposition}\label{prop4} When boundary conditions are imposed in time, then one has to incorporate these conditions in the configuration space of the discrete Lagrangian $\mathsf{N}_d$. This yields the configuration space $M^{N+1}_0 \subset M^{N+1}$. In this case, the DEL equations for $\mathsf{N}_d$ on $M^{N+1}_0$ are equivalent to the DCEL for $ \mathcal{L} _d $. The discrete variational principles yield, in addition, the discrete zero traction boundary condition.

Therefore, the solution $ \varphi ^j _a$, $j=0,...,N$, $a=0,...,A$ of the DCEL equations \eqref{cov_DEL} with discrete zero traction boundary condition \eqref{resulting_BC_first} provides a symplectic-in-space discrete flow $( \boldsymbol{\varphi} _{a-1}, \boldsymbol{\varphi}_a ) \mapsto  ( \boldsymbol{\varphi}_a, \boldsymbol{\varphi}_{a+1})$ relative to the discrete symplectic form $\Omega _{\mathsf{N}_d}= \sum_{j=1}^{N-1}\Omega ^3_ {\mathcal{L}_d }(\triangle_a ^j ) =-  \sum_{j=1}^{N-1}(\Omega ^1 _{ \mathcal{L}_d }(\triangle_a ^j ) + \Omega ^2 _ {\mathcal{L}_d }(\triangle_a ^j )) $ on $M^{N+1}_0  \times M^{N+1}_0$.
\end{proposition} 

Note that the discrete symplectic form is on $M_0^{A+1}\times M_0^{A+1}$, not on $M^{A+1}\times M^{A+1}$.

\paragraph{(C) Boundary conditions in time and space.} 
Of course, one has a similar relation between the DEL equations for $\mathsf{N}_d$ and the DCEL equations for $ \mathcal{L} _d $ in the case when both spatial and temporal boundary conditions are assumed. In this case, one has to choose, as before, the discrete configuration space $M^N_0$ and to consider the DEL equations for $\mathsf{N}_d$ in \eqref{stat_spatial}, without the second boundary conditions. In this case, the DEL equations for $\mathsf{N}_d$ read
\[
D_1 \mathcal{L} _a ^j  +D_2 \mathcal{L} _a ^{j-1} + D_3 \mathcal{L} _{a-1} ^j=0, \quad j=1,...,N-1,\;\; a=1,...,A-1,
\]
and coincide with the DCEL equations \eqref{cov_DEL}. 
As before, we have the following result.

\begin{proposition}\label{prop5}   When boundary conditions are imposed in space and time, the DEL equations for $\mathsf{L}_d$ on $M^{N+1}_0$ are equivalent to the DCEL equations for $ \mathcal{L} _d $.

Therefore, the solution $ \varphi ^j _a$, $j=0,...,N$, $a=0,...,A$ of the DCEL equations \eqref{cov_DEL} provides a symplectic-in-space discrete flow $( \boldsymbol{\varphi} _{a-1}, \boldsymbol{\varphi}_a ) \mapsto  ( \boldsymbol{\varphi}_a, \boldsymbol{\varphi}_{a+1})$ relative to the discrete symplectic form $\Omega _{\mathsf{N}_d}= \sum_{j=1}^{N-1}\Omega ^3_ {\mathcal{L}_d }(\triangle_a ^j ) =-  \sum_{j=1}^{N-1}(\Omega ^1 _{ \mathcal{L}_d }(\triangle_a ^j ) + \Omega ^2 _ {\mathcal{L}_d }(\triangle_a ^j )) $ on $M^{N+1}_0  \times M^{N+1}_0$.
\end{proposition}

\begin{remark}[On the discrete Lagrange-d'Alembert principles]{\rm We have seen that the discrete spacetime covariant Hamilton principle can be written as a classical discrete Hamilton principle for $\mathsf{L}_d$ or $\mathsf{N}_d$, see \eqref{discrete_action_L}, \eqref{discrete_action_N}. In the same way, in the presence of external forces, the covariant discrete Lagrange-d'Alembert principle \eqref{Discrete_CLdA} can be written as a classical discrete Lagrange-d'Alembert principle for $\mathsf{L}_d$ or $\mathsf{N}_d$, as follows
\begin{align} \label{CDLdA_beam_time} 
\delta\sum_{j=0}^{N-1}  \mathsf{L}_d&(\boldsymbol{\varphi}^j, \boldsymbol{\varphi}^{j+1})
  + \sum_{j=0}^{N-1} \left[ \mathsf{F}_{ \mathsf{L}_d}^-(\boldsymbol{\varphi}^j, \boldsymbol{\varphi}^{j+1}) \cdot \delta \boldsymbol{\varphi}^j + \mathsf{F}_{\mathsf{L}_d}^+(\boldsymbol{\varphi}^j, \boldsymbol{\varphi}^{j+1}) \cdot \delta \boldsymbol{\varphi}^{j+1} \right]   = 0,\\
\text{with}\;\;\;\; &\mathsf{F}_{\mathsf{L}_d}^-(\boldsymbol{\varphi}^j, \boldsymbol{\varphi}^{j+1}) \cdot \delta \boldsymbol{\varphi}^j =   \sum_{a=0}^{A-1 } \Big[  F ^1_{d}(j^1 \varphi _d (\triangle_a^j))\cdot \delta \varphi_a^j +  F ^3_{d}(j^1 \varphi _d (\triangle_a^j)) \cdot \delta \varphi_{a+1}^j \Big], \nonumber \\
&\mathsf{F}_{\mathsf{L}_d}^+(\boldsymbol{\varphi}^j, \boldsymbol{\varphi}^{j+1}) \cdot \delta \boldsymbol{\varphi}^{j+1} =     \sum_{a=0}^{A-1 } F ^2_{d}(j^1 \varphi _d (\triangle_a^j)) \cdot \delta \varphi_a^{j+1}, \nonumber
\end{align}    
and
\begin{align} \label{CDLdA_beam_space} 
\delta  \sum_{a=0}^{A-1} \mathsf{N}_d&(\boldsymbol{\varphi}_a, \boldsymbol{\varphi}_{a+1})
  + \sum_{a=0}^{A-1 } \left[\mathsf{F}_{\mathsf{N}_d}^-(\boldsymbol{\varphi}_a, \boldsymbol{\varphi}_{a+1}) \cdot \delta \boldsymbol{\varphi}_a + \mathsf{F}_{\mathsf{N}_d}^+(\boldsymbol{\varphi}_a, \boldsymbol{\varphi}_{a+1}) \cdot \delta \boldsymbol{\varphi}_{a+1} \right] = 0,\\
\text{with}\;\;\;\; &\mathsf{F}_{\mathsf{N}_d}^-(\boldsymbol{\varphi}_a, \boldsymbol{\varphi}_{a+1}) \cdot \delta \boldsymbol{\varphi}_a = \sum_{j=0}^{N-1} \Big[ F ^1_{d}(j^1 \varphi _d (\triangle_a^j)) \cdot \delta \varphi_a^j   + F ^2_{d}(j^1 \varphi _d (\triangle_a^j)) \cdot  \delta \varphi_{a}^{j+1}   
\Big],\nonumber  \\
&\mathsf{F}_{\mathsf{N}_d}^+(\boldsymbol{\varphi}_a, \boldsymbol{\varphi}_{a+1}) \cdot \delta \boldsymbol{\varphi}_{a+1} =    \sum_{j=0}^{N-1} F ^3_{d}(j^1 \varphi _d (\triangle_a^j)) \cdot \delta \varphi_{a+1}^{j}. \nonumber
\end{align}}
\end{remark} 

\subsection{Discrete momentum maps }\label{cov_VS_usual_momap} 

Suppose that the discrete covariant Lagrangian density $ \mathcal{L} _d : J ^1 Y _d \rightarrow \mathbb{R}  $ is invariant under the action of a Lie group $G$ on $M$. The associated discrete classical Lagrangians $\mathsf{L}_d: M^{A+1}\times M^{A+1} \rightarrow \mathbb{R}  $ and $\mathsf{N}_d: M^{N+1}\times M^{N+1} \rightarrow \mathbb{R}  $ associated to the "temporal evolution" and "spatial evolution", respectively, inherit this $G$-invariance. Indeed, both $\mathsf{L}_d$ and $\mathsf{N}_d$ are $G$-invariant under the diagonal action of $G$ on $M^{A+1}$ and $M^{N+1}$, respectively.

The associated discrete momentum maps are $\mathsf{J}^{\pm}_{\mathsf{L}_d} : M^{A+1}\times M^{A+1} \rightarrow \mathfrak{g}^*$ and $\mathsf{J}^{\pm}_{\mathsf{N}_d} : M^{N+1}\times M^{N+1} \rightarrow \mathfrak{g}^*$, given by
\begin{equation}\label{def_all_J}
\begin{aligned}
\left< \mathsf{ J}_{\mathsf{L}_d}^+( \boldsymbol{\varphi} ^j  , \boldsymbol{\varphi} ^{j+1}), \zeta  \right> & = \left\langle D _2 \mathsf{L}_d( \boldsymbol{\varphi} ^j  , \boldsymbol{\varphi} ^{j+1}), \zeta _{M^{A+1}} ( \boldsymbol{\varphi} ^{j+1}) \right\rangle  \\
\left< \mathsf{ J}_{\mathsf{L}_d}^-( \boldsymbol{\varphi} ^j  , \boldsymbol{\varphi} ^{j+1}), \zeta  \right> &= \left\langle -D _1 \mathsf{L}_d( \boldsymbol{\varphi} ^j  , \boldsymbol{\varphi} ^{j+1}), \zeta  _{M^{A+1}} ( \boldsymbol{\varphi} ^j )  \right\rangle\\
\left< \mathsf{ J}_{\mathsf{N}_d}^+( \boldsymbol{\varphi} _a   , \boldsymbol{\varphi} _{a+1}),  \zeta  \right> & = \left\langle D _2  \mathsf{N}_d( \boldsymbol{\varphi} _a   , \boldsymbol{\varphi} _{a+1}), \zeta _{M^{N+1}} ( \boldsymbol{\varphi} _{a+1}) \right\rangle  \\
\left< \mathsf{ J}_{\mathsf{N}_d}^-( \boldsymbol{\varphi} _a   , \boldsymbol{\varphi} _{a+1}), \zeta  \right> &= \left\langle -D _1 \mathsf{N}_d( \boldsymbol{\varphi}_a   , \boldsymbol{\varphi} _{a+1}), \zeta  _{M^{N+1}} ( \boldsymbol{\varphi} _a  )  \right\rangle,
\end{aligned}
\end{equation} 
for all $ \zeta \in \mathfrak{g}  $.

From  the definition of $\mathsf{L}_d$ and $\mathsf{N}_d$ in terms of $ \mathcal{L} _d $, we have the relations
\begin{equation}\label{}
\begin{aligned}
\mathsf{ J}_{\mathsf{L}_d}^+( \boldsymbol{\varphi} ^j  , \boldsymbol{\varphi} ^{j+1})  & = \sum_{a=0}^{A-1}J^2_{\mathcal{L}_d}( j ^1 \varphi _d ( \triangle _a ^j ))\\
\mathsf{ J}_{\mathsf{L}_d}^-( \boldsymbol{\varphi} ^j  , \boldsymbol{\varphi} ^{j+1})&=  -\sum_{a=0}^{A-1}J^1_{\mathcal{L}_d}( j ^1 \varphi _d ( \triangle _a ^j ))+J^3_{\mathcal{L}_d}( j ^1 \varphi _d ( \triangle _a ^j ))\\
\mathsf{ J}_{\mathsf{N}_d}^+( \boldsymbol{\varphi} _a   , \boldsymbol{\varphi} _{a+1})  & = \sum_{j=0}^{N-1}J^3_{\mathcal{L}_d}( j ^1 \varphi _d ( \triangle _a ^j ))\\
\mathsf{ J}_{\mathsf{N}_d}^-( \boldsymbol{\varphi} _a   , \boldsymbol{\varphi} _{a+1})&= -\sum_{j=0}^{N-1}J^1_{\mathcal{L}_d}( j ^1 \varphi _d ( \triangle _a ^j ))+J^2_{\mathcal{L}_d}( j ^1 \varphi _d ( \triangle _a ^j )),
\end{aligned}
\end{equation}  
between the various discrete momentum maps. The $G$-invariance of $\mathcal{L}_d$ implies \eqref{local_DNT}, which consistently implies $\mathsf{ J}_{\mathsf{L}_d}^+=\mathsf{ J}_{\mathsf{L}_d}^-$ and $\mathsf{ J}_{\mathsf{N}_d}^+=\mathsf{ J}_{\mathsf{N}_d}^-$.  

\paragraph{Covariant versus evolutionary Noether theorem.} In the next lemma, we relate the expression $\mathscr{J}_{B,C}^{K,L}$ in \eqref{DCN} with the discrete momentum maps $\mathsf{J}^\pm_{\mathsf{L}_d}$ and $\mathsf{J}^\pm_{\mathsf{N}_d}$. This follows from a direct computation.

\begin{lemma}\label{important_lemma}  When $B=0$ and $C=A-1$, or $K=0$ and $L=N-1$, we have, respectively
\begin{equation}\label{JL}
\begin{aligned}
\mathscr{J}_{0,A-1}^{K,L}( \varphi _d )=&\sum_{j=K+1}^L \left( J ^1_{\mathcal{L}_d}( j ^1 \varphi _d (\triangle _0 ^j )+J ^2_{\mathcal{L}_d}(  j ^1 \varphi _d (\triangle _0 ^{j-1}) )+ J ^3_{\mathcal{L}_d}( j ^1 \varphi _d (\triangle _{A-1} ^j )) \right)\\
&+ \mathsf{J}^+_{\mathsf{L}_d}( \boldsymbol{\varphi} ^L , \boldsymbol{\varphi} ^{L+1})- \mathsf{J}^-_{\mathsf{L}_d}( \boldsymbol{\varphi} ^K , \boldsymbol{\varphi} ^{K+1})
\end{aligned} 
\end{equation}
\begin{equation}\label{JN}   
\begin{aligned}
\mathscr{J}_{B,C}^{0,N-1}( \varphi _d )=&\sum_{a=B+1}^C \left( J ^1_{\mathcal{L}_d}( j ^1 \varphi _d (\triangle _a ^0 )+J ^2_{\mathcal{L}_d}(  j ^1 \varphi _d (\triangle _a ^{N-1}) )+ J ^3_{\mathcal{L}_d}( j ^1 \varphi _d (\triangle _{a-1} ^0 )) \right)\\
&+ \mathsf{J}^+_{\mathsf{N}_d}( \boldsymbol{\varphi} _C , \boldsymbol{\varphi} _{C+1})- \mathsf{J}^-_{\mathsf{N}_d}( \boldsymbol{\varphi} _B , \boldsymbol{\varphi} _{B+1}).
\end{aligned}
\end{equation} 
\end{lemma} 

\medskip

Recall from Theorem \ref{Cov_Noether}, that if $ \mathcal{L} _d $ is $G$-invariant and if $ \varphi _d $ verifies the DCEL equations $D_1 \mathcal{L} _a ^j + D_2 \mathcal{L}_a^{j-1}+ 
D_3 \mathcal{L}_{a-1}^j=0$, $j=1,...,N-1$, $a=1,...,A-1$, then $\mathscr{J}_{B,C}^{K,L}( \varphi _d )=0$ for all $0\leq B<C\leq A-1$, $0\leq K<L\leq N-1$. As we see from the Lemma \ref{important_lemma}, at this point, the discrete covariant Noether theorem $\mathscr{J}_{B,C}^{K,L}( \varphi _d )=0$ does not imply the discrete Noether theorem for $\mathsf{J}_{\mathsf{L}_d}$ and $\mathsf{J}_{\mathsf{N}_d}$. This is due to the fact that the DEL equations for $\mathsf{L}_d$ (or for $\mathsf{N}_d$) imply (but are not equivalent to) the DCEL equations for $ \mathcal{L} _d $. To analyze this situation further, we have to take into account the boundary conditions involved.

\paragraph{(A) Boundary condition in time.} In this case, the discrete equations are given by the DEL equations $D_1\mathsf{L}^j_d+D_2\mathsf{L}^{j-1}_d =0$, $j=1,...,N-1$. They are equivalent to the DCEL equations together with the zero traction boundary conditions
\begin{align*} 
&D_1 \mathcal{L} _a ^j + D_2 \mathcal{L}_a^{j-1}+ 
D_3 \mathcal{L}_{a-1}^j=0,\;\;j=1,...,N-1, \;\;a=1,...,A-1,\\
&D_1 \mathcal{L} _0 ^j+ D_2 \mathcal{L}_0 ^{j-1}=0,\;\; D_3 \mathcal{L} _{A-1}^j=0 ,\;\; j=1,...,N-1.
\end{align*}
The first of these equations implies $\mathscr{J}_{0,A-1}^{K,L}( \varphi _d )=0$, while the second and third equations imply that 
the first term of the right hand side of \eqref{JL} vanishes. So, we get
\[
\mathscr{J}_{0,A-1}^{K,L}( \varphi _d )= \mathsf{J}_{\mathsf{L}_d}( \boldsymbol{\varphi} ^L , \boldsymbol{\varphi} ^{L+1})- \mathsf{J}_{\mathsf{L}_d}( \boldsymbol{\varphi} ^K , \boldsymbol{\varphi} ^{K+1})=0,
\]
where we used $\mathsf{J}_{\mathsf{L}_d}^+=\mathsf{J}_{\mathsf{L}_d}^-$ because $\mathsf{L}_d$ is $G$-invariant.
This shows that the covariant discrete Noether theorem implies the discrete Noether theorem by choosing the special case $B=0$, $C=A-1$.

Recall that, when using the discrete Lagrangian $\mathsf{N}_d$, we have to restrict to the space $M_0^{N+1}$. The equations above are equivalent to $D_1\mathsf{N}^a_d+D_2\mathsf{N}^{a-1}_d =0$, $a=1,...,A-1$, $D_1\mathsf{N}_d^0=0$, and $D_2\mathsf{N}_d^{A-1}=0$. Note that in this case, the Noether theorem for the Lagrangian $\mathsf{N}_d$ does not apply since $G$ does not act on $M_0^{N+1}$. We can, nevertheless consider the expressions $\mathsf{J}^\pm_{\mathsf{N}_d}$. Using Lemma \ref{important_lemma} and the discrete covariant Noether theorem $\mathscr{J}_{B,C}^{0,N-1}( \varphi _d )= 0$, we see explicitly how the Noether theorem fails for $\mathsf{J}^\pm_{\mathsf{N}_d}$, namely,
\begin{equation}\label{non_cons_JN}
\begin{aligned} 
&\mathsf{J}^+_{\mathsf{N}_d}(\boldsymbol{\varphi}_C , \boldsymbol{\varphi}_{C+1})- \mathsf{J}^-_{\mathsf{N}_d}(\boldsymbol{\varphi}_B , \boldsymbol{\varphi}_{B+1})\\
&=- \sum_{a=B+1}^C \left(J^1_{\mathcal{L}_d}(
j^1 \varphi_d (\triangle_a^0 )+ J ^2_{\mathcal{L}_d}(
j^1 \varphi _d (\triangle_a^{N-1}) )+ J ^3_{\mathcal{L}_d}( 
j^1 \varphi_d (\triangle_{a-1}^0)) \right).
\end{aligned}
\end{equation} 

\paragraph{(B) Boundary condition in space.} The same discussion holds when the configuration is prescribed at the spatial boundary and when zero momentum boundary conditions in time are used, by exchanging the role of $\mathsf{L}_d$ and $\mathsf{N} _d $. In this case, we have
\[
\mathscr{J}_{B,C}^{0,N-1}( \varphi _d )=\mathsf{J}^+_{\mathsf{N}_d}( \boldsymbol{\varphi} _C , \boldsymbol{\varphi} _{C+1})- \mathsf{J}^-_{\mathsf{N}_d}( \boldsymbol{\varphi} _B , \boldsymbol{\varphi} _{B+1})=0,
\]
and $\mathscr{J}_{0,A-1}^{K,L}( \varphi _d )=0$ implies
\begin{equation}\label{non_cons_JL}
\begin{aligned}
&\mathsf{J}^+_{\mathsf{L}_d}( \boldsymbol{\varphi} ^L , \boldsymbol{\varphi} ^{L+1})- \mathsf{J}^-_{\mathsf{L}_d}( \boldsymbol{\varphi} ^K , \boldsymbol{\varphi} ^{K+1})\\
&=-\sum_{j=K+1}^L \left( J ^1_{\mathcal{L}_d}( j ^1 \varphi _d (\triangle _0 ^j )+J ^2_{\mathcal{L}_d}(  j ^1 \varphi _d (\triangle _0 ^{j-1}) )+ J ^3_{\mathcal{L}_d}( j ^1 \varphi _d (\triangle _{A-1} ^j )) \right). 
\end{aligned}
\end{equation} 

\paragraph{(C) Boundary condition in both space and time.} In this case, the equations are given by the DEL equations $D_1\mathsf{L}^j_d+D_2\mathsf{L}^{j-1}_d =0$, $j=1,...,N-1$ or, equivalently, $D_1\mathsf{N}^a_d+D_2\mathsf{N}^{a-1}_d =0$, $a=1,...,A-1$, defined on $M^{A+1}_0$ and $M^{N+1}_0$, respectively. They are both equivalent to $D_1 \mathcal{L} _a ^j + D_2 \mathcal{L}_a^{j-1}+ 
D_3 \mathcal{L}_{a-1}^j=0$, $j=1,...,N-1$, $a=1,...,A-1$. In this case, Noether's theorem for the Lagrangians $\mathsf{L}_d$ and $\mathsf{N}_d$ does not apply, since $G$ does not act on $M_0^{A+1}$ and $M_0^{N+1}$. However, the covariant Noether theorem does apply, so that 
$\mathscr{J}_{B,C}^{K,L}( \varphi _d )= 0$. We can also see directly
how the Noether's theorems fail for $\mathsf{L}_d$ and $\mathsf{N}_d$, namely
\begin{align*}
&\mathsf{J}^+_{\mathsf{L}_d}( \boldsymbol{\varphi} ^L , \boldsymbol{\varphi} ^{L+1})- \mathsf{J}^-_{\mathsf{L}_d}( \boldsymbol{\varphi} ^K , \boldsymbol{\varphi} ^{K+1})\\
&=-\sum_{j=K+1}^L \left( J ^1_{\mathcal{L}_d}( j ^1 \varphi _d (\triangle _0 ^j )+J ^2_{\mathcal{L}_d}(  j ^1\varphi _d (\triangle _0 ^{j-1}) )+ J ^3_{\mathcal{L}_d}( j ^1 \varphi _d (\triangle _{A-1} ^j )) \right)
\end{align*}
\begin{align*} 
&\mathsf{J}^+_{\mathsf{N}_d}( \boldsymbol{\varphi} _C , \boldsymbol{\varphi} _{C+1})- \mathsf{J}^-_{\mathsf{N}_d}( \boldsymbol{\varphi} _B , \boldsymbol{\varphi} _{B+1})\\
&=- \sum_{a=B+1}^C \left(J^1_{\mathcal{L}_d}(j ^1 \varphi_d 
(\triangle_a^0 )+J ^2_{\mathcal{L}_d}(j ^1 \varphi _d (\triangle _a ^{N-1}) )+ J ^3_{\mathcal{L}_d}( j ^1 \varphi _d (\triangle _{a-1} ^0 )) \right).
\end{align*}
The situation can be summarized as follows.

\begin{theorem}\label{summary_Noether}  Let $ \mathcal{L} _d : J ^1 Y _d \rightarrow \mathbb{R}  $ be a discrete covariant Lagrangian density and consider the associated discrete Lagrangians $\mathsf{L}_d:M^{A+1} \times M^{A+1} \rightarrow \mathbb{R}  $ and $\mathsf{N}_d:M^{N+1} \times M^{N+1} \rightarrow \mathbb{R}$. Consider a Lie group action of $G$ on $M$ and the associated discrete covariant momentum maps $J^k_{ \mathcal{L} _d }$ and discrete momentum maps $ \mathsf{J}^\pm_{\mathsf{L}_d}$, $ \mathsf{J}^\pm_{\mathsf{N}_d}$.
Suppose that the discrete covariant Lagrangian density $ \mathcal{L} _d : J ^1 Y _d \rightarrow \mathbb{R}  $ is invariant under the action of a Lie group $G$ on $M$. While the discrete covariant Noether theorem $\mathscr{J}_{B,C}^{K,L}(\varphi_d ) = 0$ $($Theorem \ref{Cov_Noether}$)$ is always verified, independently on the imposed boundary conditions, the validity of the discrete Noether theorems for $\mathsf{J}_{\mathsf{L}_d}$ and $\mathsf{J}_{\mathsf{N}_d}$ depends on the boundary conditions.

If the configuration is prescribed at the temporal extremities and zero traction boundary conditions are used, then the discrete momentum map $\mathsf{J}_{\mathsf{L}_d}$ is conserved. Conservation of $\mathsf{J}_{\mathsf{N}_d}^\pm$ does not hold in this case, as illustrated by  formula \eqref{non_cons_JN}.

If the configuration is prescribed at the spatial extremities and zero momentum boundary conditions are used, then the discrete momentum map $\mathsf{J}_{\mathsf{N}_d}$ is conserved. Conservation of $\mathsf{J}_{\mathsf{L}_d}^\pm$ does not hold in this case, as illustrated by  formula \eqref{non_cons_JL}. 
\end{theorem}
 
\section{Multisymplectic variational integrators on Lie groups}\label{mult_var_int}

In this section, we consider the particular case when the configuration 
field $\varphi$ takes values in a Lie group $G$. Completely analogous
to the continuous case treated in \S\ref{Cov_EL_LG}, the discrete equations also admit a formulation that uses the trivialization of the
tangent bundle of the Lie group. The resulting equations present
clear advantages in the discrete setting, since one can take advantage
of the vector space structure of the Lie algebra via the use of a 
time difference map.

In \S\ref{DCELELG}, we present the discrete covariant Hamilton 
principle, the discrete Legendre transform, the discrete Cartan 
forms, the DCEL equations, the discrete covariant momentum maps, 
and the discrete covariant Noether theorem, in their trivialized 
formulation. In \S\ref{subsec_TSDE}, we quickly describe the 
symplectic properties of the time-evolutionary and space-evolutionary 
discrete flows in the trivialized form on Lie groups, following 
the results obtained in \S\ref{symplectic_properties} and 
\S\ref{cov_VS_usual_momap}.

\subsection{Discrete covariant Euler-Lagrange equations on Lie groups}\label{DCELELG} 

Let us now consider the case when the fiber $M=G$ is a Lie group. We use the notation $\varphi _a ^j =g _a ^j $. Recall that the discrete version of the first jet bundle $J ^1 (X \times G)$ is $J^1(X _d \times G)=X ^\triangle _d \times G ^3 $. Note also that we have the isomorphism
\[
X^\triangle _d \times G ^3 \ni ( \triangle_a^j, g _a ^j , g _a^{j+1}, g _{a+1}^{j}) \mapsto \left( \triangle_a^j,  g _a ^j , (g _a ^j) ^{-1} g _a^{j+1}, (g _a ^i) ^{-1}  g_{a+1}^{j}\right)\in X^\triangle _d \times G ^3 .
\]

In order to discretize the relations $ \xi = g ^{-1} \dot g $ and $ \eta= g ^{-1} g'$, we shall fix a local diffeomorphism $ \tau : \mathfrak{g}  \rightarrow G $ in a neighborhood of the identity, such that $ \tau (0)=e$. Examples for $ \tau $ are provided by the exponential map or the Cayley transform. The approach will involve the right trivialized 
derivative ${\rm d}^R\tau $ of $ \tau $ defined by 
\begin{equation}
\label{r_log_der_def}
{\rm d}^R\tau _\xi : \mathfrak{g}  \rightarrow \mathfrak{g}  , \quad {\rm d}^R\tau_\xi( \eta ):= \left( T_ \xi \tau ( \eta ) \right) \tau ( \xi ) ^{-1},
\end{equation}
where $T_ \xi \tau : \mathfrak{g}  \rightarrow T_{ \tau (\xi) } G$ is the derivative of $ \tau $.
The right trivialized derivative of $ \tau ^{-1} $ is defined by
\[
{\rm d}^R \tau ^{-1} _\xi : \mathfrak{g}  \rightarrow \mathfrak{g}  , \quad {\rm d}^R \tau ^{-1} _\xi ( \eta ):= T_g \tau ^{-1} ( \eta g),
\]
where $g:= \tau (\xi)$. It is readily checked that $ {\rm d }^R\tau ^{-1}_ \xi = ({\rm d} ^R \tau_ \xi ) ^{-1} $.

Using the local diffeomorphism $ \tau $, the relations $ \xi = g ^{-1} \dot g $ and $ \eta= g ^{-1} g'$ are discretized as
\begin{equation}\label{delta_xi_discrete} 
\begin{aligned}
\xi^j_a: &=  \tau^{-1}\left((g_a^j )^{-1}g_a^{j+1}\right)/\Delta t\in\mathfrak{g}, \\
\eta_a^j : &=  \tau^{-1}\left((g^j_a)^{-1}g^j_{a+1} \right)/\Delta s\in\mathfrak{g}.
\end{aligned}
\end{equation} 
From these definitions, we can define the discrete Lagrangian $ \bar{\mathcal{L}}_d :X_d^\triangle \times G \times \mathfrak{g} \times \mathfrak{g} \rightarrow \mathbb{R}  $ by
\[
\bar{ \mathcal{L} }_d(\triangle_a^j, g _a ^j, \xi_a^j , \eta_a^j):=\mathcal{L}_d ( \triangle_a^j, g _a ^j , g _a^{j+1}, g _{a+1}^{j}).
\]
Note that $X_d^\triangle \times G \times \mathfrak{g} \times \mathfrak{g} $ is thought of as the discretization of the trivialized first jet bundle $(T^\ast X\otimes \mathfrak{g}) \times G$, see \eqref{trivialization}, and the discrete Lagrangian $\bar{ \mathcal{L} }_d$ is the discretization of the trivialized Lagrangian $\bar { \mathcal{L} }=\bar{ \mathcal{L} }(g, \xi , \eta )$ defined in \S\ref{Cov_EL_LG}.

We have the following relations between the partial derivatives of $ \mathcal{L} _d $ and $\bar{ \mathcal{L} } _d $.
\begin{align}\label{link_PD} 
(g_a^j)^{-1}D_1 \mathcal{L} _a^j & =  (g_a^j)^{-1}D_g \bar{\mathcal{L}}_a^j -  \frac{1}{\Delta t} \mu_a^j -  \frac{1}{\Delta s} \lambda_a^j,\nonumber \\
(g_a^{j+1})^{-1}D_2 \mathcal{L}_a^j & =  \frac{1}{\Delta t}  \operatorname{Ad}_{\tau(\Delta t\xi_a^{j})}^*\mu_a^{j}, \\
(g_{a+1}^j)^{-1}D_3\mathcal{L} _a^j & =   \frac{1}{\Delta s} \operatorname{Ad}_{\tau(\Delta s \eta_{a}^{j})}^*\lambda_{a}^{j} .\nonumber
\end{align}

\subsubsection{Discrete covariant Hamilton's principle}\label{subsec_DCHP}

The discrete covariant Hamilton's principle reads
\begin{align}\label{DEL_multisymplectic} 
  \begin{split}
\delta \bar{\mathfrak{S}}_d(g_d) = & \, \delta\sum_{j=0}^{N-1} \sum_{a=0}^{A-1} \bar{\mathcal{L}}_d( \triangle_a^j, g _a ^j , \xi _a ^j , \eta _a ^j )  = 0.
\end{split}
\end{align}
Using the definitions \eqref{delta_xi_discrete}, we obtain the variations 
\begin{equation} \label{xi_eta_variations}
\begin{aligned}
\delta \xi^j_a & = {\rm d}^R\tau^{-1}_{\Delta t\xi^j_a } \left(-\zeta ^j _a +
\operatorname{Ad}_{\tau(\Delta t\xi^j_a )}\zeta ^{j+1}_a \right)/\Delta t, \\
\delta \eta ^j_a & = {\rm d}^R\tau^{-1} _{\Delta s\eta ^j_a}\left(-\zeta ^j _a +
\operatorname{Ad}_{\tau(\Delta s\eta ^j_a )}\zeta ^{j}_{a+1} \right)/\Delta s,
\end{aligned}
\end{equation}
where we defined $\zeta  ^j _a:= ( g ^j_a  ) ^{-1} \delta g _a ^j$ and we used \eqref{r_log_der_def}.

For simplicity, we will use the notation $\bar{\mathcal{L}}_a^j:=\bar{\mathcal{L}}_d( \triangle_a^j, g _a ^j , \xi _a ^j , \eta _a ^j )$. Defining the discrete momenta
\begin{align*}
&\mu_a^j:=\left( \operatorname{d}^R \tau^{-1} _{\Delta t\xi _a^j}\right) ^\ast D_\xi \bar{\mathcal{L}} _a^j ,
 \qquad \lambda_a^j:=\left( \operatorname{d}^R \tau^{-1} _{\Delta s \eta_a^j}\right)^*D_\eta \bar{\mathcal{L}}_a^j
\end{align*}
and applying the covariant discrete Hamilton principle we get
{\scriptsize \begin{equation}\label{DEL_principle}
\begin{aligned}
&\delta \bar{\mathfrak{S}}_d(g_d)  = \sum_{j=0}^{N-1} \sum_{a=0}^{A-1} \left( D_g \bar{\mathcal{L}}_a^j\cdot \delta g_a^j + D_\xi \bar{\mathcal{L}}_a^j\cdot \delta \xi_a^j + D_\eta \bar{\mathcal{L}}_a^j\cdot \delta \eta_a^j\right)  \\
& =  \sum_{j=0}^{N-1} \sum_{a=0}^{A-1} \left[ (g_a^j)^{-1}D_g \bar{\mathcal{L}}_a^j -  \frac{1}{\Delta t} \mu_a^j -  \frac{1}{\Delta s} \lambda_a^j \right] \cdot \zeta_a^j +  \frac{1}{\Delta t}  \operatorname{Ad}_{\tau(\Delta t\xi_a^{j})}^*\mu_a^{j} \cdot \zeta_a^{j+1}+  \frac{1}{\Delta s} \operatorname{Ad}_{\tau(\Delta s \eta_{a}^{j})}^*\lambda_{a}^{j} \cdot \zeta_{a+1}^j \\
& = \sum_{j=1}^{N-1} \sum_{a=1}^{A-1} \left[ (g_a^j)^{-1}D_{g} \bar{\mathcal{L}}_a^j + \frac{1}{\Delta t} \left( \operatorname{Ad}_{\tau(\Delta t\xi_a^{j-1})}^*\mu_a^{j-1} -\mu_a^j \right) + \frac{1}{\Delta s} 
\left(\operatorname{Ad}_{\tau(\Delta s \eta_{a-1}^{j})}^*\lambda_{a-1}^{j} - \lambda_{a}^j \right) \right]\cdot \zeta_a^j \\
& + \sum_{j=1}^{N-1} \left\{ \left[ (g_0^j)^{-1}D_{g} \bar{\mathcal{L}}_0^j + \frac{1}{\Delta t} \left( \operatorname{Ad}_{\tau(\Delta t\xi_0^{j-1})}^*\mu_0^{j-1} -\mu_0^j \right) - \frac{1}{\Delta s} \lambda_{0}^j  \right]\cdot  \zeta_0^j + \frac{1}{\Delta s} 
 \operatorname{Ad}_{\tau(\Delta s \eta_{A-1}^{j})}^*\lambda_{A-1}^{j} \cdot \zeta_A^j \right\} \\
 & +\sum_{a=1}^{A-1}  \left\{ \left[ (g_a^0)^{-1}D_{g} \bar{\mathcal{L}}_a^0 - \frac{1}{\Delta t} \mu_a^0  + \frac{1}{\Delta s} 
\left(\operatorname{Ad}_{\tau(\Delta s \eta_{a-1}^{0})}^*\lambda_{a-1}^{0} - \lambda_{a}^0 \right) \right]\cdot  \zeta_a^0 + \frac{1}{\Delta t}  \operatorname{Ad}_{\tau(\Delta t\xi_a^{N-1})}^*\mu_a^{N-1} \cdot  \zeta_a^N \right\} \\
 & +  \left[ (g_0^0)^{-1}D_{g} \bar{\mathcal{L}}_0^0 - \frac{1}{\Delta t} \mu_0^0 - \frac{1}{\Delta s} \lambda_{0}^0 \right] \cdot \zeta_0^0 + \left[ \frac{1}{\Delta t}  \operatorname{Ad}_{\tau(\Delta t\xi_0^{N-1})}^*\mu_0^{N-1} \right] \cdot \zeta_0^N  \\
 & + \left[ \frac{1}{\Delta s} 
 \operatorname{Ad}_{\tau(\Delta s \eta_{A-1}^{0})}^*\lambda_{A-1}^{0} \right] \cdot  \zeta_A^0.
\end{aligned}
\end{equation}}
This can be also obtained directly from \eqref{discrete_variations} by using \eqref{link_PD} 

For later use, we now list the stationarity conditions obtained in the case there is no boundary conditions imposed on the variations
\begin{align}
&\begin{aligned}
&\frac{1}{\Delta t} \left(\mu_a^j-\operatorname{Ad}_{\tau(\Delta t\xi_a^{j-1})}^*\mu_a^{j-1}\right) + \frac{1}{\Delta s} 
\left(\lambda_{a}^j-\operatorname{Ad}_{\tau(\Delta s \eta_{a-1}^{j})}^*\lambda_{a-1}^{j}\right) = (g_a^j)^{-1}D_{g} \bar{\mathcal{L}}_a^j,\\
& \qquad \qquad \text{for all $j=1,...,N-1,\;\; a=1,...,A-1$},\qquad 
\end{aligned}\label{1} \\
&\begin{aligned}
&\frac{1}{\Delta t} \left( \mu_0^j-\operatorname{Ad}_{\tau(\Delta t\xi_0^{j-1})}^*\mu_0^{j-1} \right) + \frac{1}{\Delta s} \lambda_{0}^j  =(g_0^j)^{-1}D_{g} \bar{\mathcal{L}}_0^j, \\
&\qquad \qquad \text{and}\quad \frac{1}{\Delta s} \operatorname{Ad}_{\tau(\Delta s \eta_{A-1}^{j})}^*\lambda_{A-1}^{j}  = 0,\qquad \text{for all $j=1,...,N-1$},\qquad 
\end{aligned} \label{2}\\
&\begin{aligned}
&\frac{1}{\Delta t} \mu_a^0  + \frac{1}{\Delta s} 
\left(\lambda_{a}^0-\operatorname{Ad}_{\tau(\Delta s \eta_{a-1}^{0})}^*\lambda_{a-1}^{0}  \right) =(g_a^0)^{-1}D_{g} \bar{\mathcal{L}}_a^0,\\
&\qquad  \qquad\text{and}\quad \frac{1}{\Delta t}  \operatorname{Ad}_{\tau(\Delta t\xi_a^{N-1})}^*\mu_a^{N-1} =0,\qquad \text{for all $a=1,...,A-1$},
\end{aligned}\label{3}\\
&\begin{aligned}
&\frac{1}{\Delta t} \mu_0^0 + \frac{1}{\Delta s} \lambda_{0}^0=(g_0^0)^{-1}D_{g} \bar{\mathcal{L}}_0^0,\\
&\qquad \frac{1}{\Delta t}  \operatorname{Ad}_{\tau(\Delta t\xi_0^{N-1})}^*\mu_0^{N-1}=0\quad\text{and}\quad  \frac{1}{\Delta s} \operatorname{Ad}_{\tau(\Delta s \eta_{A-1}^{0})}^*\lambda_{A-1}^{0}=0\qquad.
\end{aligned}\label{4}
\end{align}

Equations \eqref{1} will be referred to as the \textit{Lie group DCEL equations}.

\begin{remark}[Discrete Cartan forms]{\rm
In terms of the trivialized discrete Lagrangian $ \bar{ \mathcal{L} }_d $, the discrete Cartan forms \eqref{DCFs} are computed to be 
\begin{equation}\label{DCartan_1form}
\begin{aligned}
\theta_{\bar{\mathcal{L}}_d}^1( \triangle_a^j, g _a ^j , \xi _a ^j , \eta _a ^j ) & = \left\langle  (g_a^j)^{-1}D_g \bar{\mathcal{L}}_a^j -  \frac{1}{\Delta t} \mu_a^j -  \frac{1}{\Delta s} \lambda_a^j ,    (g_a^j)^{-1} d  g_a^j\right\rangle , \\
\theta_{\bar{\mathcal{L}}_d}^2( \triangle_a^j, g _a ^j , \xi _a ^j , \eta _a ^j ) & = \left\langle \frac{1}{\Delta t}  \operatorname{Ad}_{\tau(\Delta t\xi_a^{j})}^*\mu_a^{j} , (g_a^{j+1})^{-1}d g_a^{j+1}\right\rangle  \\
&= \left\langle \frac{1}{\Delta t}g ^j _a \mu ^j _a , d g ^j _a \right\rangle+ \left\langle  D_\xi \bar{\mathcal{L}} _a^j , d \xi _a ^j \right\rangle  \\
\theta_{\bar{\mathcal{L}}_d}^3( \triangle_a^j, g _a ^j , \xi _a ^j , \eta _a ^j ) & = \left\langle  \frac{1}{\Delta s} \operatorname{Ad}_{\tau(\Delta s \eta_{a}^{j})}^*\lambda_{a}^{j} ,(g_{a+1}^j)^{-1}d g_{a+1}^j\right\rangle \\
&= \left\langle \frac{1}{\Delta t}g ^j _a \lambda  ^j _a , d g ^j _a \right\rangle+ \left\langle  D_\eta  \bar{\mathcal{L}} _a^j , d \eta  _a ^j \right\rangle.
\end{aligned}
\end{equation}
We note the relation
\begin{equation}\label{relation_discrete_forms}
\theta ^k _{\bar{ \mathcal{L}} _d }=( \phi _\tau ) _\ast \Theta _{ \mathcal{L} _d }^k , \quad k=1,2,3,
\end{equation} 
where $ \phi _\tau :X_d^\triangle \times G \times G \times G \rightarrow X_d^\triangle \times G \times \mathfrak{g}  \times \mathfrak{g}  $ is the local diffeomorphism defined by $ \phi _\tau ( \triangle _a ^i ,g _a ^i , g^{i+1} _a , g ^i _{a+1})=(  \triangle _a ^i,g ^i _a , \xi ^i _a , \zeta ^i _a )$ and $ \Theta ^k _{ \mathcal{L} _d }$ are the discrete one-forms defined in \eqref{DCFs}. From the relations \eqref{relation_discrete_forms} and $\bar{\mathcal{L}} _d \circ \phi _\tau = \mathcal{L} _d $ and the formula \eqref{one_form_relationship}, we get
\begin{equation} \label{one_form_relationship_LA}
\mathbf{ d} \bar{\mathcal{L}}_d = \theta_{\bar{\mathcal{L}}_d }^{1} + \theta_{\bar{\mathcal{L}}_d }^{2} + \theta_{\bar{\mathcal{L}}_d }^{3}.
\end{equation}}
\end{remark}

\medskip

Given a vector field $V$ on $X_d^\triangle \times G$ and its first jet extension $j^1V$ on the discrete jet bundle $X_d^\triangle \times G \times G \times G$, we define the vector field $\overline{ j ^1 V}$   
induced on $X_d^\triangle \times G \times \mathfrak{g}   \times \mathfrak{g}  $ by $ \phi _\tau $. Similarly, given a discrete section $ g _d $ and its first jet extension $ j ^1 g _d :X_d^\triangle \rightarrow X_d^\triangle \times G \times G \times G$ we define $\overline{ j ^1 g _d }:= \phi _\tau \circ j ^1 g_d $. With these notations, we can write the formulas
\begin{align*}
\left\langle (g_a^j)^{-1}D_g \bar{\mathcal{L}}_a^j -  \frac{1}{\Delta t} \mu_a^j -  \frac{1}{\Delta s} \lambda_a^j , ( g _a ^j ) ^{-1} \delta g _a ^j \right\rangle  & = \left[ ( \overline{ j^1g_d})^* \mathbf{i}_{\overline{ j^1V}\,} \theta_{\bar{\mathcal{L}}_d}^1 \right](\triangle_a^j),
\\
\left\langle \frac{1}{\Delta t}  \operatorname{Ad}_{\tau(\Delta t\xi_a^{j})}^*\mu_a^{j} , ( g _a ^{j+1}) ^{-1} \delta  g _a ^{j+1} \right\rangle& = \left[( \overline{ j^1g_d})^* \mathbf{i}_{\overline{ j^1V}\,}  \theta_{\bar{\mathcal{L}}_d}^2 \right](\triangle_a^j), 
\\
\left\langle  \frac{1}{\Delta s} \operatorname{Ad}_{\tau(\Delta s \eta_{a}^{j})}^*\lambda_{a}^{j} , (g_{a+1} ^j ) ^{-1} \delta g_{a+1} ^j \right\rangle  & = \left[ ( \overline{ j^1g_d})^* \mathbf{i}_{\overline{ j^1V}\,}  \theta_{\bar{\mathcal{L}}_d}^3 \right](\triangle_a^j).
\end{align*}
from which we deduce, as in \eqref{reformulated_discrete_action}, that \eqref{DEL_principle} can be written as
\begin{equation}\label{DEL_principle_Cartan_forms}
\begin{aligned}
\delta \bar{\mathfrak{S}}_d(g_d) =&\sum_{j=1}^{N-1} \sum_{a=1}^{A-1} \left\langle  (g_a^j)^{-1}D_{g} \bar{\mathcal{L}}_a^j + \frac{1}{\Delta t} \left( \operatorname{Ad}_{\tau(\Delta t\xi_a^{j-1})}^*\mu_a^{j-1} -\mu_a^j \right) \right. \\
& \left. \hspace{3.5cm} + \frac{1}{\Delta s} 
\left(\operatorname{Ad}_{\tau(\Delta s \eta_{a-1}^{j})}^*\lambda_{a-1}^{j} - \lambda_{a}^j \right) , ( g _a ^j) ^{-1} \delta  g _a ^j \right\rangle \\
 &+ \sum_{\{\triangle \in X_d^ \triangle\mid \triangle \cap \partial U \neq \varnothing\}} \left( \sum_{k \in \{1,2,3\}; \triangle^{(k)}\in \partial U} \left[ (\overline{j^1g_d} )^* (\mathbf{i}_{\overline{j^1V}\,} \theta_{\bar{\mathcal{L}}_d}^k) \right](\triangle ) \right).
\end{aligned} 
\end{equation} 

\color{black}

\paragraph{(A) Spacetime boundary conditions.} When the values of the discrete configuration $ g _a ^j $ are prescribed at the spacetime boundary, then the covariant Hamilton principle only yield the Lie-group DCEL equations \eqref{1}.

\paragraph{(B) Boundary conditions in time.} When the values of the discrete configuration $ g _a ^j $ are prescribed for $j=0$ and $j=N$,
then the covariant Hamilton principle yields the Lie-group DCEL equations  \eqref{1} together with the \textit{discrete zero traction boundary conditions \eqref{2}}.

\paragraph{(C) Boundary conditions in space.} When the values of the discrete configuration $g_a^j $ are prescribed for $a=0$ and $a=A$, then the covariant Hamilton principle yields the Lie-group DCEL equations  \eqref{1} together with the \textit{discrete zero momentum boundary conditions \eqref{3}}.

\color{black}
\subsubsection{Discrete Legendre transforms}

The discrete covariant Legendre transforms have been defined in \eqref{DCLegendre_transf}. Their expression in terms of the discrete trivialized Lagrangian $\bar{\mathcal{L}}_d$ are $\mathbb{F}^k\bar{\mathcal{L}}_d : X_d^\triangle \times G \times \mathfrak{g} \times \mathfrak{g} \rightarrow G \times \mathfrak{g}^*$ $k=1,2,3$,
given by
\begin{equation}\label{DLegendre}
\begin{aligned}
\mathbb{F}^1\bar{\mathcal{L}}_d ( \triangle_a^j,g _a ^j , \xi _a ^j , \eta _a ^j ) & = \left(g_a^j, \, (g_a^j)^{-1}D_g \bar{\mathcal{L}}_a^j -  \frac{1}{\Delta t} \mu_a^j -  \frac{1}{\Delta s} \lambda_a^j \right), \\
\mathbb{F}^2\bar{\mathcal{L}}_d (\triangle_a^j, g _a ^j , \xi _a ^j , \eta _a ^j ) & = \left( g_a^{j+1},\, \frac{1}{\Delta t}  \operatorname{Ad}_{\tau(\Delta t\xi_a^{j})}^*\mu_a^{j} \right), \\
\mathbb{F}^3\bar{\mathcal{L}}_d (\triangle_a^j, g _a ^j , \xi _a ^j , \eta _a ^j ) & = \left( g_{a+1}^j,\, \frac{1}{\Delta s} \operatorname{Ad}_{\tau(\Delta s \eta_{a}^{j})}^*\lambda_{a}^{j} \right).
\end{aligned}
\end{equation}
Note that it is related to the  discrete Legendre transforms $\mathbb{F}^k  \mathcal{L} _d$ via the formula
\[
\rho _L \circ  \mathbb{F}^k  \mathcal{L} _d= \mathbb{F}^k\bar{\mathcal{L}}_d \circ \phi _\tau,
\]
where $ \rho _L: T^*G \rightarrow G \times \mathfrak{g}  ^\ast $, $ \rho _L( \alpha _g ):=( g , g ^{-1} \alpha _g )$. Also, as before,
the DCEL equations \eqref{1} can be written as
\[
\mathbb{F}^1\bar{\mathcal{L}}_d (\triangle_a^j, g _a ^j , \xi _a ^j , \eta _a ^j ) + \mathbb{F}^2\bar{\mathcal{L}}_d (\triangle_a^j, g _a ^{j-1} , \xi _a ^{j-1} , \eta _a ^{j-1} )+\mathbb{F}^3\bar{\mathcal{L}}_d (\triangle_a^j, g _{a-1} ^j , \xi _{a-1} ^j , \eta _{a-1}^j )=0.
\]

\subsubsection{Discrete momentum maps}

We now consider symmetries given by a subgroup $H$ of the Lie group fiber $G$. We assume that $H$ acts on the left by translation, i.e. $ \Phi :H \times G \rightarrow G$, $ \Phi _h(g)=hg$. The infinitesimal generator associated to $ \zeta \in \mathfrak{h}  $ is $\zeta _G(g)= \zeta g$. Using the formulas \eqref{discrete_momap} adapted to this special case and written in terms of the trivialized discrete Lagrangian $\bar{ \mathcal{L} } _d $, we get the discrete momentum maps $J^k _{ \bar{\mathcal{L}} _d }:X_d^\triangle \times G \times \mathfrak{g} \times \mathfrak{g} \rightarrow \mathfrak{h}^*$, $k=1,2,3$,
\begin{equation}\label{discr_momap_LA}
\begin{aligned}
J^1_{\bar{\mathcal{L}}_d}(\triangle_a^j, g _a ^j , \xi _a ^j , \eta _a ^j )&= i ^\ast \mathrm{Ad}_{(g_a^j)^{-1}}^* \left( (g_a^j)^{-1}D_g \bar{\mathcal{L}}_a^j -  \frac{1}{\Delta t} \mu_a^j -  \frac{1}{\Delta s} \lambda_a^j \right), \\
J^2_{\bar{\mathcal{L}}_d}(\triangle_a^j, g _a ^j , \xi _a ^j , \eta _a ^j )&= i ^\ast\frac{1}{\Delta t}  \operatorname{Ad}_{(g_a^j)^{-1}}^*\mu_a^{j} , \\
J^3_{\bar{\mathcal{L}}_d}(\triangle_a^j, g _a ^j , \xi _a ^j , \eta _a ^j )&=  i ^\ast \frac{1}{\Delta s} \operatorname{Ad}_{(g_a^j)^{-1}}^*\lambda_{a}^{j} ,
\end{aligned}
\end{equation} 
where $i ^\ast : \mathfrak{g}  ^\ast \rightarrow \mathfrak{h}  ^\ast $ denotes the dual map to the Lie algebra inclusion $i: \mathfrak{h}  \rightarrow \mathfrak{g}  $. We have the formula
\begin{equation}\label{formulas_J} 
\left\langle J ^k _{ \bar{\mathcal{L}}_d  }, \xi \right\rangle =\mathbf{i} _{\overline{ \xi} _{ J ^1 (X _d \times G)}} \theta ^k _{ \bar{\mathcal{L}} _d  },\quad \xi  \in \mathfrak{h}\quad\text{and}\quad J^k_{\bar {\mathcal{L} }_d} \circ \phi _\tau = J^k_{\mathcal{L} _d},
\end{equation} 
where $\xi _{ J ^1 (X _d \times G)}$ is the infinitesimal generator of the $H$-action induced on $J ^1 Y _d = X _d ^\triangle \times G \times G \times G$ and $ \overline{ \xi} _{ J ^1 (X _d\times G) }$ is the vector field induced on $X _d ^\triangle \times G \times \mathfrak{g}   \times \mathfrak{g}  $.

Let us assume that $\bar{\mathcal{L}} _d $ is $H$-invariant, that is $\bar{ \mathcal{L} } _d (\triangle _a ^j , h g ^j _a , \xi _a ^j , \eta _a ^j )=\bar{ \mathcal{L} } _d (\triangle _a ^j ,  g ^j _a , \xi _a ^j , \eta _a ^j )$, for all $ h \in H$.
This implies the infinitesimal $H$-invariance $\left\langle D_g \bar{\mathcal{L}}_a ^j, \zeta g _a ^j  \right\rangle =0$, for all $ \zeta \in\mathfrak{h}  $,  i.e. $ i ^\ast  D_g \bar{\mathcal{L}}_a ^j(g _a ^j ) ^{-1} =0$. From the expressions \eqref{discr_momap_LA}, this can be equivalently written as
\[
\left( J^1_{\bar{\mathcal{L}}_d} + J^2_{\bar{\mathcal{L}}_d}+ J^3_{\bar{\mathcal{L}}_d}\right) (\triangle_a^j, g _a ^j , \xi _a ^j , \eta _a ^j )=0,
\]
which is the statement of the local discrete Noether theorem.
We now state its global version.

\begin{theorem}\label{theorem_DCN_LA}  Suppose that the discrete Lagrangian $\bar{\mathcal{L}}_d: X_d^\triangle \times G \times \mathfrak{g}\times \mathfrak{g} \rightarrow \mathbb{R}$ is invariant under the left action of the Lie group $H$ on $G$. Suppose that $ g _d   $ is a solution of the DCEL equations for $ \bar{\mathcal{L} }_d $. Then, for all $0\leq B<C\leq A-1$, $0\leq K<L\leq N-1$, we have the conservation law
\begin{equation}\label{DCN_G}  
\bar{\mathscr{J}}_{B,C}^{K,L}(g _d  ) = 0,
\end{equation}
where $\bar{\mathscr{J}}_{B,C}^{K,L}$ is given by \eqref{Def_DN}, with $J ^k_{\mathcal{L}_d}( j ^1 \varphi _d (\triangle _a ^K )$ replaced by $J ^k_{\bar{\mathcal{L}}_d}( j ^1 g_d (\triangle _a ^K )$.
\end{theorem}
\textbf{Proof}. From the $H$-invariance of $\bar{\mathcal{L}} _d$ we have $\mathfrak{S}_d(h\cdot g_d) = \mathfrak{S}_d(g_d)$, so, the derivative of this this expression with respect to $h$ vanishes. Using this fact in \eqref{DEL_principle_Cartan_forms} together the fact that $ g _d $ is a solution of the DCEL equations, we get
\begin{align*}
0= \sum_{\{\triangle \in X_d^ \triangle\mid \triangle \cap \partial X_d \neq \varnothing\}} \left( \sum_{k \in \{1,2,3\}; \triangle^{(k)}\in \partial X_d} \left[  (\overline{j^1g_d} )^* (\mathbf{i}_{\bar{\xi }_{J^1(X_d\times G)}} \theta_{\bar{\mathcal{L}}_d}^k) \right](\triangle ) \right),
\end{align*}
for all $ \xi \in  \mathfrak{h}  $. More generally this can be done for a rectangular subdomain $U$ as in Theorem \ref{Cov_Noether}. The global Noether theorem follows from the first formula in \eqref{formulas_J}.
$\qquad\blacksquare$

\subsubsection{Covariant Euler-Lagrange equations with forces}\label{forces}

Recall that when forces are present, one has to use the principle \eqref{Discrete_CLdA}. Given the discrete forces
{\footnotesize \[
F ^1 _d (\triangle _a ^j , g _a ^j , g^{j+1}_a, g^j_{a+1})\in T^*_{ g _a ^j }G, \quad F ^2_d  (\triangle _a ^j , g _a ^j , g^{j+1}_a, g^j_{a+1})\in T^*_{ g _a ^{j+1} }G, \quad F ^3_d  (\triangle _a ^j , g _a ^j , g^{j+1}_a, g^j_{a+1})\in T^*_{ g _{a+1} ^j }G,
\]}
we define their trivialization
\begin{align*} 
\bar F ^1_d  ( \triangle _a ^j ,g _a ^j , \xi _a ^j , \eta _a ^j ):&=(g _a ^j) ^{-1} F ^1_d  (\triangle _a ^j , g _a ^j , g^{j+1}_a, g^j_{a+1})\\
\bar F ^2_d  ( \triangle _a ^j ,g _a ^j  , \xi _a ^j , \eta _a ^j ):&= (g _a ^{j+1}) ^{-1} F ^2 _d  (\triangle _a ^j , g _a ^j , g^{j+1}_a, g^j_{a+1})\\
\bar F ^3_d  ( \triangle _a ^j ,g _a ^j  , \xi _a ^j , \eta _a ^j ):&= (g _{a+1} ^j) ^{-1} F ^3 _d  (\triangle _a ^j , g _a ^j , g^{j+1}_a, g^j_{a+1}).
\end{align*}
Using these definitions, the variational principle \eqref{Discrete_CLdA} reads
\begin{equation}\label{Discrete_CLdA_LG}
\begin{aligned} 
&\delta\sum_{j=0}^{N-1} \sum_{a=0}^{A-1}\bar{\mathcal{L}} _d(\triangle_a^j, g_a^j, \xi_a^j, \eta_a^j)
  + \sum_{j=0}^{N-1} \sum_{a=0}^{A-1}\Big[ \bar F ^1_d  ( \triangle _a ^j ,g _a ^j , \xi _a ^j , \eta _a ^j ) \zeta_a^j \\
& \qquad\qquad\qquad   \left. 
+ \bar F ^2_d  ( \triangle _a ^j ,g _a ^{j} , \xi _a ^j , \eta _a ^j )\zeta_{a}^{j+1}  
+ \bar F ^3_d  ( \triangle _a ^j ,g _{a} ^j , \xi _a ^j , \eta _a ^j )\zeta_{a+1}^{j}\right] = 0
\end{aligned}
\end{equation}    
and yield the forced DCEL equations
\begin{align*}
& (g_a^j)^{-1}D_{g} \bar{\mathcal{L}}_a^j + \frac{1}{\Delta t} \left( \operatorname{Ad}_{\tau(\Delta t\xi_a^{j-1})}^*\mu_a^{j-1} -\mu_a^j \right) + \frac{1}{\Delta s} 
\left(\operatorname{Ad}_{\tau(\Delta s \eta_{a-1}^{j})}^*\lambda_{a-1}^{j} - \lambda_{a}^j \right) \\
&  \hspace{0 cm} + \bar F_d^1(\triangle_a ^j, g _a ^j , \xi _a ^j , \eta _a ^j )+ \bar F_d^2(\triangle _a ^{j-1},g_a ^{j-1}, \xi _a ^{j-1} , \eta _a ^{j-1}) + \bar F_d^3(\triangle_{a-1}^j,g_{a-1}^j, \xi _{a-1} ^j , \eta _{a-1} ^j )=0.
\end{align*}
 
\paragraph{Discrete forced Noether theorem.} In order to obtain the discrete forced Noether theorem, we have to assume  the discrete forces are orthogonal to the $H$ action, see \eqref{orthogonal_forces}. In our case, this reads
\footnotesize{
\[
\left\langle \operatorname{Ad}^*_{(g _a ^j ) ^{-1} }    \bar F ^1_d  ( \triangle _a ^j ,g _a ^j , \xi _a ^j , \eta _a ^j )+\operatorname{Ad}^*_{(g _a ^{j+1} ) ^{-1} }    \bar F ^2_d  ( \triangle _a ^j ,g _a ^j , \xi _a ^j , \eta _a ^j ) +\operatorname{Ad}^*_{(g _{a+1} ^j ) ^{-1} }    \bar F ^1_d  ( \triangle _a ^j ,g _a ^j , \xi _a ^j , \eta _a ^j ) , \zeta \right\rangle =0,
\]}
\normalsize
for all $ \zeta \in \mathfrak{h}  $.

In this particular case, the discrete forced Noether Theorem \ref{Cov_Noether_force} reads as follows.
 
\begin{theorem}  Suppose that the discrete Lagrangian $\bar{\mathcal{L}}_d: X_d^\triangle \times G \times \mathfrak{g}\times \mathfrak{g} \rightarrow \mathbb{R}$ is invariant under the left action of the Lie group $H$ on $G$ and suppose that the discrete forces $\bar{F}_d^k : X_d^\triangle \times G \times \mathfrak{g}\times \mathfrak{g} \rightarrow \mathfrak{g}  ^\ast $, $k= 1,2,3$ are orthogonal to this action. Suppose that $ g _d   $ is a solution of the DCEL equations for $ \bar{\mathcal{L} }_d $ with forces. Then, for all $0\leq B<C\leq A-1$, $0\leq K<L\leq N-1$, we have the conservation law
\begin{equation}\label{DCN_G_F}  
(\bar{\mathscr{J}}^F)_{B,C}^{K,L}(g _d  ) = 0,
\end{equation}
where $(\bar{\mathscr{J}}^F)_{B,C}^{K,L}$ is given by \eqref{Def_DN}, with $J ^k_{\mathcal{L}_d}( j ^1 \varphi _d (\triangle _a ^K )$ replaced by $J ^{k,F}_{\bar{\mathcal{L}}_d}( j ^1 g_d (\triangle _a ^K )$, where $J^{k,F} _{ \bar{\mathcal{L}} _d }:X_d^\triangle \times G \times \mathfrak{g} \times \mathfrak{g} \rightarrow \mathfrak{h}^*$, $k=1,2,3$, are defined by
\begin{equation*} 
\begin{aligned}
J^{1,F}_{\bar{\mathcal{L}}_d}(\triangle_a^j, g _a ^j , \xi _a ^j , \eta _a ^j )&= i ^\ast \mathrm{Ad}_{(g_a^j)^{-1}}^* \left( (g_a^j)^{-1}D_g \bar{\mathcal{L}}_a^j -  \frac{1}{\Delta t} \mu_a^j -  \frac{1}{\Delta s} \lambda_a^j + \bar F ^1_d  ( \triangle _a ^j ,g _a ^j , \xi _a ^j , \eta _a ^j ) \right), \\
J^{2,F}_{\bar{\mathcal{L}}_d}(\triangle_a^j, g _a ^j , \xi _a ^j , \eta _a ^j )&= i ^\ast\frac{1}{\Delta t}  \operatorname{Ad}_{(g_a^j)^{-1}}^* \left( \mu_a^{j}+  \bar F ^2_d  ( \triangle _a ^j ,g _a ^j , \xi _a ^j , \eta _a ^j ) \right), \\
J^{3,F}_{\bar{\mathcal{L}}_d}(\triangle_a^j, g _a ^j , \xi _a ^j , \eta _a ^j )&=  i ^\ast \frac{1}{\Delta s} \operatorname{Ad}_{(g_a^j)^{-1}}^* \left( \lambda_{a}^{j} +  \bar F ^3_d  ( \triangle _a ^j ,g _a ^j , \xi _a ^j , \eta _a ^j ) \right) . 
\end{aligned}
\end{equation*} 
\end{theorem}

\subsubsection{The $G$-invariant case and discrete covariant Euler-Poincar\'e equations}\label{G_invariant}

If the given Lagrangian density $ \mathcal{L} :J^1(X \times G) \rightarrow \mathbb{R}  $ is $G$-invariant, it induces the expression $\ell=\ell(\xi , \eta ): L(TX,\mathfrak{g}  )\rightarrow \mathbb{R}  $, as recalled in \S\ref{Cov_EL_LG}. The CEL equations for $ \mathcal{L}$ are equivalent to the covariant Euler-Poincar\'e equations for $\ell$.

In the case of a $G$-invariant Lagrangian, we shall choose a discrete Lagrangian $ \mathcal{L} _d $ that inherits the same invariance.
Consider the left action of $G$ on itself by left translation.  This action naturally lifts to $J^1 Y_d= X^\triangle_d \times G \times G \times G$. Then the discrete covariant Lagrangian $\mathcal{L} _d : J^1Y_d \rightarrow \mathbb{R}  $ is $G$-invariant if and only if its trivialized expression $\bar { \mathcal{L} }_d(\triangle _a ^j , g_a^j, \xi _a ^j ,\eta _a^j )$, defined through a local diffeomorphism $ \tau : \mathfrak{g}  \rightarrow G$,  does not depend on $ g _a ^j $. We thus obtain a discrete reduced Lagrangian 
\[
\ell_d(\triangle_a^j, \xi_a^j,\eta_a^j) : X_d ^\triangle\times \mathfrak{g} \times \mathfrak{g} \rightarrow \mathbb{R},
\]
that approximates the reduced Lagrangian $\ell$: 
\[
\ell_d(\triangle_a^j, \xi_a^j, \eta_a^j) \simeq \int_{\square_a ^j }\ell(t,s,\xi(t,s), \eta(t,s))ds \, dt.
\]

From the results obtained previously, it is straightforward to obtain the stationarity conditions associated to the covariant Euler-Poincar\'e principle. It suffices to set $D_g \bar{ \mathcal{L} }=0$ in \eqref{1}--\eqref{4}.

For example, if there are only temporal boundary conditions, we get the equations
\begin{equation}\label{1EP}
\frac{1}{\Delta t} \left(\mu_a^j-\operatorname{Ad}_{\tau(\Delta t\xi_a^{j-1})}^*\mu_a^{j-1}\right) + \frac{1}{\Delta s} 
\left(\lambda_{a}^j-\operatorname{Ad}_{\tau(\Delta s \eta_{a-1}^{j})}^*\lambda_{a-1}^{j}\right) = 0
\end{equation} 
for all $j=1,...,N-1$, $a=1,...,A-1$ with the natural (zero traction like) boundary conditions
\begin{equation}\label{2EP}
\frac{1}{\Delta t} \left( \mu_0^j-\operatorname{Ad}_{\tau(\Delta t\xi_0^{j-1})}^*\mu_0^{j-1} \right) + \frac{1}{\Delta s} \lambda_{0}^j  =0, \\
\quad\text{and}\quad \frac{1}{\Delta s} \operatorname{Ad}_{\tau(\Delta s \eta_{A-1}^{j})}^*\lambda_{A-1}^{j}  = 0,
\end{equation} 
for all $j=1,...,N-1$. Equations \eqref{1EP} are called the \textit{discrete covariant Euler-Poincar\'e equations}.

\subsection{Time and space discrete evolutions}\label{subsec_TSDE} 

\subsubsection{Symplectic properties of the time and space discrete evolutions}\label{symp_time_space}

As in \S\ref{symplectic_properties}, given the discrete Lagrangian $\mathcal{L} _d ( \triangle _a ^j , g _a ^j , g^{j+1}_a, g^j_{a+1})$, we can associate the discrete Lagrangians $\mathsf{L}_d=\mathsf{L}_d( \mathbf{g} ^j , \mathbf{g} ^{j+1}): G^{A+1} \times G^{A+1} \rightarrow \mathbb{R}$ and $\mathsf{N}_d=\mathsf{N}_d( \mathbf{g} _a, \mathbf{g} _{a+1}):G^{N+1} \times G^{N+1} \rightarrow \mathbb{R}$, where $ \mathbf{g} ^j =(g _0 ^j , ..., g _A ^j)$ and $ \mathbf{g} _a=(g^0 _a, ...,g ^N_a)$.
These Lagrangians are associated to the temporal and spatial discrete evolutions respectively.

As above, we shall fix a local diffeomorphism $ \tau : \mathfrak{g}  \rightarrow G $ in a neighborhood of the identity, such that $ \tau (0)=e$. To $\mathsf{L}_d$ and $\mathsf{N}_d$ are naturally associated  the discrete Lagrangian $\bar{ \mathsf{L}}_d: G^{A+1} \times \mathfrak{g}  ^{A+1} \rightarrow \mathbb{R}  $ and $\bar{ \mathsf{N}}_d: G^{N+1} \times \mathfrak{g}  ^{N+1} \rightarrow \mathbb{R}  $ defined by
\[
\bar{ \mathsf{L}}_d( \mathbf{g} ^j , \boldsymbol{\xi} ^j ):= \mathsf{L}_d( \mathbf{g} ^j , \mathbf{g} ^{j+1}) \quad\text{and}\quad \bar{ \mathsf{N}}_d( \mathbf{g} _a , \boldsymbol{\eta }_a ):= \mathsf{N}_d( \mathbf{g} _a, \mathbf{g} _{a+1}),
\]
where $\boldsymbol{\xi} ^j:= \frac{1}{\Delta t}\tau ^{-1} \left( ( \mathbf{g} ^j ) ^{-1}  \mathbf{g} ^{j+1} \right) \in \mathfrak{g}  ^{A+1}$ and $\boldsymbol{\eta } _a:= \frac{1}{\Delta s}\tau ^{-1} \left( ( \mathbf{g} _a ) ^{-1}  \mathbf{g} _{a+1} \right) \in \mathfrak{g}  ^{N+1}$.

\paragraph{Discrete time evolution.} The discrete Lagrangians $\bar{\mathsf{L}}_d$ and $\bar{\mathcal{L} }_d $ are related as
\begin{equation}\label{def_L_d}  
\bar{\mathsf{L}}_d(\mathbf{g}^j, \boldsymbol{\xi}^j) :=\sum_{a=0}^{A-1} \bar{\mathcal{L}}_d (\triangle_a^j, g_a^j, \xi_a^j, \eta_a^j).
\end{equation}
The variations of $ \boldsymbol{\xi} ^j $ are $\delta \boldsymbol{\xi}^j  = {\rm d}^R\tau^{-1} (\Delta t\boldsymbol{\xi}^j ) \left(-\boldsymbol{\zeta}  ^j  +
\operatorname{Ad}_{\tau(\Delta t \boldsymbol{\xi}^j )}\boldsymbol{\zeta}  ^{j+1} \right)/\Delta t$,
where $ \boldsymbol{\zeta} ^j = (\mathbf{g} ^j ) ^{-1} \delta \mathbf{g} ^j $.
Applying the discrete Hamilton principle to the discrete action
\begin{equation} \label{disc_action_time}
\mathfrak{S}_d(\mathbf{g}_d, \boldsymbol{\xi}_d) = \sum_{j=0}^{N-1}  \bar{\mathsf{L}}_d(\mathbf{g}^j, \boldsymbol{\xi}^j),
\end{equation}
without assuming any boundary conditions, we get the conditions
\begin{equation} \label{stat_dynamic_3}
\begin{aligned}
&\frac{1}{\Delta t} \left( \boldsymbol{\mu}^j- \operatorname{Ad}_{\tau(\Delta t \boldsymbol{\xi}^{j-1})}^* \boldsymbol{\mu}^{j-1} \right) =(\mathbf{ g}^j)^{-1}D_{\mathbf{ g}} \bar{\mathsf{L}}^j, \\
&\frac{1}{\Delta t} \boldsymbol{\mu}^0 =(\mathbf{ g}^0)^{-1}D_{\mathbf{ g}} \bar{\mathsf{L}}^0 , \qquad  \frac{1}{\Delta t}  \operatorname{Ad}_{\tau(\Delta t \boldsymbol{\xi}^{N-1})}^* \boldsymbol{\mu}^{N-1}  =0,
\end{aligned}
\end{equation}
where we defined
\[
\boldsymbol{\mu}^j:=\left( \operatorname{d}^R \tau^{-1}  (\Delta t \boldsymbol{\xi}^j)\right) ^\ast D_{\boldsymbol{\xi}} \bar{\mathsf{L}}^j.
\]
As in \S\ref{DTE} we can show, by using \eqref{def_L_d}, that the conditions \eqref{stat_dynamic_3} are equivalent to the DCEL equations \eqref{1} together with the boundary conditions \eqref{2}, \eqref{3}, \eqref{4}. The various boundary conditions treated in \S\ref{DTE} can be treated similarly here.

The discrete Cartan one-forms $\Theta^\pm_{\bar{\mathsf{L}}_d}$ on $G^{A+1}\times \mathfrak{g}^{A+1}$ are computed as
{\footnotesize\begin{align}
\Theta^-_{\bar{\mathsf{L}}_d}(\mathbf{g}^j, \boldsymbol{\xi}^j) & = \left\langle \frac{1}{\Delta t} \boldsymbol{\mu}^j- (\mathbf{ g}^j)^{-1}D_{\mathbf{ g}} \bar{\mathsf{L}}^j ,(\mathbf{g} ^j ) ^{-1}  d \mathbf{g} ^j\right\rangle \nonumber
\\
& =-  \sum_{a=1}^{A-1} \left\langle  (g_a^j)^{-1}D_g \bar{\mathcal{L}}_a^j -  \frac{1}{\Delta t} \mu_a^j -  \frac{1}{\Delta s}\left[ \lambda_a^j   - \operatorname{Ad}_{\tau(\Delta s \eta_{a-1}^{j})}^*\lambda_{a-1}^{j} \right],( g _a ^j ) ^{-1} d g _a ^j  \right\rangle \nonumber
\\
& \qquad - \left\langle  (g_0^j)^{-1}D_g \bar{\mathcal{L}}_0^j -  \frac{1}{\Delta t} \mu_0^j -  \frac{1}{\Delta s} \lambda_0^j, ( g _0 ^j ) ^{-1} d g _0 ^j \right\rangle \nonumber  \\
& \qquad  -  \left\langle  \frac{1}{\Delta s} \operatorname{Ad}_{\tau(\Delta s \eta_{A-1}^{j})}^*\lambda_{A-1}^{j}, ( g _A ^j ) ^{-1} d g _A ^j   \right\rangle   \nonumber
\\
& =- \sum_{a=0}^{A-1} \left( \theta_{\bar{\mathcal{L}}_d}^1( \triangle_a^j, g _a ^j , \xi _a ^j , \eta _a ^j ) + \theta_{\bar{\mathcal{L}}_d}^3( \triangle_a^j, g _a ^j , \xi _a ^j , \eta _a ^j ) \right) \label{L_d_min_one_form}
\\
\Theta^+_{\bar{\mathsf{L}}_d}(\mathbf{g}^j, \boldsymbol{\xi}^j) & =  \left\langle \frac{1}{\Delta t}  \operatorname{Ad}_{\tau(\Delta t \boldsymbol{\xi}^{j})}^* \boldsymbol{\mu}^{j} , (\mathbf{g} ^{j+ 1} ) ^{-1}  d \mathbf{g} ^{j+ 1} \right\rangle \nonumber \\
& =  \sum_{a=0}^{A-1}  \left\langle \frac{1}{\Delta t}  \operatorname{Ad}_{\tau(\Delta t\xi_a^{j})}^*\mu_a^{j} ,( g _a ^{j+ 1} ) ^{-1} d g _a ^{j+ 1}  \right\rangle 
= \sum_{a=0}^{A-1} \theta_{\bar{\mathcal{L}}_d}^2 ( \triangle_a^j, g _a ^j , \xi _a ^j , \eta _a ^j ),\label{L_d_plus_one_form}
\end{align}}
where we used \eqref{DCartan_1form} to write them in terms of the discrete Cartan form $ \theta _{\bar{ \mathcal{L} }_d }$ associated to $\bar {\mathcal{L}} _d $.

One can now proceed as in Propositions \ref{prop1}, \ref{prop1.5}, \ref{prop2} and obtain that the solutions $( g _a ^j , \xi _a ^j , \eta _a ^j )$ has the following property: the discrete flow map (discrete temporal evolution)
\[
( \mathbf{g} ^j , \boldsymbol{\xi} ^j ) \mapsto ( \mathbf{g} ^{j+1}, \boldsymbol{\xi} ^{j+1})
\]
is symplectic relative to the discrete Lagrangian symplectic form $ \Omega _{\bar{\mathsf{L}_d}}:= \mathbf{d} \Theta^-_{\bar{\mathsf{L}}_d}=\mathbf{d} \Theta^+_{\bar{\mathsf{L}}_d}$.

\paragraph{Discrete spatial evolution.} The discrete Lagrangians $\bar{\mathsf{N}}_d$ and $\bar{\mathcal{L} }_d $ are related as
\begin{equation}\label{def_N_d} 
\bar{\mathsf{N}}_d(\mathbf{g}_a, \boldsymbol{\eta}_a) :=\sum_{j=0}^{N-1} \bar{\mathcal{L}}_d (\triangle_a^j, g_a^j, \xi_a^j, \eta_a^j).
\end{equation} 
The variations of $ \boldsymbol{\eta} _a $ are $\delta \boldsymbol{\eta}_a  = {\rm d}^R\tau^{-1} (\Delta s \boldsymbol{\eta}_a )\left(-\boldsymbol{\zeta} _a +
\operatorname{Ad}_{\tau(\Delta s \boldsymbol{\eta}_a )}\boldsymbol{\zeta} _{a+1} \right)/\Delta s$, where $ \boldsymbol{\zeta} _a = ( \mathbf{g} _a ) ^{-1} \delta \mathbf{g} _{a}$. Applying the discrete Hamilton principle to the discrete action
\begin{equation} \label{disc_action_space}
\mathfrak{S}_d(\mathbf{g}_d, \boldsymbol{\eta}_d) = \sum_{a=0}^{A-1}  \bar{\mathsf{N}}_d(\mathbf{g}_a, \boldsymbol{\eta}_a),
\end{equation}
without assuming any boundary conditions, we get the conditions
\begin{equation} \label{stat_dynamic_4}
\begin{aligned}
&\frac{1}{\Delta s} \left( \boldsymbol{\lambda}_a-\operatorname{Ad}_{\tau(\Delta s \boldsymbol{\eta}_{a-1})}^* \boldsymbol{\lambda}_{a-1} \right) =(\mathbf{ g}_a)^{-1}D_{\mathbf{ g}} \bar{\mathsf{N}}_a,\\
&\frac{1}{\Delta s} \boldsymbol{\eta}_0=(\mathbf{ g}_0)^{-1}D_{\mathbf{ g}} \bar{\mathsf{N}}_0, \qquad  \frac{1}{\Delta s}  \operatorname{Ad}_{\tau(\Delta s \boldsymbol{\eta}_{A-1})}^* \boldsymbol{\lambda}_{A-1}  =0.
\end{aligned}
\end{equation}
where
\[
\boldsymbol{\lambda}_a:=\left( \operatorname{d}^R \tau^{-1}  (\Delta s \boldsymbol{\eta}_a)\right) ^\ast D_{\boldsymbol{\eta}} \bar{\mathsf{N}}_a.
\]
As in \S\ref{DTE} we can show, by using \eqref{def_N_d}, that the conditions \eqref{stat_dynamic_4} are equivalent to the DCEL equations \eqref{1} together with the boundary conditions \eqref{2}, \eqref{3}, \eqref{4}. The various boundary conditions treated in \S\ref{DTE} can be treated similarly here.

The discrete Cartan one-forms $\Theta^\pm_{\bar{\mathsf{N}}_d}$ on $G^{N+1}\times \mathfrak{g}^{N+1}$ are computed as
\begin{align}
\Theta^-_{\bar{\mathsf{N}}_d}(\mathbf{g}_a , \boldsymbol{\eta}_a ) & =  \left\langle \frac{1}{\Delta s} \boldsymbol{\lambda}_a-  (\mathbf{ g}_a)^{-1}D_{\mathbf{ g}}\bar{\mathsf{N}}_a ,(\mathbf{g}_a) ^{-1} d\mathbf{g}_a \right\rangle  \nonumber
\\
& =- \sum_{j=1}^{N-1} \left\langle  (g_a^j)^{-1}D_g \bar{\mathcal{L}}_a^j -  \frac{1}{\Delta s} \lambda_a^j  - \frac{1}{\Delta t} \left[  \mu_a^j  - \operatorname{Ad}_{\tau(\Delta t\xi_a^{j-1})}^*\mu_a^{j-1} \right], ( g _a ^j ) ^{-1} d g _a ^j \right\rangle  \nonumber
\\
& \qquad - \left\langle   (g_a^0)^{-1}D_g \bar{\mathcal{L}}_a^0 -  \frac{1}{\Delta s} \lambda_a^0 - \frac{1}{\Delta t} \mu_a^0 ,  ( g _a ^0 ) ^{-1} d g _a ^0 \right\rangle \nonumber \\
& \qquad -   \left\langle  \frac{1}{\Delta t} \operatorname{Ad}_{\tau(\Delta t\xi_a^{N-1})}^*\mu_a^{N-1},( g _a ^N ) ^{-1} d g _a ^N \right\rangle \nonumber
\\
& =- \sum_{a=0}^{A-1} \left( \Theta_{\bar{\mathcal{L}}_d}^1( \triangle_a^j, g _a ^j , \xi _a ^j , \eta _a ^j ) + \Theta_{\bar{\mathcal{L}}_d}^2( \triangle_a^j, g _a ^j , \xi _a ^j , \eta _a ^j ) \right) \label{N_d_min_one_form}
\\
\Theta^+_{\bar{\mathsf{N}}_d}(\mathbf{g}_a , \boldsymbol{\eta}_a ) & = \left\langle \frac{1}{\Delta s}  \operatorname{Ad}_{\tau(\Delta s \boldsymbol{\eta}_{a})}^* \boldsymbol{\lambda}_{a} , (\mathbf{g} _{a+ 1} ) ^{-1}  d \mathbf{g} _{a+ 1} \right\rangle \nonumber\\
& \qquad =  \sum_{j=0}^{N-1}  \left\langle \frac{1}{\Delta s} \operatorname{Ad}_{\tau(\Delta s \eta_{a}^{j})}^*\lambda_{a}^{j}, (g_{a+1} ^j ) ^{-1} dg_{a+1} ^j \right\rangle   =  \sum_{j=0}^{N-1} \Theta_{\bar{\mathcal{L}}_d}^3 ( \triangle_a^j, g _a ^j , \xi _a ^j , \eta _a ^j ).  \label{N_d_plus_one_form}
\end{align}
One can now proceed as in Propositions \ref{prop3}, \ref{prop4}, \ref{prop5} and obtain that the solutions $( g _a ^j , \xi _a ^j , \eta _a ^j )$ has the following property: the discrete flow map (discrete spatial evolution)
\[
( \mathbf{g} _a  , \boldsymbol{\eta }_a ) \mapsto ( \mathbf{g}_{a+1}, \boldsymbol{\eta} _{a+1})
\]
is symplectic relative to the discrete Lagrangian symplectic form $ \Omega _{\bar{\mathsf{N}}_d}:= \mathbf{d} \Theta^-_{\bar{\mathsf{N}}_d}=\mathbf{d} \Theta^+_{\bar{\mathsf{N}}_d}$.

\subsubsection{Discrete Lagrangian momentum maps of the time and space discrete evolutions} \label{L_d_N_d_momentum_map}

Consider as above the left subgroup action of $H$ on $G$. We assume that the discrete covariant density $ \bar{\mathcal{L} }_d$ is $H$-invariant. As in \S\ref{cov_VS_usual_momap}, the associated Lagrangians $\bar{ \mathsf{L}}_d $ and $\bar{ \mathsf{N}}_d $ are $H$-invariant under the diagonal actions. The expression of discrete Lagrangian momentum maps $\mathsf{J}^{\pm}_{\bar{\mathsf{L}}_d} : G^{A+1}\times \mathfrak{g}^{A+1} \rightarrow \mathfrak{h}^*$, $\mathsf{J}^{\pm}_{\bar{\mathsf{N}}_d} : G^{N+1}\times \mathfrak{g}^{N+1} \rightarrow \mathfrak{h}^*$ are easily obtained from the expressions of $\mathsf{J}^{\pm}_{\mathsf{L}_d} : G^{A+1}\times G^{A+1} \rightarrow \mathfrak{h}^*$, $\mathsf{J}^{\pm}_{\mathsf{N}_d} : G^{N+1}\times G^{N+1} \rightarrow \mathfrak{h}^*$, see \eqref{def_all_J}. They are related to the discrete covariant momentum maps $J^k_{\bar{\mathcal{L} } _d }$ via the following formulas
\begin{equation} \label{time_momentum_map}
\begin{aligned}
\mathsf{J}_{\bar{\mathsf{L}}_d}^+(\mathbf{g}^j, \boldsymbol{\xi}^j) & = \sum_{a=0}^{A-1}\left( J^2_{\bar{\mathcal{L}}_d}(\triangle_a^j, g _a ^j , \xi _a ^j , \eta _a ^j )  \right)\\
\mathsf{J}_{\bar{\mathsf{L}}_d}^-(\mathbf{g}^j, \boldsymbol{\xi}^j) & = -\sum_{a=0}^{A-1}\left( J^1_{\bar{\mathcal{L}}_d}(\triangle_a^j, g _a ^j , \xi _a ^j , \eta _a ^j )  + J^3_{\bar{\mathcal{L}}_d}(\triangle_a^j, g _a ^j , \xi _a ^j , \eta _a ^j ) \right), \\
\mathsf{J}_{\bar{\mathsf{N}}_d}^+(\mathbf{g}_a, \boldsymbol{\eta}_a) & = \sum_{j=0}^{N-1}\left( J^3_{\bar{\mathcal{L}}_d}(\triangle_a^j, g _a ^j , \xi _a ^j , \eta _a ^j )  \right),\\
\mathsf{J}_{\bar{\mathsf{N}}_d}^-(\mathbf{g}_a, \boldsymbol{\eta}_a) & = -\sum_{j=0}^{N-1}\left( J^1_{\bar{\mathcal{L}}_d}(\triangle_a^j, g _a ^j , \xi _a ^j , \eta _a ^j )  + J^2_{\bar{\mathcal{L}}_d}(\triangle_a^j, g _a ^j , \xi _a ^j , \eta _a ^j ) \right).
\end{aligned}
\end{equation}
The results obtained in \S\ref{cov_VS_usual_momap}, especially Lemma \ref{important_lemma} and Theorem \ref{summary_Noether} are of course still valid in the case of a Lie group and can be written in a trivialized form, using the quantities trivialized quantities $\bar{\mathscr{J}}_{B,C}^{K,L}(g _d  )$, 
$\mathsf{J}_{\bar{\mathsf{L}}_d}^\pm(\mathbf{g}^j, \boldsymbol{\xi}^j) $, $\mathsf{J}_{\bar{\mathsf{N}}_d}^\pm(\mathbf{g}_a, \boldsymbol{\eta}_a) $ together with the discrete covariant Noether Theorem in trivialized form (Theorem \ref{theorem_DCN_LA}). This exercise is left to the reader.

\section{Numerical example : the three dimensional geometrically exact beam model} \label{sect_num_ex}

In this section, we illustrate the results obtained in this paper with the example of a geometrically exact beam (\cite{Re1972, Si1985, SiMaKr1988}). We will take advantage of the multisymplectic character of the integrator to simulate the motion of the beam knowing the time evolution of position and the strain of one of the extremity. This unusual boundary value problem can be treated simply by our integrator for which the time and space discretization are discretized in the same way.
By switching the space and time variables, this boundary problem reduces to a standard boundary problem with given position and velocity at initial time.

In geometrically exact models, the instantaneous configuration of a beam is described by its line of centroids as a map $\mathbf{r}:[0,L]\rightarrow\mathbb{R}^3$ and the orientation of all its cross-sections at points $\mathbf{r}(s)$, where $s \in [0,L]$, by a moving orthonormal basis $\{\textbf{d}_1(s), \textbf{d}_2(s), \textbf{d}_3(s)\}$. The attitude of this moving basis is described by a map $\Lambda:[0,L]\rightarrow SO(3)$ satisfying $\mathbf{d}_I(s)=\Lambda(s)\mathbf{E}_I$,  $I=1,2,3$, where $\{\textbf{E}_1, \textbf{E}_2, \textbf{E}_3\}$ is a fixed orthonormal basis.
\begin{figure}[H]
\centering
\includegraphics[width=1.6 in]{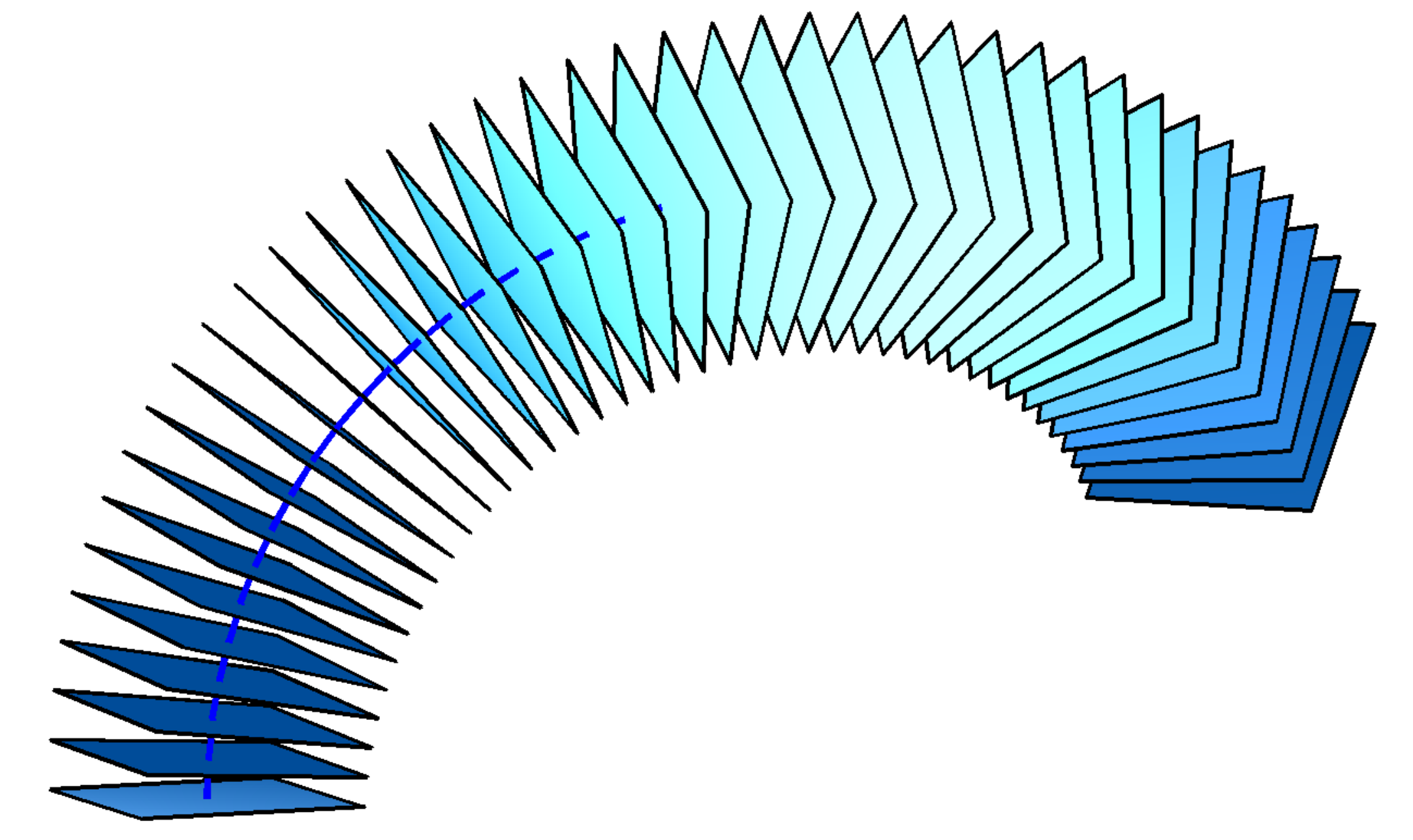}
\vspace{-3pt}
\caption{\footnotesize The geometrically exact model of beam is defined by the position $\mathbf{r}(t,s) \in \mathbb{R}^3$ of the line of centroids, and by the orientation $\Lambda(t,s) \in SO(3)$ of the cross sections. }
\end{figure}

\paragraph{Covariant formulation.} The convective covariant formulation of geometrically exact beams has been developed in \cite{ElGa-BaHoPuRa2010}, see especially \S6 and
\S7 of this paper. In this approach, the maps $ \Lambda , \mathbf{r} $ 
are interpreted as space-time dependent fields
\[
(t,s) \in \mathbb{R}  \times [0,L] \longmapsto g(t,s):=( \Lambda(t,s) , \mathbf{r}(t,s) )\in G=SE(3),
\]
taking values in the special Euclidean group. The fiber bundle of the problem is therefore given by
$X \times G \rightarrow X$, with $X= \mathbb{R}  \times [0,L] \ni (t,s)$ and the approach fits into the framework of \S\ref{Cov_EL_LG}. The convected variables $ \xi (t,s)= g (t,s)^{-1} \partial _t g(t,s)$ and $ \eta (t,s)= g (t,s)^{-1} \partial _s g(t,s)$ are here given by the convected angular and linear velocities and strains, i.e., 
\[
\xi = g ^{-1} \partial _t g= ( \Lambda ^{-1} \dot \Lambda , \Lambda ^{-1} \dot{\mathbf{r}} )=( \omega , \boldsymbol{\gamma}  ), \quad  \eta = g ^{-1} \partial _s g=  ( \Lambda ^{-1} \Lambda ', \Lambda ^{-1} \mathbf{r}')=( \Omega , \boldsymbol{\Gamma} ).
\]
The Lagrangian density of geometrically exact beams reads
\begin{equation}\label{Lagrangian}
\begin{aligned} 
\mathcal{L} (g, \xi , \eta )&= \frac{1}{2}\langle \mathbb{J}\xi, \xi\rangle - \frac{1}{2} \langle \mathbb{C}\ (\eta - \mathbf{E}_6), 
 (\eta-\mathbf{E}_6) \rangle- \Pi(g)= :K( \xi )- \Phi ( \eta )- \Pi(g),
\end{aligned}
\end{equation} 
where $K( \xi )$, $ \Phi ( \eta )$, and $ \Pi (g)$ are, respectively, the kinetic energy density, the strain energy density, and the external (such as gravitational) potential energy density.
Here $\mathbb{J}$ is a $6 \times 6$ diagonal whose diagonal elements are composed of the principal moments of inertia and the mass of the cross-section; the linear strain tensor $\mathbb C$ is a $6 \times 6$ diagonal matrix, whose diagonal elements depend on the cross-sectional area, the principal moments of inertia of the cross-sections, the Young's modulus, and the Poisson's ratio; and $ \mathbb{E}  _6 =(0,0,0,0,0,1)\in \mathbb{R}  ^6 $.
Both $ \mathbb{J}  $ and $ \mathbb{C}  $ are assumed to be independent of $(t,s)$.

The covariant Hamilton principle reads
\begin{equation}\label{trivialized_HP}
 \delta \int_0^T \!\!\!\int_0^L \left( K(\xi) - \Phi(\eta)  - \Pi (g) \right) \operatorname{d}\!s \operatorname{d}\!t =0,
\end{equation}
for arbitrary variations $ \delta g$ of $g$, vanishing at the boundary.
It yields the trivialized CEL equations,
\begin{equation}\label{Triv_CEL_BEam} 
\frac{\partial }{\partial t}\frac{\partial K}{\partial \xi }  - \mathrm{ad}_\xi^* \frac{\partial K}{\partial \xi }= \frac{\partial }{\partial s}\frac{\partial \Phi }{\partial \eta}- \mathrm{ad}_\eta^*\frac{\partial \Phi }{\partial \eta}   - g^{-1} \frac{\partial \Pi}{\partial g},
\end{equation}
see \eqref{Triv_CEL}. We refer to 
\cite{ElGa-BaHoPuRa2010} for a detailed derivation of these equations for geometrically exact models.

\paragraph{Time and space evolutions.} Following the theory developed in \S\ref{symplectic_properties}, we now define the Lagrangians associated to the temporal and spatial evolutions. We first do this at the continuous level.
The Lagrangian associated to the time evolution is $\bar{\mathsf{L}}(g, \xi )=\int_0^L \left( K(\xi) - \Phi(g ^{-1} \partial _s g) - \Pi(g)\right) \operatorname{d}\!s$ and the associated energy is
\begin{equation}\label{energy}
\mathsf{E}_{ \mathsf{L}}(g, \xi )= \int_0^L \left\langle \frac{\partial K}{\partial \xi }, \xi \right\rangle \operatorname{d}\!s-  \mathsf{L}(g, \xi )= \int_0^L \left( K(\xi) + \Phi(g ^{-1} \partial _s g) + \Pi(g)\right) \operatorname{d}\!s.
\end{equation}
The Lagrangian associated to the spatial evolution is
\[
\mathsf{N}(g, \eta  )=\int_0^T \left( K(g ^{-1} \partial _t g) - \Phi(\eta ) - \Pi(g)\right) \operatorname{d}\!t.
\]
One can also associate to $\mathsf{N}$ an energy function $\mathsf{E}_{\mathsf{N}}$ defined by the same formula, namely,
\begin{equation} \label{"energy"}
\begin{aligned}
\mathsf{E}_{\mathsf{N}}(g, \eta  ) &=
-\int_0^T\left\langle \frac{\partial \Phi }{\partial \eta  }, \eta  \right\rangle \operatorname{d}\!t- \mathsf{N}(g, \eta  ) \\
& = \int_0^T \left( -K(g ^{-1} \partial _t g)  - \left\langle \mathbb{C}  ( \eta - \mathbf{E} _6), \mathbf{E} _6 \right\rangle -\Phi(\eta ) + \Pi(g)\right) \operatorname{d}\!t.
\end{aligned}
\end{equation}
This energy function does not correspond to the physical energy.

Of course $\mathsf{E}_{\mathsf{L}}$ is conserved along the solutions of the EL equations for $\mathsf{L}$ on $ \mathcal{F} ([0,L], SE(3))$ and $\mathsf{E}_{\mathsf{N}}$ is conserved along the solution of the EL equations for $\mathsf{N}$ on $ \mathcal{F} ([0,T], SE(3))$. One has to remember also that the EL equations for $\mathsf{L}$ and $\mathsf{N}$ both imply not only the CEL equations for $\mathcal{L}$ but also boundary conditions: zero traction boundary condition in the case of the EL equations associated to  $\mathsf{L}$ and zero momentum boundary condition in space for the EL equations associated to $\mathsf{N}$, given respectively by
\begin{equation}\label{Conditions_de_bord} 
\frac{\partial \mathcal{L} }{\partial g'}(t,0)= \frac{\partial \mathcal{L} }{\partial g'}(t,L)=0,\;\forall \,t\quad\text{and}\quad \frac{\partial \mathcal{L} }{\partial \dot g}(0,s)= \frac{\partial \mathcal{L} }{\partial \dot g}(T,s)=0,\;\forall \,s.
\end{equation}

\paragraph{Numerical tests.} We shall use the multisymplectic integrator on Lie groups obtained in \S\ref{DCELELG} from the discrete covariant variational principle. Let us consider a geometrically exact beam of length $0.8 \,m$, and with cross-section given by a square of side $a=0.01\,m$. We assume that there are no exterior forces and that $\Pi(g)=0$ so that $\mathcal{L} $ is $SE(3)$ invariant. 

We choose the spacetime $X=[0,T] \times [0,L]$, with time of simulation $T=2\,s$, and length $L=0.8 \,m$. The space and time steps are $\Delta s=0.02\,m$ and $\Delta t =0.04\,s$. 
The spacetime is discretized as in \S\ref{DMFT} as $X_d =\{(j,a)\in \mathbb{Z}\times \mathbb{Z} \ | \ j=0,...,N-1, \ a=0,...,A-1\}$, 
where $N-1$ and $A-1$ correspond to $T$ and $L$, respectively. Recall that for all $(j,a) \in X_d$, we consider the triangles $\triangle_a^j= ((j,a),(j+1,a),(j,a+1))$ that also involve the nodes $(N,a)$ for all $a=0,..., A-1$, and $(j,A)$ for all $j=0,...,N-1$.

The construction of the discrete Lagrangian density $\mathcal{L}_d( \Delta _a ^j , \xi _a ^j , \eta _a ^j )$ as well as the detailed derivation of the associated discrete scheme obtained from the formula \eqref{1} are described in \cite{DeGBKoRa2013}.

In this example we consider the space evolution through the multisymplectic variational integrator, followed by the time reconstruction.

\medskip

\noindent \textbf{1) \textit{Space evolution}}: The problem treated here corresponds to the following situation. We assume that at the initial time $t=0$ and at the final time $t=T$, the velocity of the beam is zero. This corresponds to zero momentum boundary conditions. However the configuration of the beam at $t=0$ and $t=T$ is unknown. We assume however that we know the evolution (for all $ t\in [0,T]$) of one of the extremity, say $s=0$, as well as the evolution if its strain (for all $ t\in [0,T]$). The approach described in this paper, that makes use of both the temporal and spatial evolutionary description at both the continuous and discrete level, is especially well designed to discretize this problem in a structure preserving way.

Note that we do not impose zero traction boundary conditions, given here by 
\begin{equation}
\label{cont_boun_cond}
\left\{
\begin{array}{l}
\vspace{0.2cm}\displaystyle\left. (\Gamma-\mathbf{E}_3)\right|_{s=0}=0\\
\vspace{0.2cm}\displaystyle\left. (\Gamma-\mathbf{E}_3)\right|_{s=L}=0\\
\Omega  (0)=\Omega  (L)=0.
\end{array}
\right.
\end{equation}
at the two extremities of the beam.

The initial conditions are given by 
the configurations $\mathbf{g}_0= (g_0^0,...,g_0^N)$ and the strain $\boldsymbol{\eta}_0= (\eta_0^0,...,\eta_0^N)$ at the extremity $s=s_0$. 
In this example we choose the following configuration and strain:
\[
g_0^0=(\mathrm{Id}, (0,0,0)),\quad g_0^{j+1} = g_0^j\,\mathrm{cay}(\Delta t \xi _0 ^j ),\quad \text{for all $j=0,...,N-1$},
\]
where $ \xi _0 ^j = (0, -0.85, 0 ; 0,-0.1, 0)$, for all $j=0,...,N-1$, and
\[
\eta_0^j= \frac{1}{\Delta s} \mathrm{cay}^{-1}\left((g_0^j)^{-1} g^j_1 \right),\quad \text{for all $j=0,...,N-1$},
\]
where $g^0_1  =(\mathrm{Id}, (0,0,\Delta s))$ and $g^{j+1}_1= g^j_1 \mathrm{cay}(\Delta t \xi ^j _1 )$, with $\xi ^j _1  = (0.06, -0.849, -0.04;-0.03, -0.1, 0)$, for all $j=0,...,N-1$.

\begin{figure}[H]
\centering
\begin{center}
\begin{tabular}{cc}
\includegraphics[width=1.7 in]{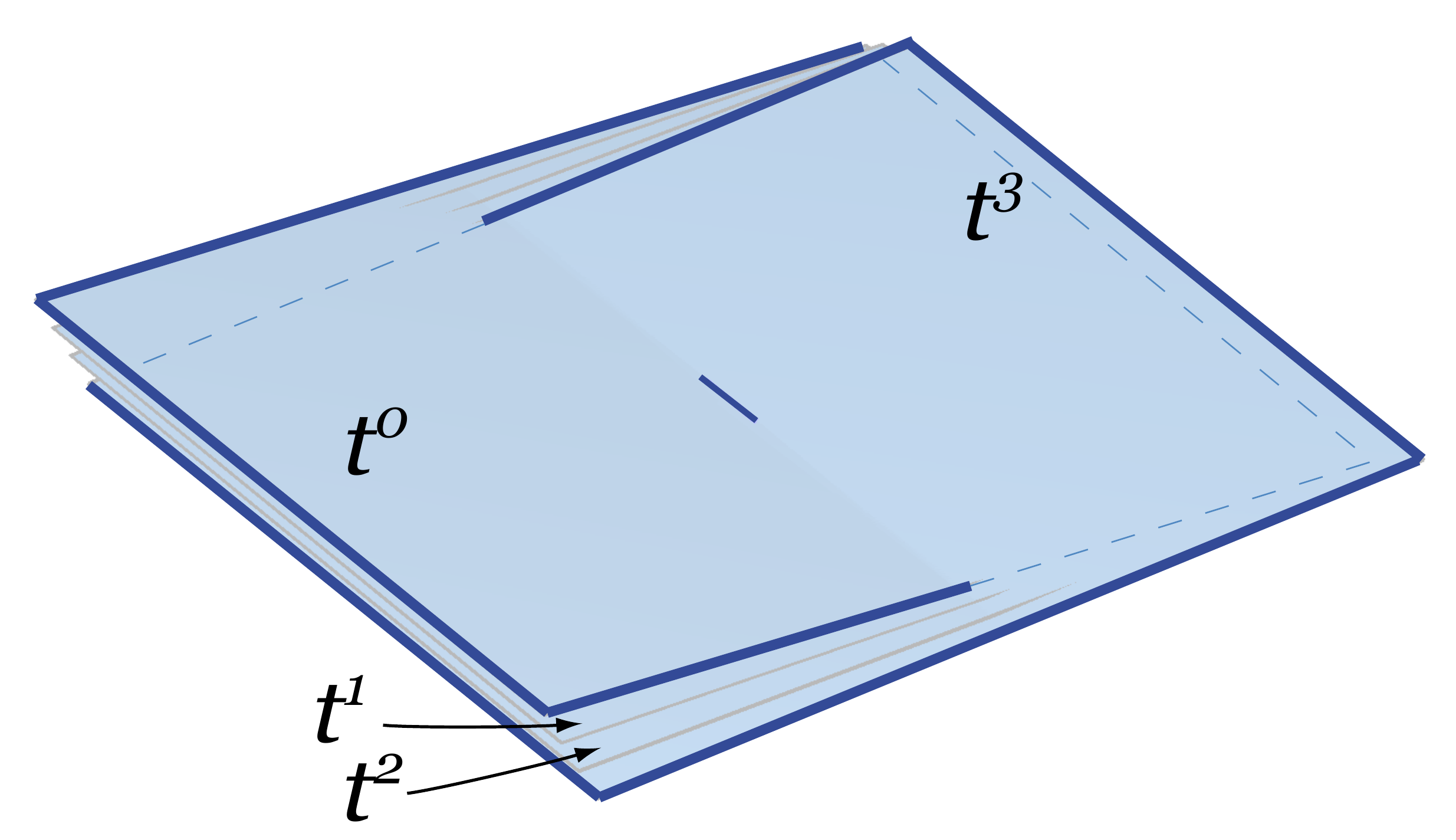}
\end{tabular}
\caption{\footnotesize Initial conditions $\mathbf{g}_0$ (enlarged), when $j\in\{0,1,2,3\}$.}\label{initial_condition}
\end{center}
\end{figure}

For the problem treated here, the boundary conditions thus given by \eqref{3}, which are the discretization of the right hand side condition in \eqref{Conditions_de_bord}, i.e., discrete zero momentum boundary conditions. Note also that since $ \Pi (g)=0$, $D_g \bar { \mathcal{L} } _a ^j =0$. So the discrete scheme is
\begin{align*}
&\begin{aligned}
&\frac{1}{\Delta t} \left(\mu_a^j-\operatorname{Ad}_{\tau(\Delta t\xi_a^{j-1})}^*\mu_a^{j-1}\right) + \frac{1}{\Delta s} 
\left(\lambda_{a}^j-\operatorname{Ad}_{\tau(\Delta s \eta_{a-1}^{j})}^*\lambda_{a-1}^{j}\right) = 0,\\
& \qquad \qquad \text{for all $j=1,...,N-2,\;\; a=1,...,A-2$},\qquad 
\end{aligned}  \\
& \begin{aligned}
&- \frac{1}{\Delta t} \mathrm{Ad}^*_{\tau(\Delta t \xi_a^{N-2})} \mu_a^{N-2} +   \frac{1}{\Delta s} \left(   \lambda_a^{N-1} -   \mathrm{Ad}^*_{\tau(\Delta s \eta_{a-1}^{N-1})} \lambda_{a-1}^{N-1}  \right)  =0, \\
&  \qquad \qquad  \text{for all $a =1,... A-2$,}
 \end{aligned} \\
&\begin{aligned}
&\frac{1}{\Delta t} \mu_a^0  + \frac{1}{\Delta s} 
\left(\lambda_{a}^0-\operatorname{Ad}_{\tau(\Delta s \eta_{a-1}^{0})}^*\lambda_{a-1}^{0}  \right) =0,\\
&\qquad  \qquad\text{and}\quad \frac{1}{\Delta t}  \operatorname{Ad}_{\tau(\Delta t\xi_a^{N-1})}^*\mu_a^{N-1} =0,\qquad \text{for all $a=1,...,A-2$},
\end{aligned} 
\end{align*}
where $\mu_a^j:=\left( {\rm d}\tau^{-1}_{\Delta t \xi^j_a }\right)^*\partial_\xi K(\xi_a^j)$ and $\lambda_a^j:= -\left({\rm d}\tau^{-1}_{\Delta s \eta ^j_a } \right)^*\partial_\eta \Phi(\eta_a^j)$.

\begin{figure}[H]
\centering
\begin{center}
\begin{tabular}{cc}
\includegraphics[width=1.3 in]{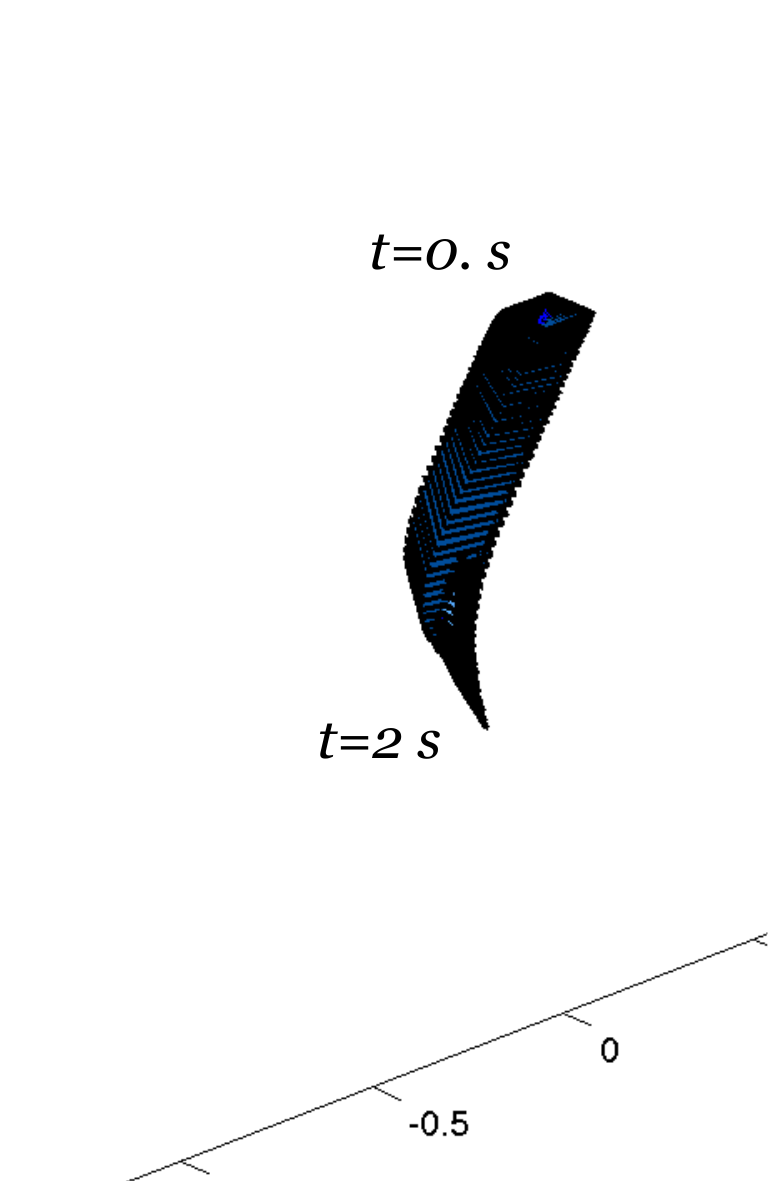}  \includegraphics[width=1.3 in]{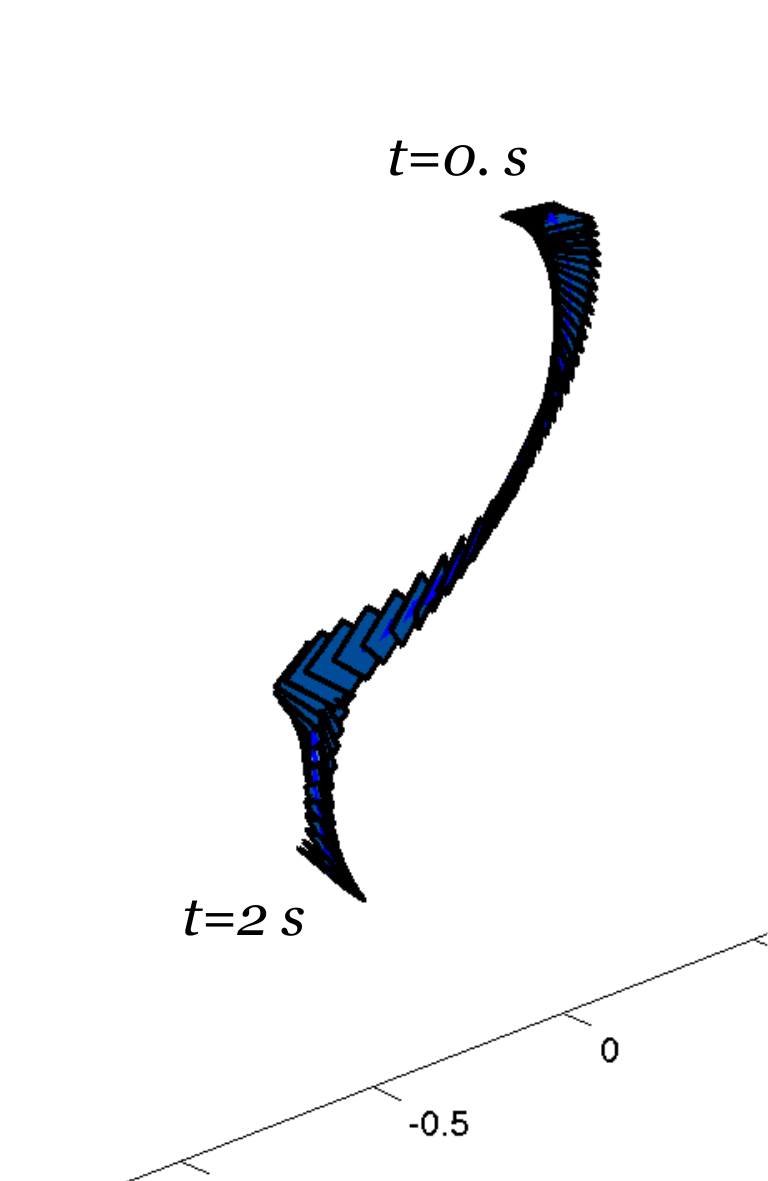}  
\includegraphics[width=1.3 in]{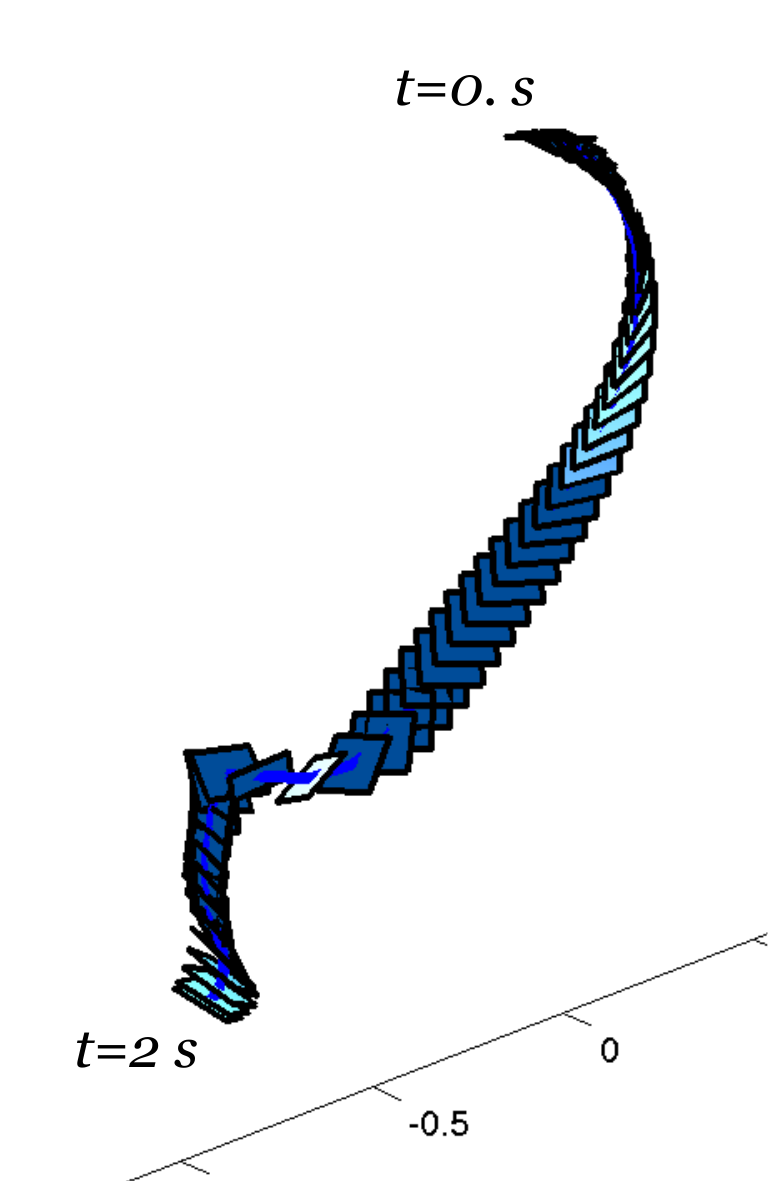}     \includegraphics[width=1.3 in]{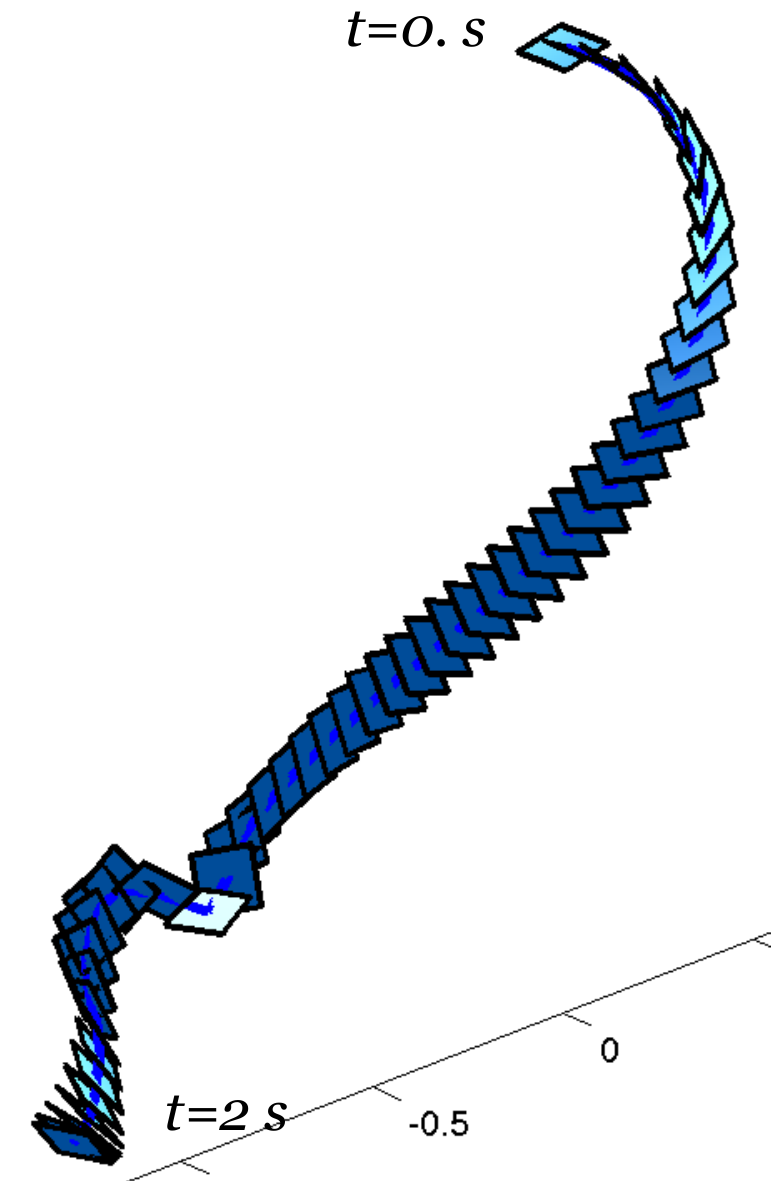}
\end{tabular}
\caption{\footnotesize Each figure represents the time evolution $\mathbf{g}  _a =\{g_a ^j, j=1,...,N\} $, of a given node $a$ of the beam.
The chosen nodes correspond to $s=0.26\,m$, $0.46\,m$, $0.62\,m$, $0.8\,m$.}\label{space_evolution}  
\end{center}
\end{figure}

This variational integrator produces the following displacement ``in space'' $\mathbf{g}_1,..., \mathbf{g}_A$ of the trajectories ``in time'' $\mathbf{g}_a =(g_a^0,..., g_a^N)$ of the beam sections (see Figure \ref{space_evolution}).
\color{black}

\begin{figure}[H]
\centering
\begin{center}
\begin{tabular}{cc}
\includegraphics[width=2.8 in]{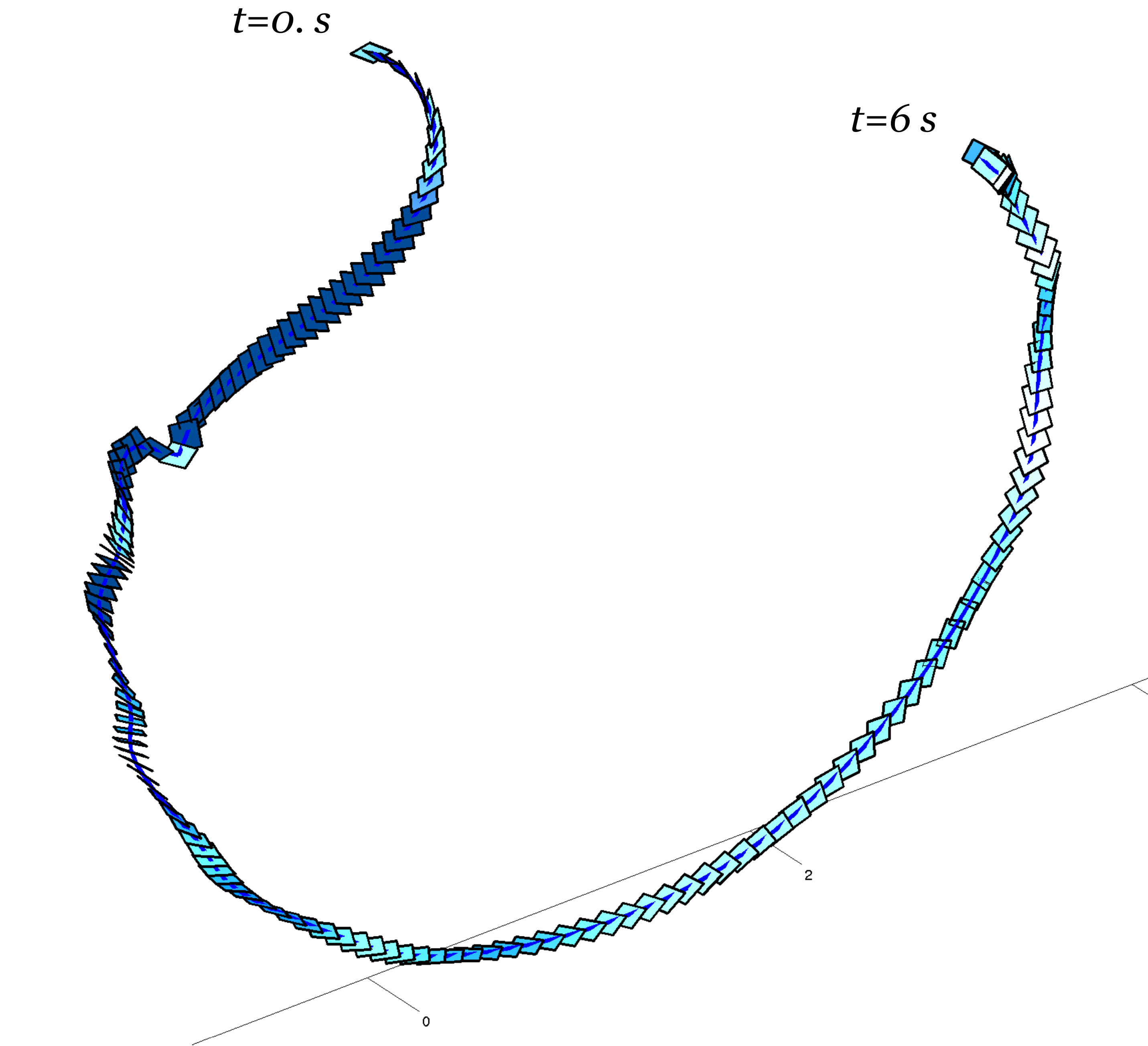} 
\end{tabular}
\caption{\footnotesize This figure represents the time evolution $\mathbf{g}  _a =\{g_a ^j, j=1,...,N\} $, of a given node $a$ of the beam, when $N=6s$.
The chosen nodes correspond to $s=0.8\,m$.}   
\end{center}
\end{figure}

\noindent\textit{Energy behavior.} The above DCEL equations together with the boundary conditions are equivalent to the DEL equations for $\bar{\mathsf{N}}_d( \mathbf{g} _a , \boldsymbol{\eta} _a )$ in \eqref{def_N_d}, see the discussion in \S\ref{cov_VS_usual_momap}. In particular, the solution of the discrete scheme define a discrete symplectic flow in space $( \mathbf{g} _a  , \boldsymbol{\eta }_a ) \mapsto ( \mathbf{g}_{a+1}, \boldsymbol{\eta} _{a+1})$ relative to the discrete symplectic form $ \Omega _{\bar{\mathsf{N}}_d}$ (see the end of \S\ref{symp_time_space}). As a consequence, the energy $\mathsf{E}_{\mathsf{N}_d}$ of to the Lagrangian $\mathsf{N}_d$ associated to the spatial evolution description is approximately conserved.

\medskip

\noindent\textit{Momentum map conservation.} Recall that the Lagrangian density $ \mathcal{L} $ is $SE(3)$ invariant, so the covariant Noether theorem is verified. At the discrete level, since the discrete Lagrangian density is also $SE(3)$-invariant (\cite{DeGBKoRa2013}), we get the discrete covariant Noether theorem
$\bar{\mathscr{J}}_{B,C}^{K,L}(g _d  ) = 0$, see \eqref{DCN_G}. Since the discrete Lagrangian $\mathsf{N}_d$ is $SE(3)$-invariant, the discrete momentum maps coincide: $\mathsf{J}^+_{\bar{\mathsf{N}}_d}=\mathsf{J}^-_{\bar{\mathsf{N}}_d}= \mathsf{J}_{\bar{\mathsf{N}}_d}$, and we have
\[
\mathsf{J} _{\bar{\mathsf{N}}_d}(\mathbf{g}_a, \boldsymbol{\eta}_a)
=\sum_{j=0}^{N-1} \Delta t \mathrm{Ad}^*_{(g_a^j)^{-1}} \lambda_a^j,
\]
see  \eqref{time_momentum_map} and \eqref{discr_momap_LA}. In view of the boundary conditions used here, from the discussion in \S\ref{cov_VS_usual_momap} it follows that the discrete momentum maps $\mathsf{J} _{\mathsf{N}_d}$ is conserved. This can be seen as a consequence of the covariant discrete Noether theorem $\bar{\mathscr{J}}_{B,C}^{0,N-1}( g _d )=0$.

\medskip

The discrete energy behavior and the conservation of the discrete momentum map $\mathsf{J}_{\bar{\mathsf{N}}_d}=(\mathsf{J}^1,...,\mathsf{J}^6) \in \mathbb{R}  ^6 $ are illustrated in Figure \ref{energy_momentum_space} below.
\begin{figure}[H]
\centering
\begin{center}
\begin{tabular}{cc}
 \includegraphics[width=1.6 in]{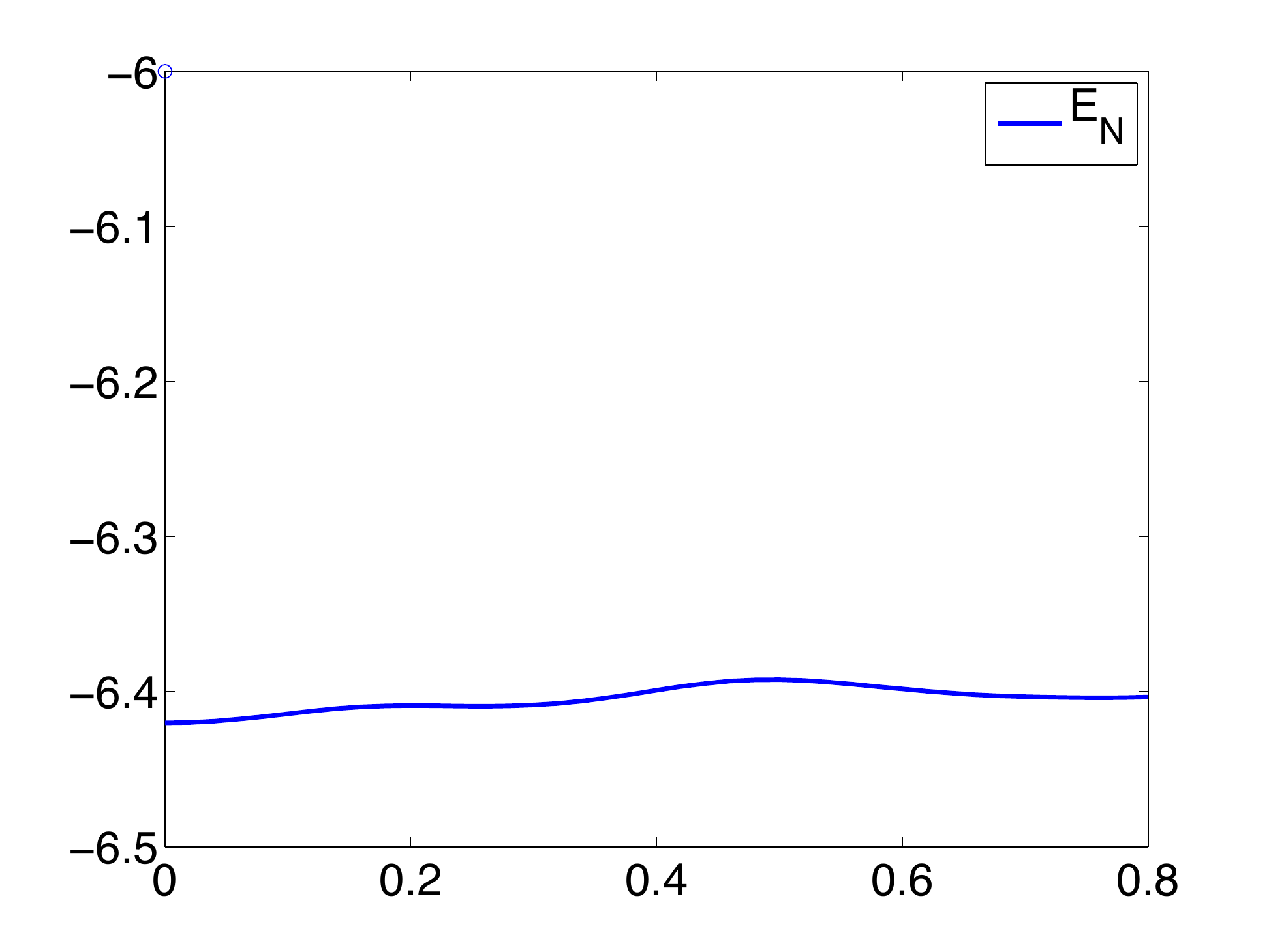}  \qquad   \includegraphics[width=1.6 in]{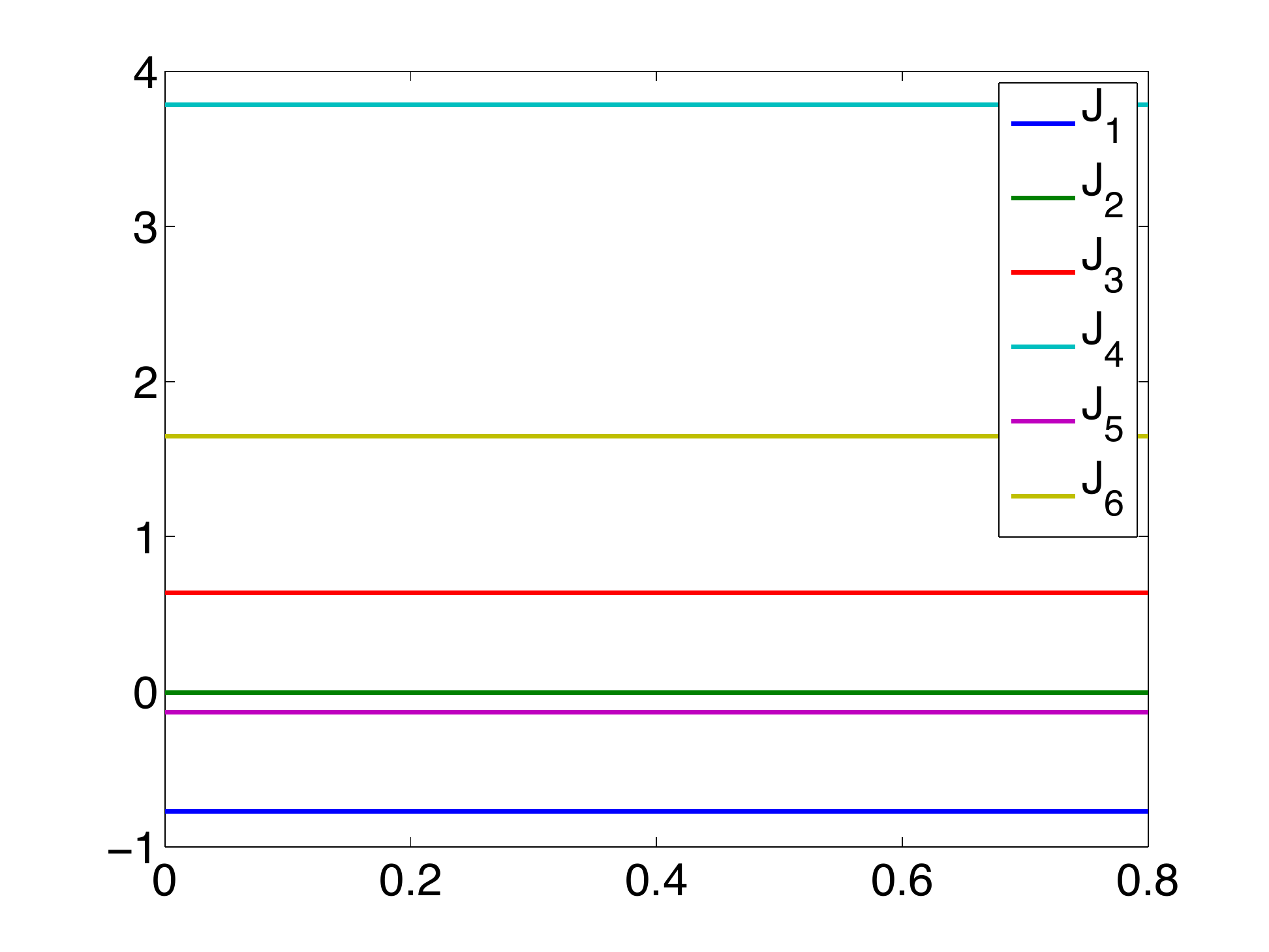}    
\end{tabular}
\vspace{-3pt}
\caption{\footnotesize The discrete energy function
$ \mathsf{E}_{\mathsf{N}_d}$ (left) is approximately conserved, and the discrete momentum maps 
$\mathsf{J}_{\mathsf{N}_d}$ (right) are exactly preserved. }\label{energy_momentum_space} 
\end{center}  
\end{figure}
The covariant discrete Noether theorem has also been numerically verified on the solutions of the discrete scheme.
We checked for example that $\bar{\mathscr{J}}_{0,A-1}^{0,N-1}( g _d )=0$. Recall that this follows from Lemma \ref{important_lemma} and the discussion after it. Indeed, we can write
\[
\begin{aligned}
\bar{\mathscr{J}}_{0,A-1}^{0,N-1}( g _d )=&\sum_{a=1}^{A-1} \left( J ^1_{\bar{\mathcal{L}}_d}( j ^1 g _d (\triangle _a ^0 )+J ^2_{\bar{\mathcal{L}}_d}(  j ^1 g _d (\triangle _a ^{N-1}) )+ J ^3_{\bar{\mathcal{L}}_d}( j ^1 g _d (\triangle _{a-1} ^0 )) \right)\\
&+ \mathsf{J}_{\bar{\mathsf{N}}_d}( \mathbf{g} _{A-1} ,\boldsymbol{\eta} _{A-1})- \mathsf{J}_{\bar{\mathsf{N}}_d}( \mathbf{g} _0 ,\boldsymbol{\eta} _{0}).
\end{aligned}
\]
The first line vanishes because of the boundary condition and the second line vanishes from the discrete Noether theorem.

\medskip

One can also consider the discrete Lagrangian $\bar{\mathsf{L}}_d$. However, as explained in \S\ref{DTE} (see \eqref{part1}) , the DEL equations for $\bar{\mathsf{L}}_d$ yield, in addition to the DCEL, zero traction boundary conditions, that are not verified here. As we have see, one can include boundary conditions in space by restricting $\bar{\mathsf{L}}_d$ to a subspace determined by these conditions.  However, as explained in \S\ref{DTE}, (B), these conditions have to be time independent, which is not the case in the problem considered here. So the equations of motion cannot be written as DEL equations for $\bar{\mathsf{L}}_d$ and the energy $\mathsf{E}_{\mathsf{L}_d}$ is not expected to be conserved. Of course, the same discussion holds at the continuous as well.
For similar reasons, the discrete momentum maps associated to $\mathsf{L}_d$ is not conserved, but verify the formula \eqref{non_cons_JL}.

\medskip 

\noindent \textbf{2) \textit{Reconstruction}}: The initial conditions are given by the set of configurations $\mathbf{g}_1,..., \mathbf{g}_A$ obtained through the space evolution (see Fig.\,\ref{space_evolution}). Thus we can immediately reconstruct the time advancement $\mathbf{g}^1,..., \mathbf{g}^N$ of the configuration of the beam, where $\mathbf{g}^j=(g_0^j,...,g_A^j)$.

\begin{figure}[H]
\centering
\begin{center}
\begin{tabular}{cc}
\includegraphics[width=0.9 in]{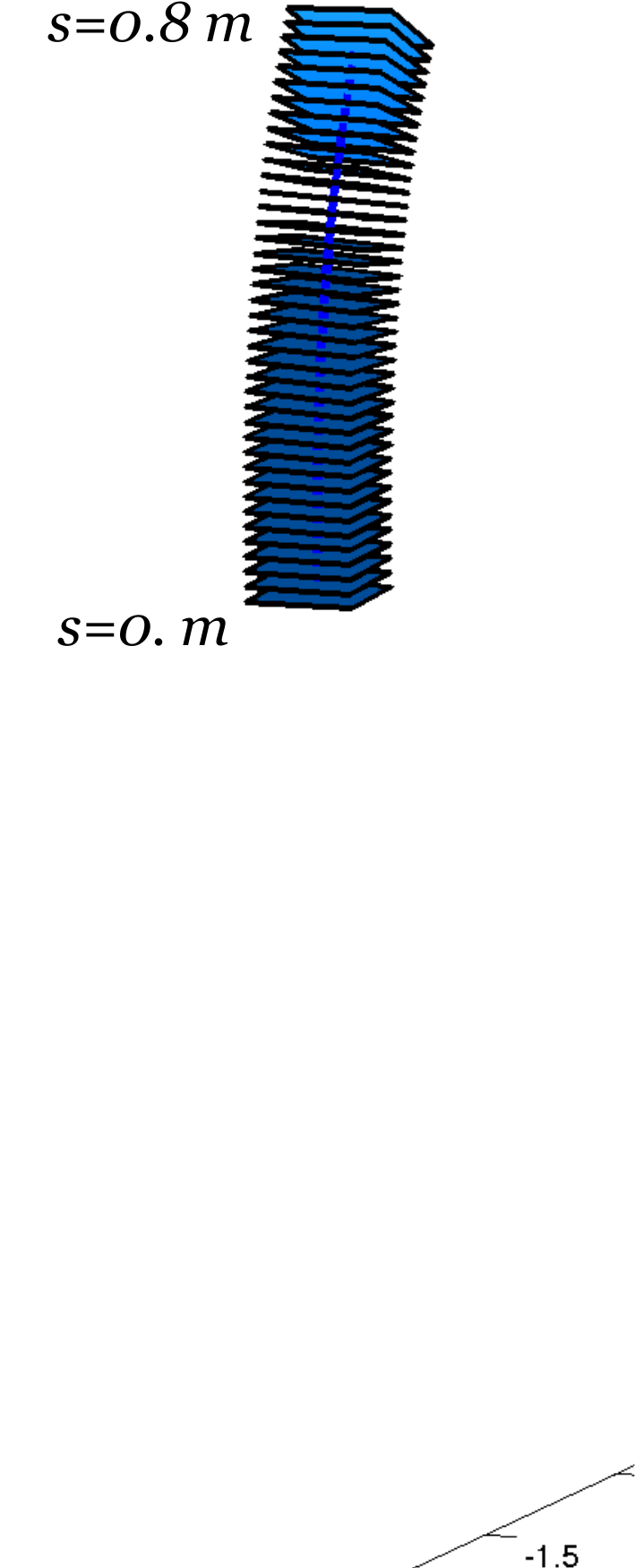}  \includegraphics[width=0.9 in]{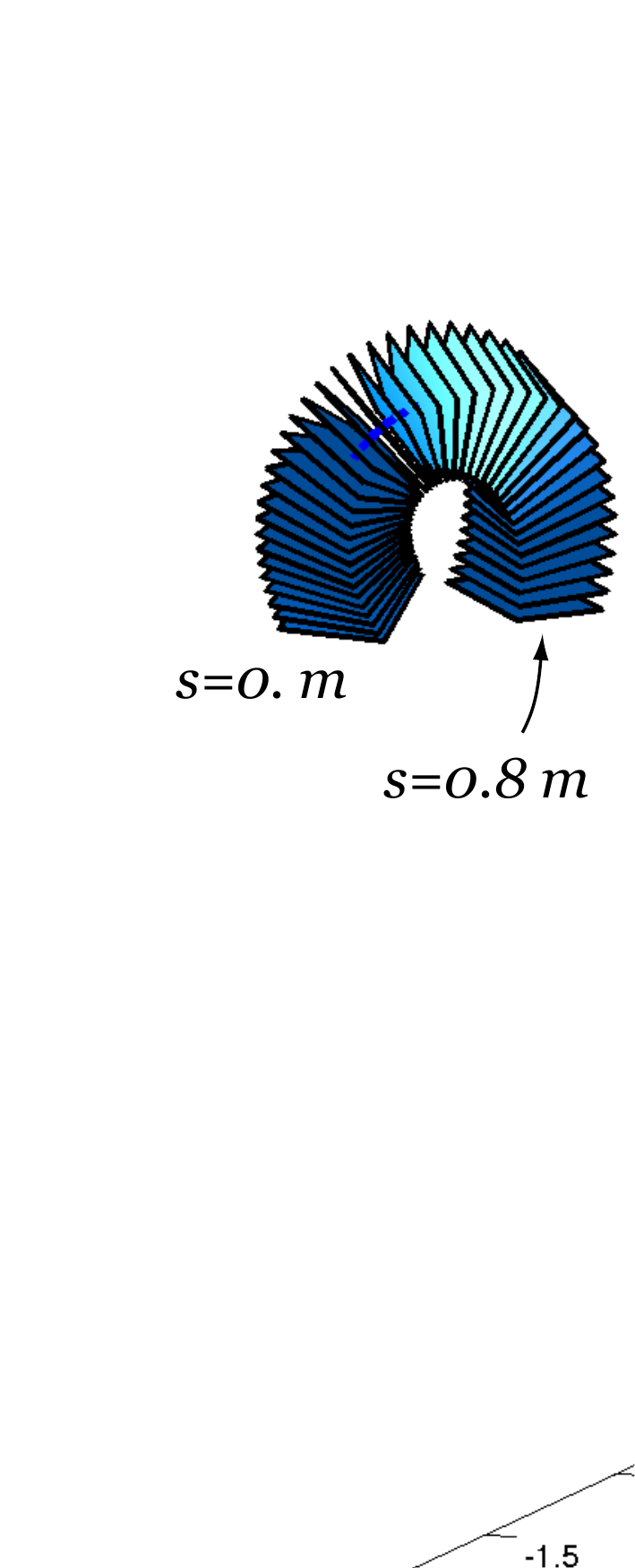}  
\includegraphics[width=0.9 in]{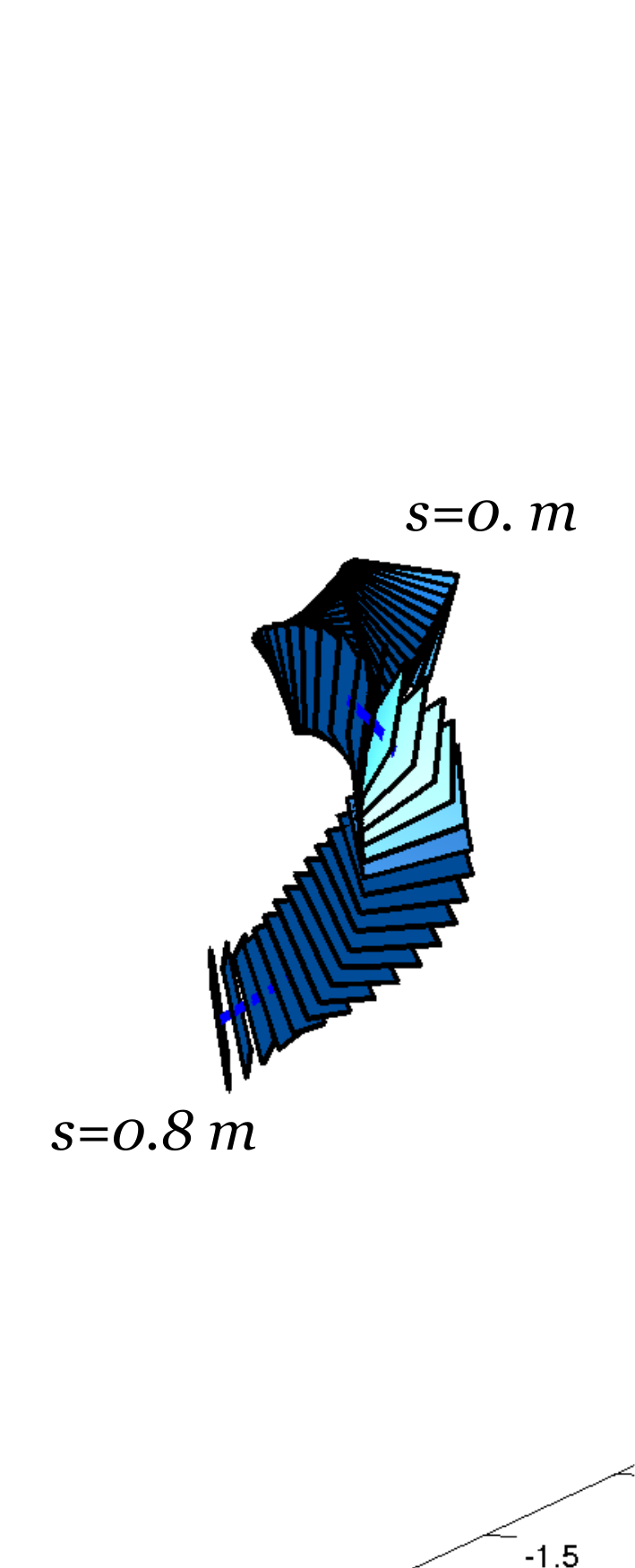}   
 \includegraphics[width=0.9 in]{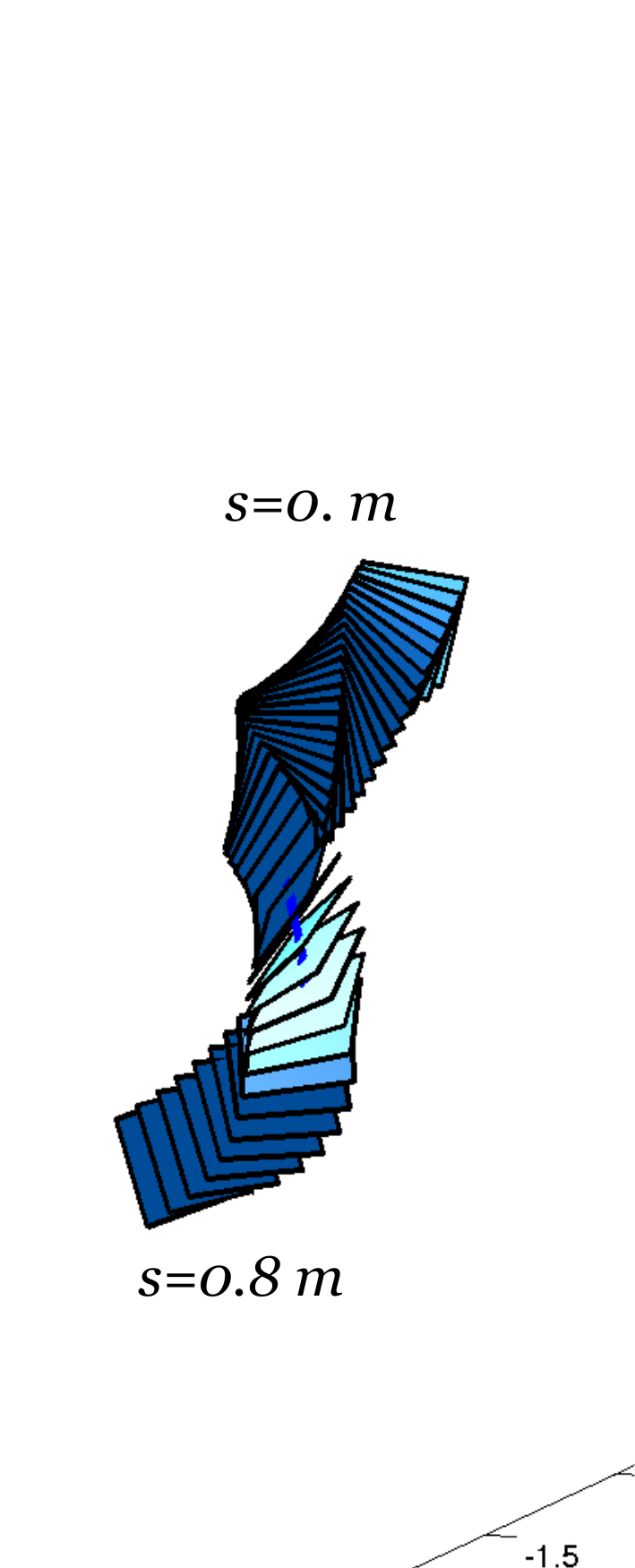}   \includegraphics[width=0.9 in]{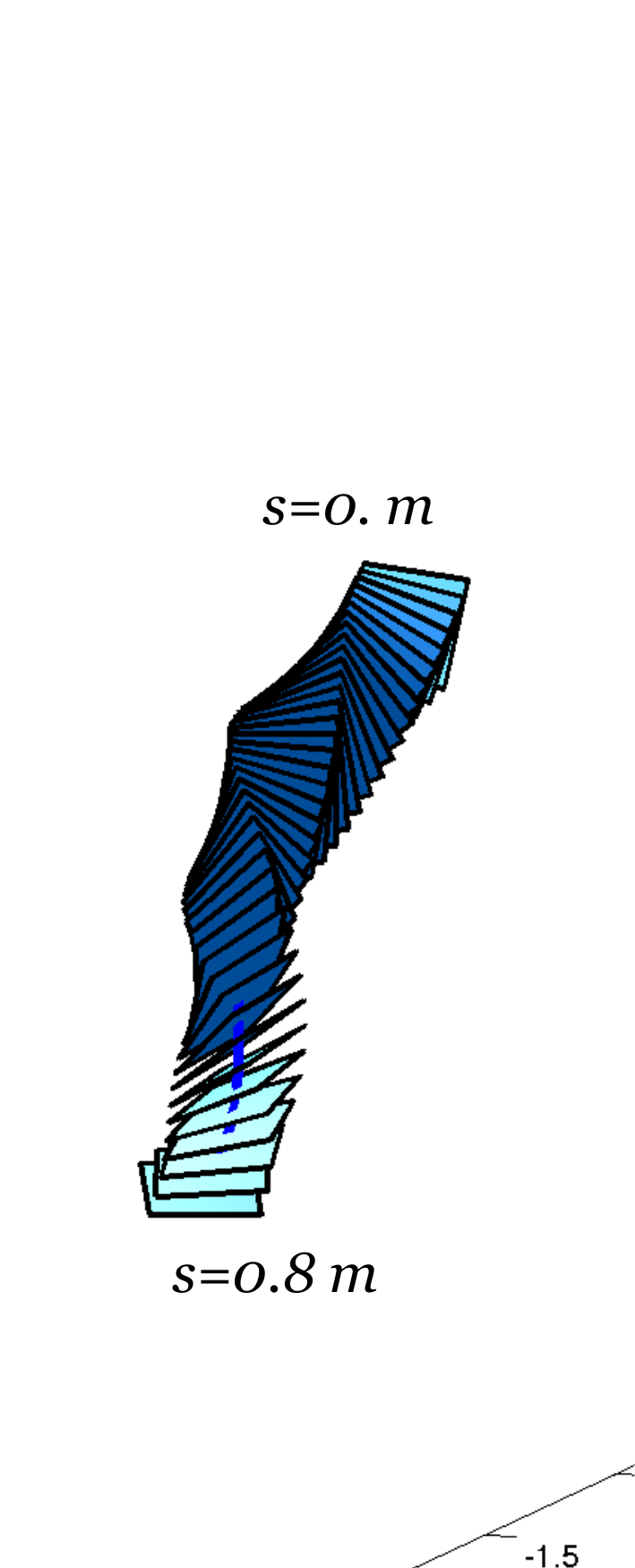} \includegraphics[width=0.9 in]{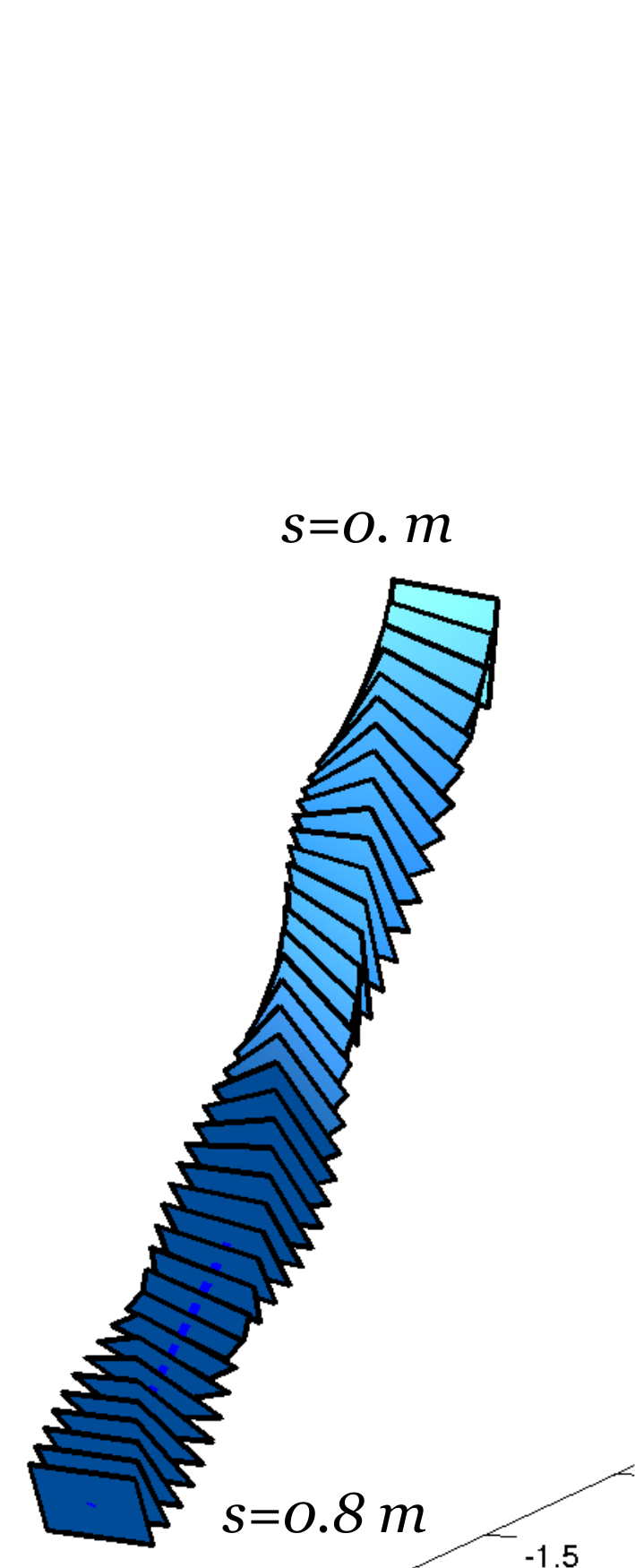}
\end{tabular}
\caption{\footnotesize  
From left to right: reconstruction of the trajectories in space of the sections, at times $t=0.1\,s, 0.6\,s, 1.25\,s, 1.4\,s, 1.45\,s, 1.95\,s$. }\label{reconstruction_space_time}  
\end{center}
\end{figure}

\medskip

\paragraph{Acknowledgment.}
We thank Marin Kobilarov and Julien Nembrini for their help with the numerical implementation.

{\footnotesize

\bibliographystyle{new}

}

\end{document}